\documentclass[10pt,a4paper]{article}
\usepackage{amsmath}
\usepackage{amsfonts}
\usepackage{amssymb}
\usepackage{graphicx}
\usepackage{enumerate}
\usepackage{color}
\newtheorem{notation}{Notation}[section]
\newtheorem{definition}{Definition}[section]
\newtheorem{observation}{Observation}[section]
\newtheorem{example}{Example}[section]
\newtheorem{lemma}{Lemma}[section]
\newtheorem{theorem}{Theorem}[section]
\newtheorem{proposition}{Proposition}[section]
\newtheorem{construction}{Construction}[section]
\newtheorem{splitproof}{Disjunction Sequence}[section]

\newcommand{\Z}{\mathbb Z}
\newcommand{\R}{\mathbb R}
\newcommand{\Q}{\mathbb Q}

\newcommand{\conv}{\textup{conv}}
\newcommand{\inter}{\textup{int}}

\newcommand{\cqfd}
{%
\mbox{}%
\nolinebreak%
\hfill%
$\square$
\medbreak%
\par%
}
\begin{document}
\title{Split Rank of Triangle and Quadrilateral Inequalities}
\author{Santanu Dey$^1$ \and Quentin Louveaux$^2$}
\footnotetext{${}^1$ Center for Operations Research and Econometrics, Universit\'e catholique
de Louvain, Belgium, email: santanu.dey@uclouvain.be\\This text presents research results of the Belgian Program on Interuniversity Poles of Attraction initiated by the Belgian State, Prime Minister's Office, Science Policy
Programming. The scientific responsibility is assumed by the authors.}
\footnotetext{${}^2$ Montefiore Institute, Universit\'e de Li\`ege, Belgium, email: q.louveaux@ulg.ac.be}
\maketitle
\begin{abstract}
A simple relaxation of two rows of a simplex tableau is a mixed integer set consisting of two equations with two free integer variables and non-negative continuous variables. Recently Andersen et al.~\cite{andersen:2007} and Cornu\'ejols and Margot~\cite{coma07} showed that the facet-defining inequalities of this set are either split cuts or intersection cuts obtained from lattice-free triangles and quadrilaterals. Through a result by Cook et al.~\cite{cook90}, it is known that one particular class of facet-defining triangle inequality does not have a finite split rank. In this paper, we show that all other facet-defining triangle and quadrilateral inequalities have a finite split-rank. The proof is constructive and given a facet-defining triangle or quadrilateral inequality we present an explicit sequence of split inequalities that can be used to generate it.
\end{abstract}

\section{Introduction}\label{intro}
Recently Andersen et al.~\cite{andersen:2007} and Cornu\'ejols and Margot~\cite{coma07} analyzed the facet-defining inequalities of the convex hull of the following mixed integer set:
\begin{eqnarray}\label{finiteprob1}
P(R,f) := \{(x,s)\in (\mathbb{Z}^2\times \mathbb{R}^k_{+})\,|\,f + Rs = x\},
\end{eqnarray}
where $f\in \mathbb{Q}^2\setminus \mathbb{Z}^2$ and $R = [r^1, r^2, ..., r^k] \in \mathbb{Q}^{2\times k}$. These inequalities are either split cuts or intersection cuts (Balas~\cite{balas:71}) (the so called \emph{triangle} and \emph{quadrilateral} inequalities).

The motivation for studying $P(R,f)$ is the following: Given two rows of a simplex tableau corresponding to integer basic variables that are at fractional values, $P(R,f)$ is obtained by relaxing the non-basic integer variables to be continuous variables and by relaxing the basic non-negative integer variables to be free integer variables. As $P(R,f)$ can be obtained as a relaxation of any mixed integer program, valid inequalities for the convex hull of $P(R,f)$ can be used as a source of cutting planes for general mixed integer programs. Empirical experiments with some classes of related cutting planes by Espinoza~\cite{espIPCO08} show that these new inequalities may be useful computationally. Various extensions to the basic relaxation $P(R,f)$ have also been recently studied where the inequalities are related to triangles and quadrilaterals; see for example Dey and Wolsey~\cite{deyIPCO:2008}, Andersen et al.~\cite{andersen:2009}, Dey and Wolsey~\cite{deywolseysfree}, Basu et al.~\cite{basuetallsfree}, Conforti et. al~\cite{conforticornuejolszambelli}, Fukasawa and G\"unl\"uk~\cite{fukasawagunluk09}.

The aim of this paper is to obtain a better understanding of the triangle and quadrilateral inequalities vis-$\grave{a}$-vis split inequalities. The motivation comes from the following well-known fact: One particular class of facet-defining triangle inequality for (\ref{finiteprob1}) does not have a finite split rank, i.e., it cannot be obtained by repeated application of split cuts (Cook et al.~\cite{cook90}). This leads to the following natural question: Which facet-defining inequalities for (\ref{finiteprob1}) have a finite split rank? \emph{We prove that the split rank of all the facet-defining inequalities of $\textup{conv}(P(R,f))$ is finite except for the particular class of triangle inequalities discussed in Cook et al.~\cite{cook90}.} For all facet-defining inequalities of the convex hull of (\ref{finiteprob1}) that have a finite split rank, we present an explicit sequence of split inequalities that can be used to generate them.

The paper is organized as follows. In Section \ref{presection}, we present some necessary definitions, the characterization of facet-defining inequalities for the convex hull of $P(R,f)$, and introduce the notation used in the rest of the paper. In Section \ref{mainresultsection}, we formally present the main result and provide an outline of its proof. The rest of the paper is devoted to the various steps of this proof. In particular, in Section \ref{fourone} and Section \ref{propsection} we present some general properties of split ranks that allow us to condense the analysis of inequalities for sets with at most four continuous variables. In Sections \ref{twovarsec}, \ref{threevarsec}, and \ref{fourvarsec}, we present split rank results for facet-defining inequalities of sets with two, three, and four continuous variables respectively.

\section{Preliminaries}\label{presection}
We assume that $P(R,f) \neq \emptyset$. If $R = [ r^1, ...,r^i, ..., r^k]$, then we say $r^i \in R$. We assume that if $r \in R$, then $r \neq (0,0)$. We begin this section with a definition of split rank. We then present a characterization of facet-defining inequalities for $\textup{conv}(P(R,f))$.

\subsection{Split Rank}
Consider a general mixed integer set $\mathcal{Q}:= \{(x,y) \in \mathbb{Z}^p \times \mathbb{R}^q\,|\, Gx + Hy \leq b \}$ where $G \in \mathbb{Q}^{m \times p}$, $H \in \mathbb{Q}^{m \times q}$, and $b \in \mathbb{Q}^{m \times 1}$. Let $\mathcal{Q}^{0} := \{(x,y) \in \mathbb{R}^p \times \mathbb{R}^q\,|\, Gx + Hy \leq b \}$ denote the linear programming relaxation of $\mathcal{Q}$. Given a vector $\pi \in \mathbb{Z}^p$ and $\pi_{0} \in \mathbb{Z}$, any vector $x\in \mathbb{Z}^p$ satisfies the \emph{split disjunction} defined as $(\pi^Tx \leq \pi_{0}) \,\vee \,(\pi^Tx \geq \pi_{0} + 1)$. An inequality that is valid for $\mathcal{Q}^{0}_{\pi,\pi_{0}}:= \textup{conv}( (\mathcal{Q}^0 \cap \{(x,y)| \pi^Tx \leq \pi_{0}\}) \cup (\mathcal{Q}^0 \cap \{(x,y)| \pi^Tx \geq \pi_{0} +1\}))$ is called a \emph{split inequality} (Balas~\cite{balassplitcut}).

The concept of split rank follows from the concept of \emph{split closure} of a mixed integer program introduced in Cook et al.~\cite{cook90}.
\begin{definition}[Split closure]
Given the linear programming relaxation $\mathcal{Q}^{0} := \{(x,y) \in \mathbb{R}^p \times \mathbb{R}^q\,|\, Gx + Hy \leq b \}$ of $\mathcal{Q} = \{(x,y) \in \mathbb{Z}^p \times \mathbb{R}^q\,|\,Gx + Hy \leq b \}$, the first split closure $\mathcal{Q}^1$ is defined as $\cap_{\pi \in \mathbb{Z}^p, \pi_{0} \in \mathbb{Z}}\mathcal{Q}^{0}_{\pi,\pi_{0}}$.
\end{definition}
The first split closure of a mixed integer set is a polyhedron (Cook et al.~\cite{cook90}) (see Andersen et al.~\cite{andcorli05}, Vielma~\cite{vielma07} and Dash et al.~\cite{dashGL:2007} for alternative proofs of this result.) Balas and Saxena~\cite{balassaxena08} and Dash et al.~\cite{dashGL:2007} conducted empirical studies of the strength of the first split closure. Cornu\'ejols and Li~\cite{cornuejols01} compare the closure with respect to 18 different classes of general purpose cuts. Recently Basu et al.~\cite{BasBonCorMar08} have made a comparison of the first split closure of $P(R,f)$ with the closure based on \emph{triangle} and \emph{quadrilateral} inequalities. Andersen et al.~\cite{AndWagWei2009} have generalized these results for sets with more rows.

The split closure procedure applied to the polyhedron $\mathcal{Q}^1$ gives the second split closure $\mathcal{Q}^2$. In general, we denote the $k^{\textrm{th}}$ split closure by $\mathcal{Q}^k$.

\begin{definition}[Split rank]\label{defnsplitrank} The split rank of an inequality $\alpha^Tx + \beta^Ty \leq \gamma$ wrt $\mathcal{Q}^0$ is defined as the smallest non-negative integer $k$ such that $\alpha^Tx + \beta^Ty \leq \gamma$ is a valid inequality for $\mathcal{Q}^{k}$.
\end{definition}

The split rank of a valid inequality for $\textup{conv}(\mathcal{Q})$ depends on the `formulation', i.e., the split rank of an inequality $\alpha^Tx + \beta^Ty \leq \gamma$ wrt $\mathcal{Q}^0$ may be different from the split rank wrt $\mathcal{Q'}^0$ where $\mathcal{Q} = \mathcal{Q}' := \{(x,y) \in \mathbb{Z}^p \times \mathbb{R}^q\,|\, G'x + H'y \leq b' \}$ but $\mathcal{Q}^0 \neq \mathcal{Q}'^0$ as $(G, H, b) \neq (G', H, b')$. If $\mathcal{Q}^0$ is clear from context, then we will typically not write the phrase `wrt to $\mathcal{Q}^0$'.

Upper bounds on split rank of inequalities is known to be finite in some cases. For example, Balas~\cite{balas1979}, Nemhauser and Wolsey~\cite{nemhauser90}, Balas et al.~\cite{balas93} show that the split rank of all valid inequalities is at most $n$ for a mixed binary program with $n$ binary variables. Dash and G\"unl\"uk~\cite{dashgunlukmixingrank} prove an upper bound of $n$ on the split rank of a mixing inequality based on $n$ rows.

\subsection{Facets of $\textup{conv}(P(R,f))$}\label{onfacetofprf}
We first begin with a discussion on valid inequalities of $\textup{conv}(P(R,f))$. A set $S \subseteq \mathbb{R}^2$ is called lattice-free if $\textup{interior}(S) \cap \mathbb{Z}^2 = \emptyset$. Lattice-free convex sets can be used to construct intersection cuts for $\textup{conv}(P(R,f))$ as described in the next proposition.
\begin{proposition}[Valid inequality from lattice-free convex set]\label{set2cut} Let $R \in \mathbb{Q}^{2 \times k}$ and $f \in \mathbb{Q}^2\setminus \mathbb{Z}^2$. Let $B$ be a closed lattice-free convex set containing $f$ in its interior. Let $\partial B$ represent the boundary of $B$. Define the vector $\phi(B) \in \mathbb{R}^k_{+}$ as
\begin{eqnarray}
\label{pidefn}\phi(B)_i = \left\{\begin{array}{cl} 0 & \textup{ if } r^i \in \textup{reccesion cone of } B\\
\lambda & \textup{ if } \lambda>0 \textup{ and }f + \frac{r^i}{\lambda} \in \partial B.\end{array}\right.
\end{eqnarray}
Then the inequality
\begin{eqnarray}\label{thcut}
\sum_{i = 1}^k\phi(B)_is_i \geq 1,
\end{eqnarray}
is a valid inequality for $\textup{conv}(P(R,f))$.
\end{proposition}
Note that the computation of the vector $\phi(B)$ depends on $B$, $f$, and $R$. However, we removed a reference to $f$ and $R$ in the notation `$\phi(B)$' for simplicity.

Valid inequalities that are not a conic combination of the inequalities $s_i \geq 0$ are called non-trivial inequalities. Every non-trivial valid inequality for $\textup{conv}(P(R,f))$ induces a lattice-free set as described next. (see Andersen et al.~\cite{andersen:2007}, Borozan and Cornu\'ejols~\cite{borozan:2007}, Cornu\'ejols and Margot~\cite{coma07}, Zambelli~\cite{zambelli:2008}).

\begin{proposition}[Lattice-free convex set from valid inequality]\label{inducedefn}
All non-trivial valid inequalities for $\textup{conv}(P(R,f))$ can be written in the form $\sum_{i = 1}^k\alpha_is_i \geq 1$ where $\alpha_i \geq 0$ $\forall 1 \leq i \leq k$. Then the set
\begin{eqnarray}\label{tildebcons}
L_{\alpha} = \textup{conv}\left(\cup_{\alpha_i>0}\left\{ f + \frac{r^i}{\alpha_i}\right\}\cup f\right) + \textup{cone}\left(\cup_{\alpha_i = 0}\{r^i\}\right)
\end{eqnarray}
is lattice-free and convex.
\end{proposition}

We call the set $L_{\alpha}$ as the \emph{induced lattice-free set}. The induced lattice-free set $L_{\alpha}$ depends on the coefficients $\alpha_i$, $f$, and on the columns $r^1, ..., r^k$. However, we removed a reference to $f$ and $R$ in the notation `$L_{\alpha}$' for simplicity.\\

Observe that when $\textup{cone}\{r^1, ..., r^{k}\} = \mathbb{R}^2$, $L_{\alpha} = \textup{conv}\left(\cup_{\alpha_i>0}\left\{ f + \frac{r^i}{\alpha_i}\right\}\right) + \textup{cone}\left(\cup_{\alpha_i = 0}\{r^i\}\right)$. Starting with a lattice-free set $B$ such that $f \in \textup{int}(B)$, it can be verified that
\begin{eqnarray}
\label{subseteq}L_{\phi(B)} \subseteq B.
\end{eqnarray}

We next present necessary conditions for an inequality to be facet-defining (see Andersen et al.~\cite{andersen:2007} for a proof). See Cornu\'ejols and Margot~\cite{coma07} for sufficient conditions for an inequality to be facet-defining.

\begin{figure}
\begin{center}
\includegraphics[width=0.75\linewidth]{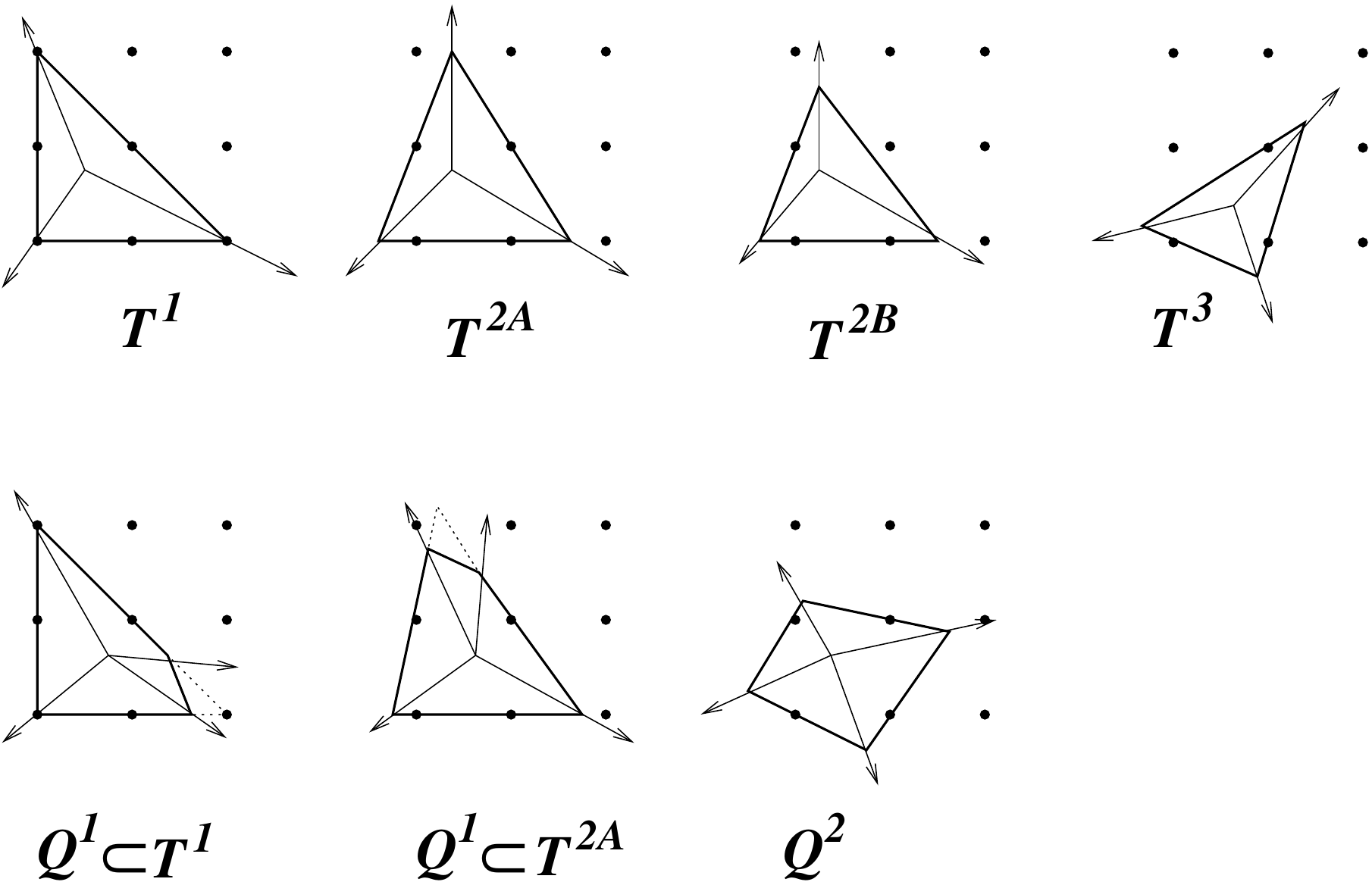}
\caption{Different cases of $L_\alpha$ (other than split sets) where $\alpha^Ts \geq 1$ is a facet-defining inequality for $\textup{conv}(P(R,f))$ and $\textup{cone}(R) = \mathbb{R}^2$.}
\label{allcases}
\end{center}
\end{figure}

\begin{proposition}\label{triadjust}
Let $\sum_{i = 1}^k\alpha_is_i \geq 1$ be a facet-defining inequality for $\textup{conv}(P([r^1, ..., r^k],f))$. If $\textup{cone}\{r^1, ..., r^k\} = \mathbb{R}^2$, then $f \in \textup{int}(L_{\alpha})$ and $L_{\alpha}$ is one of the following lattice-free sets:
\begin{enumerate}
\item Subset of Split Set: $\{(x_1, x_2)\,|\, \pi_0 \leq \pi_1 x_1 + \pi_2x_2 \leq \pi_0 + 1\}$ where $\pi_1, \pi_2, \pi_0 \in \mathbb{Z}$.
\item Type 1 triangle ($T^1$): Triangle with integral vertices and exactly one integer point in the relative interior of each side.
\item Type 2 triangle ($T^2$): Triangle with at least one non-integral vertex $v$ and the opposite side containing multiple integer points (not necessarily all in the relative interior). Let $S^1$ and $S^2$ be the two sides incident to $v$, and let $S^3$ be the third side. Then $T^2$ is further classified as:
\begin{enumerate}
\item $T^{2A}$: $S^1$ and $S^2$ contain one integer point in their relative interior.
\item $T^{2B}$: $S^1$ contains one integer point in its relative interior and $S^2$ does not contain any integer point in its relative interior. This triangle is a subset of some triangle of type $T^{2A}$.
\end{enumerate}
\item Type 3 triangle ($T^3$): Triangle with exactly three integer points on the boundary, one in the relative interior of each side and the vertices are non-integral.
\item Type 1 quadrilateral ($Q^1$): A subset of $T^{2A}$ or $T^1$ such that one side contains multiple integer points, two sides contain at least one integer point and the fourth side contains no integer point in its relative interior.
\item Type 2 quadrilateral ($Q^2$): A quadrilateral containing exactly one integer point in the relative interior of each of its sides and non-integral vertices.
\end{enumerate}
\end{proposition}
The various cases in Proposition \ref{triadjust} are illustrated in Figure \ref{allcases}.

\section{Main Result}\label{mainresultsection}
We prove the following result in this paper.
\begin{theorem}\label{mainresult}
Let $R= [ r^1, ..., r^k] \in \mathbb{Q}^{2 \times k}$ and $f \in \mathbb{Q}^2\setminus \mathbb{Z}^2$. Let $\sum_{i =1}^k\alpha_is_i \geq 1$ be a non-trivial facet-defining inequality for $\textup{conv}(P(R,f))$. The split rank of $\sum_{i =1}^k\alpha_is_i \geq 1$ is finite if and only if its induced lattice-free set $L_{\alpha}$ is not a triangle of type $T^1$.
\end{theorem}

The proof of Theorem \ref{mainresult} is technical and is presented in the rest of the paper. We next outline the various steps in the proof of Theorem \ref{mainresult}.\\\\
$\Rightarrow$ If $L_{\alpha}$ is a triangle of type $T^1$, then the inequality does not have a finite split rank. This follows from the proof in Cook et al.~\cite{cook90}. Also see Li and Richard~\cite{LiRic08}.\\\\
$\Leftarrow$ For the opposite direction, we need to show that the split rank of all facet-defining inequalities that are not split cuts (trivially) and whose induced lattice-free set is not a triangle of type $T^1$ is finite. Instead of considering only facet-defining inequalities, we analyze the split rank of the larger set of valid inequalities whose induced lattice-free set is described in Proposition \ref{triadjust} (for the case where $\textup{cone}_{r^i \in R}\{r^i\} = \mathbb{R}^2$) and the facet-defining inequalities where $\textup{cone}_{r^i \in R}\{r^i\} \neq \mathbb{R}^2$.
\begin{enumerate}
\item \textit{Restricting the proof to the case where $\textup{cone}_{r^i \in R}\{r^i\} = \mathbb{R}^2$ (Section \ref{fourone} - Section \ref{propsection})}: We can assume that the dimension of $\textup{cone}_{r \in R}\{r\}$ is 2, since otherwise the facet-defining inequalities for $\textup{conv}(P(R,f))$ are the split inequalities. We show in Proposition \ref{handlelowdim} that if $\sum_{i = 1}^k\alpha_is_i \geq 1$ is a facet-defining inequality for $\textup{conv}(P([r^1, ..., r^k],f))$ and $\textup{cone}\{r^1, r^2, ...., r^k\} \subsetneq \mathbb{R}^2$, then there exists a column $r^{k+1} \in \mathbb{R}^2$ and $\alpha_{k+1} \in \mathbb{R}_{+}$ such that $\textup{cone}\{r^1, ..., r^k, r^{k+1}\} = \mathbb{R}^2$, $\sum_{i = 1}^{k+1}\alpha_is_i \geq 1$ is facet-defining for $P([r^1, ..., r^k, r^{k+1}], f)$, and the induced lattice-free set of the inequality $\sum_{i = 1}^{k+1}\alpha_is_i \geq 1$ is not a triangle of type $T^1$. We show in Proposition \ref{generalprojection} that if $\sum_{i = 1}^{k +1}\alpha_i s_i \geq 1$ is a valid inequality of $P([r^1, ..., r^k, r^{k+1}],f)$ of split rank $\eta$ wrt $P([r^1, ..., r^k, r^{k +1}], f)^0$, then $\sum_{i = 1}^k\alpha_is_i \geq 1$ has a split rank at most $\eta$ wrt $P([r^1, ..., r^k], f)^0$.

    Thus it is sufficient to verify that the split rank of facet-defining inequalities for $\textup{conv}(P(R,f))$ is finite (except when induced lattice-free set is $T^1$) where $\textup{cone}_{r^i \in R}\{r^i\} = \mathbb{R}^2$.

\item \emph{Restricting the proof to sets with at most four continuous variables (Section \ref{propsection})}: We show in Lemma \ref{sameshape} that the split rank of an inequality $\sum_{i = 1}^{k_1}\alpha^1_i s_i \geq 1$ for $\textup{conv}(P(R^1,f^1))$ is lesser than (or equal to) the split rank of an inequality $\sum_{i = 1}^{k_2}\alpha^2_i s_i \geq 1$ for $\textup{conv}(P(R^2,f^2))$ if $L_{\alpha^1} \subseteq L_{\alpha^2}$. Specifically given $P(R,f)$ and the valid inequality for $\textup{conv}(P(R,f))$,
    \begin{eqnarray}
    \label{allkeq}\sum_{i = 1}^k \alpha_i s_i \geq 1,
    \end{eqnarray}
    let $A = \{i \in \{1, ..., k\}\,|\, f + \frac{r^i}{\alpha_i} \textup{ is a vertex of }L_{\alpha}\}$. Then consider the set
    \begin{eqnarray}
    \label{somek}x = f + \sum_{i \in A}r^is_i \quad s_i \geq 0, \quad x \in \mathbb{Z}^2
    \end{eqnarray}
    and the valid inequality
    \begin{eqnarray}
    \label{somekeq}\sum_{i \in A}\alpha_is_i \geq 1,
    \end{eqnarray}
    for (\ref{somek}). The split rank of (\ref{somekeq}) is equal to the split rank of (\ref{allkeq}) as the induced lattice-free sets for the inequalities (\ref{somekeq}) and (\ref{allkeq}) are identical. Since Proposition \ref{triadjust} shows that the induced lattice-free sets of all the facet-defining inequalities of $\conv(P(R,f))$ have at most four vertices ($|A| \leq 4$), it is sufficient to  show that the split rank of facet-defining inequalities for $\textup{conv}(P(R,f))$ is finite (except when their induced lattice-free set is $T^1$), where $\textup{cone}_{r^i \in R}\{r^i\} = \mathbb{R}^2$, $R \in \mathbb{Q}^{2 \times k}$ and $k \leq 4$.

Henceforth we call $P(R,f)$ as a $k$ variable problem when $R\in \Q^{2\times k}.$

\item \emph{Restricting the proof to `standard' triangles and quadrilaterals (Section \ref{propsection})}: We remark in  Observation \ref{transprop} that translating $f$ by an integral vector and multiplying $R$ and $f$ by a unimodular matrix $M$ does not change the split rank of a corresponding inequality. Thus the problem reduces to considering `standard' triangle and quadrilateral inequalities for problems with a maximum of four continuous variables.

\item \emph{Two variable problems (Section \ref{twovarsec})}: We prove in  Proposition \ref{2varthm} that the split rank of a facet-defining inequality for $\conv(P([r^1r^2],f))$ is at most 2. Note that $\textup{cone}\{r^1 , r^2\} \subsetneq \mathbb{R}^2$. However, this case is analyzed since this result is required for showing that the split rank of inequalities whose induced lattice-free set is a triangle (except $T^1$) is finite.
\item \emph{Three variable problems (Section \ref{threevarsec})}: The induced lattice-free set of a valid inequality not dominated by a split inequality, for the set $\conv(P([r^1r^2r^3],f))$ is a triangle when $\textup{cone}\{r^1, r^2, r^3\} = \mathbb{R}^2$. We first consider triangles of type $T^2$. As discussed in Proposition \ref{triadjust}, $T^2$ is subdivided into two classes: $T^{2A}$ and $T^{2B}$. We show in Proposition \ref{theoremt2c} that the split rank of an inequality whose induced lattice-free set is a triangle of type $T^{2B}$ is finite. This is the most technical part of the proof and is subdivided into four cases. The proof involves giving an explicit sequence of split disjunctions that yields the triangle inequality in a finite number of steps.

    It is then shown that the split rank of an inequality whose induced lattice-free set is either $T^{2A}$ or $T^3$ is at most one more than the split rank of a suitable constructed valid inequality whose induced lattice-free set is a triangle of type $T^{2B}$.
\item \emph{Four variable problems (Section \ref{fourvarsec})}: For the four variable case, if the induced lattice-free set is not a triangle, then it is a quadrilateral of type either $Q^1$ or $Q^2$. If the induced lattice-free set is a quadrilateral of type $Q^1$ and this quadrilateral is a subset of a triangle of type $T^{2A}$, then by Proposition \ref{sameshape} the split rank of the inequality is finite. We show in Proposition \ref{q1prop} that the split rank of inequalities whose induced lattice-free set is a quadrilateral of type $Q^1$ is finite, even when this quadrilateral is a proper subset of a triangle of type $T^1$. In this case the split rank is at most one more than the maximum of the split ranks of two suitably constructed inequalities with induced lattice-free set of type $T^2$. We show in Proposition \ref{q2prop} that the split rank of inequalities whose induced lattice-free set is a quadrilateral of type $Q^2$ is finite. This split rank is at most one more than the maximum of the split ranks of two suitably constructed inequalities with induced lattice-free set of type $Q^1$ or $T^{2}$.
\end{enumerate}

\section{Analyzing $L_{\alpha}$ when $\textup{cone}\{r^1, ..., r^k\} \subsetneq \mathbb{R}^2$ }\label{fourone}
Proposition \ref{triadjust} describes the shapes of $L_{\alpha}$ when $\textup{cone}\{r^1, ..., r^k\} = \mathbb{R}^2$. We now present a result to handle the case when $\textup{cone}\{r^1, ..., r^k\} \subsetneq \mathbb{R}^2$ for the proof of Theorem \ref{mainresult}. We need the following preliminary result proven in Andersen et al.~\cite{andersen:2007}.
\begin{lemma}
Let $\sum_{i=1}^k \alpha_i s_i \geq 1$ be a valid inequality for $\conv(P(R,f))$ such that $L_{\alpha}$ is not contained in any split set $\{(x_1,x_2)\in \R^2\mid \pi_0 \leq \pi^T x\leq \pi_0 + 1\}$ where $\pi\in \Z^2, \pi_0\in \Z$. Then $L_{\alpha}$ is bounded and $\alpha_i>0$ for all $i=1,\ldots, k.$
\end{lemma}
\begin{proposition}\label{handlelowdim}
Let $\sum_{i = 1}^k\alpha_is_i \geq 1$ be a facet-defining inequality for $\textup{conv}(P([r^1, ..., r^k],f))$ that is not dominated by any split inequality. If $\textup{dim}(\textup{cone}\{r^1, ..., r^k\}) = 2$ and $\textup{cone}\{r^1, ..., r^k\} \subsetneq \mathbb{R}^2$, then there exists a column $r^{k+1} \in \mathbb{R}^2$ and $\alpha_{k + 1}>0 $ such that
\begin{enumerate}[(i)]
\item $\textup{cone}\{r^1, ..., r^k, r^{k +1} \} = \mathbb{R}^2$,
\item $\sum_{i = 1}^{k +1}\alpha_is_i \geq 1$ is a facet-defining inequality for $\textup{conv}(P([r^1, ..., r^{k +1}], f))$,
\item $L':= \textup{conv}\left(\cup_{1 \leq i \leq k + 1}\left\{ f + \frac{r^i}{\alpha_i}\right\}\right)$, the induced lattice-free set of the inequality $\sum_{i = 1}^{k +1}\alpha_is_i \geq 1$ is not a triangle of type $T^1$.
\end{enumerate}
\end{proposition}
\textbf{Proof:} We present the proof for the case where $\textup{cone}\{r^1, ..., r^k\}$ is not a half-space. The other case can be similarly handled. Then WLOG let $\textup{cone}\{r^1, ..., r^k\} = \textup{cone}\{r^1, r^2\}$.

Since $L_{\alpha}$ is not contained in any split set, it is a bounded set. Moreover $\alpha_i > 0$ $\forall i$, and hence $L_{\alpha} =\textup{conv}\left(\cup_{1 \leq i \leq k}\left\{ f + \frac{r^i}{\alpha_i}\right\} \cup \left\{f\right\}\right)$. Choose any vector $r \in \mathbb{Q}^2$ such that $\textup{cone}\{r^1, r^2, r\} = \mathbb{R}^2$. Let
\begin{eqnarray}
\gamma := \textup{inf}\left\{\beta \in \mathbb{R}_{+}\,|\, \textup{conv}\left(\cup_{1 \leq i \leq k}\left\{ f + \frac{r^i}{\alpha_i}\right\}\cup \left\{f + \frac{r}{\beta}\right\}\right) \textup{ is lattice-free}\right\}.
\end{eqnarray}

\noindent Claim: $\gamma > 0$ and there exists $\tilde{\beta} \in \mathbb{R}_{+}$ such that $\gamma = \tilde{\beta}$. Assume by contradiction that $\gamma =0$. Then the set $S:=\textup{conv}\left(\cup_{1 \leq i \leq k}\left\{ f + \frac{r^i}{\alpha_i}\right\}  \right) + \textup{cone}(r)$ is lattice-free. Observe that $f\in S$. This would imply that $\textup{conv}\left(\cup_{1 \leq i \leq k}\left\{ f + \frac{r^i}{\alpha_i}\right\}\cup \left\{ f \right\}\right) + \textup{cone}(r, -r)$ is lattice-free (see Basu et al.~\cite{basuetalllatticefree}), contradicting the fact that the set $\textup{conv}\left(\cup_{1 \leq i \leq k}\left\{ f + \frac{r^i}{\alpha_i}\right\}\cup \left\{ f \right\}\right)$ is not contained in any split set. Therefore, $\gamma >0$. Now choose a suitably small $\hat{\beta} >0$ such that the set $\textup{conv}\left(\cup_{1 \leq i \leq k}\left\{ f + \frac{r^i}{\alpha_i}\right\}\cup \left\{f + \frac{r}{\hat{\beta}}\right\}\right) \textup{ is not lattice-free}$. Since $L_{\alpha}$ is bounded, we obtain that the set $\textup{conv}\left(\cup_{1 \leq i \leq k}\left\{ f + \frac{r^i}{\alpha_i}\right\}\cup \left\{f + \frac{r}{\hat{\beta}}\right\}\right)$ is bounded. Therefore there exists a finite number of integer points in its interior. Moreover if $\hat{\beta}^1 > \hat{\beta}^2$, then $\textup{conv}\left(\cup_{1 \leq i \leq k}\left\{ f + \frac{r^i}{\alpha_i}\right\}\cup \left\{f + \frac{r}{\hat{\beta}^1}\right\}\right) \subseteq \textup{conv}\left(\cup_{1 \leq i \leq k}\left\{ f + \frac{r^i}{\alpha_i}\right\}\cup \left\{f + \frac{r}{\hat{\beta}^2}\right\}\right)$. Thus, it is possible to choose $\tilde{\beta}$ such that $\gamma = \tilde{\beta}$.

\noindent Claim: $\sum_{i = 1}^{k}\alpha_is_i + \gamma s_{k+1}\geq 1$ is a  facet-defining inequality for $\textup{conv}(P([r^1, ..., r], f))$. By construction of $\gamma$, either the line segment between $f + \frac{r^1}{\alpha^1}$ and $f + \frac{r}{\gamma}$ or between $f + \frac{r^2}{\alpha^2}$ and $f + \frac{r}{\gamma}$ contains an integer point (that does not belong to $L_{\alpha}$). Let wlog $p = f + \lambda r + \lambda_1 r^1$ be this integer point where $\lambda >0$. Thus the inequality $\sum_{i = 1}^{k}\alpha_is_i + \gamma s_{k+1}\geq 1$ satisfies at equality the feasible point $(x,s) := (p, \bar{s}) \in P([r^1 \cdots r^{k+1}],f)$ where
\begin{eqnarray}
\bar{s}_i = \left\{\begin{array}{cl} 0 & \textup{if } i \neq 1, k +1\\
\lambda  & \textup{if } i = k +1\\
\lambda_1  & \textup{if } i = 1.
\end{array}\right.
\end{eqnarray}
The result follows from the fact that $\sum_{i = 1}^{k}\alpha_is_i \geq 1$ is facet-defining
for $\conv(P([r^1\cdots r^k],f)).$

Now there are two cases:
\begin{enumerate}
\item $f + \frac{r}{\gamma}$ is not integral: Then set $r^{k+1}:= r$, $\alpha_{k +1}:= \gamma$ and observe that $\textup{conv}\left(\cup_{1 \leq i \leq k + 1}\left\{ f + \frac{r^i}{\alpha_i}\right\}\right)$ is not a triangle of type $T^1$.
\item $f + \frac{r}{\gamma}$ is integral: Observe that the line segment between $f + \frac{r^1}{\alpha^1}$ and $f + \frac{r}{\gamma}$, and the line segment between $f + \frac{r^2}{\alpha^2}$ and $f + \frac{r}{\gamma}$ belong to the boundary of $\textup{conv}\left(\cup_{1 \leq i \leq k}\left\{ f + \frac{r^i}{\alpha_i}\right\} \cup \left\{f + \frac{r}{\gamma}\right\}\right)$. If any one of these line segments does not contain an integer point in its relative interior, then set $r^{k+1}:= r$, $\alpha_{k +1}:= \gamma$ and observe that $\textup{conv}\left(\cup_{1 \leq i \leq k + 1}\left\{ f + \frac{r^i}{\alpha_i}\right\}\right)$ is not a triangle of type $T^1$. If both these line segments contain an integer point in the relative interior, then let $p$ be one of these integer points. Observe that there exists a vector $r' \in \mathbb{Q}^2$ such that $\textup{cone}\{r^1, r^2, r'\} = \mathbb{R}^2$ and the ray $f + \lambda r'$, $\lambda \geq 0$ intersects the boundary of $\textup{conv}\left(\cup_{1 \leq i \leq k + 1}\left\{ f + \frac{r^i}{\alpha_i}\right\} \cup \left\{f + \frac{r}{\gamma}\right\}\right)$ at a non-integral point between the points $p$ and $f + \frac{r}{\gamma}$ (this is possible since the set $\{\bar{r} \in \mathbb{R}^2\,|\, \textup{cone}\{r^1, r^2, \bar{r}\} = \mathbb{R}^2\}$ is an open set). Now by setting $r^{k+1}:= r'$ and $\alpha_{k +1}:= \bar{\lambda}$ such that $f + \bar{\lambda}r'$ lies on the line segment between $p$ and $f + \frac{r}{\gamma}$, the result follows. \hfill $\square$
\end{enumerate}

\section{Properties of Split Rank}\label{propsection}
Section \ref{shapesection} deals with results that allow us to compare the split rank of two inequalities (for two different sets that may have some common columns $r^i$) based on the shape of the induced lattice-free set. Section \ref{standardsec} presents an operation on $P(R,f)$ under which the split ranks of related inequalities remain invariant.

\subsection{Split Rank and the Shape of Induced Lattice-free Set}\label{shapesection}
\begin{lemma}[Shape]\label{sameshape}
Let $\sum_{i = 1}^{k_1} \alpha_is_i \geq 1$ be a valid inequality for $\textup{conv}(P(R^a,f))$ with $R^a \in \mathbb{Q}^{2 \times k_1}$ and let $\sum_{i = 1}^{k_2} \beta_is_i \geq 1$ be a valid inequality for $\textup{conv}(P(R^b,f))$ with $R^b \in \mathbb{Q}^{2 \times k_2}$. We denote by $\eta_a$ and $\eta_b$ the split rank of $\sum_{i = 1}^{k_1} \alpha_is_i \geq 1$ and $\sum_{i = 1}^{k_2} \beta_is_i \geq 1$ respectively. If $\textup{cone}(R^b) = \mathbb{R}^2$ and $L_{\alpha} \subseteq L_{\beta}$, then $\eta_a \leq \eta_b$.
\end{lemma}
(Proof in Section \ref{proofshapesec}). Lemma \ref{sameshape} is straightforward to prove if $R^a$ and $R^b$ are the same set of columns, since the statement of Lemma \ref{sameshape} then implies that $\sum_{i = 1}^{k_2} \beta_is_i \geq 1$ dominates $\sum_{i = 1}^{k_1} \alpha_is_i \geq 1$. While the statement of Lemma \ref{sameshape} holds when $P(R^a,f)$ and $P(R^b,f)$ involve possibly different columns for the continuous variables, it is important to note that the two problems have same `right-hand-side' $f$.

\begin{figure}
\begin{center}
\includegraphics[width=0.35\linewidth]{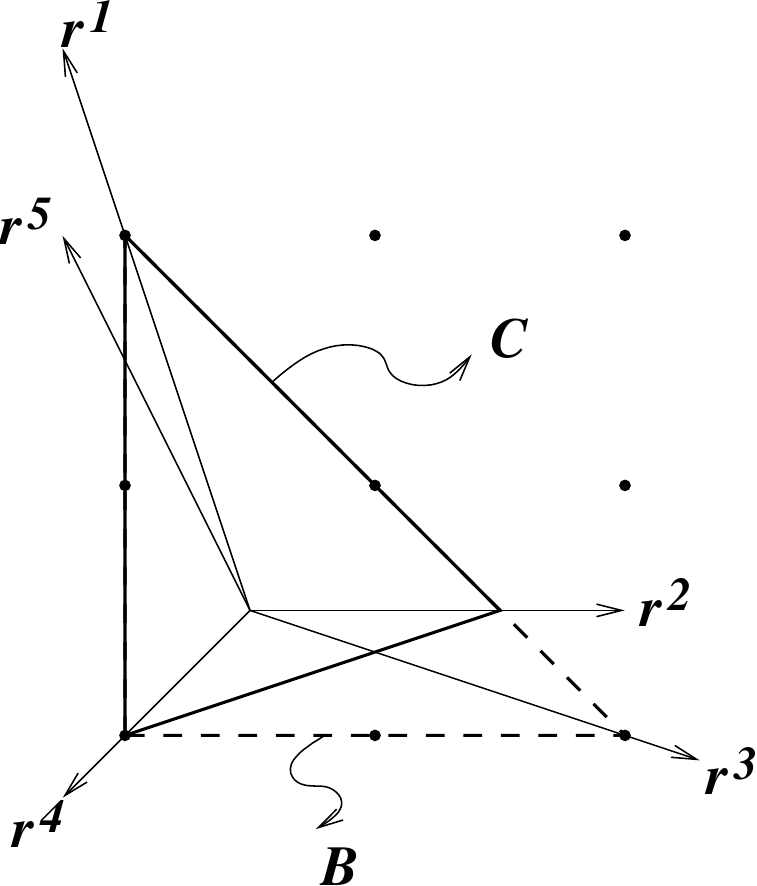}
\caption{Example \ref{toshowinduced}}
\label{FigTwoVarExample}
\end{center}
\end{figure}
\begin{example}\label{toshowinduced}
Consider the set
\begin{equation*}
\left\{x \in \mathbb{Z}^2, s \in \mathbb{R}_{+}^4 \mid  x = \left(\begin{array}{c} 0.5\\0.5 \end{array}\right) + \left(\begin{array}{c} -0.5\\1.5 \end{array}\right)s_1 + \left(\begin{array}{c} 1\\0 \end{array}\right)s_2 + \left(\begin{array}{c} 1.5\\-0.5 \end{array}\right)s_3 + \left(\begin{array}{c} -0.5\\-0.5 \end{array}\right)s_4 \right\}.
\end{equation*}
Let $B$ be the lattice-free triangle of type $T^1$ with vertices $(0,0)$, $(2,0)$, and $(0,2)$. Then using (\ref{pidefn}), $\phi(B)$ is
\begin{eqnarray}
\label{exinf}s_1 + s_2 + s_3 + s_4 \geq 1.
\end{eqnarray}
The induced lattice-free set of (\ref{exinf}) is $B$ and therefore the split rank of the inequality (\ref{exinf}) is not finite. Now consider the set where $s_3$ is dropped, i.e.,
\begin{equation*}
\left\{x \in \mathbb{Z}^2, s_1,s_2,s_4 \in \mathbb{R}_{+} \mid x = \left(\begin{array}{c} 0.5\\0.5 \end{array}\right) + \left(\begin{array}{c} -0.5\\1.5 \end{array}\right)s_1 + \left(\begin{array}{c} 1\\0 \end{array}\right)s_2 + \left(\begin{array}{c} -0.5\\-0.5 \end{array}\right)s_4\right\}.
\end{equation*}
Again using $B$ as the lattice-free triangle with vertices $(0,0)$, $(2,0)$, and $(0,2)$ we obtain the inequality $\phi(B)$
\begin{eqnarray}
\label{exinf1}s_1 + s_2 + s_4 \geq 1.
\end{eqnarray}
The induced lattice-free set of (\ref{exinf1}) is $C:=\textup{conv}\{(0,0), (1.5, 0.5), (0,2)\} \subsetneq B$ (See Figure \ref{toshowinduced}). The split rank of this inequality is finite; in fact 2. (The lower bound on the split rank is proven in Andersen \emph{et. al}~\cite{andersen:2007} and it can be verified that the inequality can be obtained by sequentially applying the disjunctions $(x_2\leq 0)\vee (x_2 \geq 1)$ and $(x_1\leq 0)\vee (x_1 \geq 1)$). Now consider the set
\begin{equation*}
\left\{x \in \mathbb{Z}^2, s_1, s_2, s_4, s_5 \in \mathbb{R}_{+} \mid x = \left(\begin{array}{c} 0.5\\0.5 \end{array}\right) + \left(\begin{array}{c} -0.5\\1.5 \end{array}\right)s_1 + \left(\begin{array}{c} 1\\0 \end{array}\right)s_2 + \left(\begin{array}{c} -0.5\\-0.5 \end{array}\right)s_4 + \left(\begin{array}{c} -0.5\\1 \end{array}\right)s_5\right\}.
\end{equation*}
Again using $B$ as the lattice-free triangle with vertices $(0,0)$, $(2,0)$, and $(0,2)$ we obtain the inequality
\begin{eqnarray}
\label{exinf2}s_1 + s_2 + s_4 + s_5 \geq 1.
\end{eqnarray}
The induced lattice-free set of (\ref{exinf2}) is again $C$. Therefore the split rank of (\ref{exinf2}) is also 2.
\end{example}

Besides illustrating the shape lemma, Example \ref{toshowinduced} also illustrates the fact that the finiteness of the split rank of an inequality $\sum_{i =1}^k\alpha_is_i \geq 1$ depends on its induced lattice-free set and not on a lattice-free convex set $B$ that is used to generate it (i.e. some $B$ such that $\phi(B) = \alpha$).

Notice that in the case of $L_{\alpha} \subsetneq L_{\beta}$, Lemma \ref{sameshape} does not imply that split rank of $\alpha$ wrt $(P(R^a,f))^0$ is strictly lesser that the split rank of $\alpha^b$ wrt $(P(R^b,f))^0$. Indeed, the following milder result implies that it is possible to have $L_{\alpha} \subsetneq L_{\beta}$ and yet have that the split rank of $\alpha$ wrt $(P(R^a,f))^0$ equal to the split rank of $\beta$ wrt $(P(R^b,f))^0$.
\begin{proposition}[General Lifting]\label{generallifting}
Let $\sum_{i = 1}^k\alpha_i s_i \geq \gamma$ be a valid inequality of $\textup{conv}(P(R,f))$ of split rank $\eta$. Then there exists $\alpha_{k+1} \geq 0$ such that $\sum_{i = 1}^k\alpha_is_i + \alpha_{k+1}s_{k+1} \geq \gamma $ is a valid inequality for $\textup{conv}(P([R \textup{  }r^{k+1}], f))$ and has a split rank of at most $\eta$ wrt $P([R \textup{  }r^{k+1}], f)^0$.
\end{proposition}

\subsubsection{Proof of Lemma \ref{sameshape} and Proposition \ref{generallifting}}\label{proofshapesec}

In Proposition \ref{generalprojection}, presented next, we analyze the split rank of an inequality when one variable is dropped from the description of the set. Proposition \ref{generalprojection} is used in the proof of Lemma \ref{sameshape} and also directly in the rest of the paper.
\begin{proposition}[Projection]\label{generalprojection}
Let $R = [r^1, ..., r^k]$. Let $\sum_{i = 1}^{k +1}\alpha_i s_i \geq \gamma$ be a valid inequality of $P([R \textup{  }r^{k+1}],f)$ of split rank $\eta$. Then $\sum_{i = 1}^k\alpha_is_i \geq \gamma $ has a split rank at most $\eta$ wrt $(P(R, f))^{0}$.
\end{proposition}
\textbf{Proof:} If $\eta = +\infty$, then the result is true. Therefore assume that $\eta$ is finite. We prove this result by proving that if $\textup{Proj}_{s, s_{k +1}}(P([R\, r^{k +1}], f))^{\eta} := \{(s,s_{k+1}) \in \mathbb{R}^k_{+}\times \mathbb{R}_{+}\,|\,As + A's_{k+1}\geq b\}$ for some $A \in \mathbb{Q}^{ g \times h}_{+}$, $A' \in \mathbb{Q}^{ g \times 1}_{+}$ and $b \in \mathbb{Q}^{g \times 1}_{+}$, then $\textup{Proj}_{s}(P(R, f))^{\eta} \subseteq \{s \in \mathbb{R}^k_{+}\,|\,As \geq b\}$. (The non-negativity of $A$ and $A'$ follows from Proposition \ref{inducedefn}). The proof is by induction on $\eta$. For $\eta = 0$ the statement is obvious. Assume that the statement is true for $\eta = 1,...,n - 1$.

Let $\textup{Proj}_{s, s_{k+1}}(P([R\, r^{k +1}], f))^{n - 1} := \{(s,s_{k+1} )\in \mathbb{R}^k_{+}\times \mathbb{R}_{+}\,|\,As + A's_{k+1}\geq b\}$. Let $\sum_{i = 1}^{k +1}\alpha_i s_i \geq \gamma$ be a valid inequality of $P([R\, r^{k+1}],f)^n$. This inequality must be dominated by a positive combination of a finite number of facet-defining inequalities $\sum_{i = 1}^{k +1}\alpha^{j}_i s_i \geq \gamma^{j}$ of $P([R\, r^{k+1}],f)^n$, where the inequality $\sum_{i = 1}^{k +1}\alpha^{j}_i s_i \geq \gamma^{j}$ is obtained by applying the disjunction $((\pi^j)^Tx \leq \pi^{j}_0) \vee ((\pi^j)^Tx \geq \pi^{j}_0 + 1)$ to $P([R\, r^{k +1}],f)^{n-1}$ ($\forall \,j$, $\pi^j \in \mathbb{Z}^2$ and $\pi^{j}_0 \in \mathbb{Z}$). 
Thus to prove $\textup{Proj}_{s}(P(R, f))^{\eta} \subseteq \{s \in \mathbb{R}^k_{+}|As \geq b\}$ it suffices to prove the following claim.

Claim: Let $\sum_{i = 1}^{k +1}\alpha_i s_i \geq \gamma$ be a valid inequality of $P([R\, r^{k+1}],f)^n$ obtained by applying the disjunction $(\pi^Tx \leq \pi_0) \vee (\pi^Tx \geq \pi_0 + 1)$ to $P([R\, r^{k+1}],f)^{n -1}$. Then $\sum_{i = 1}^k\alpha_is_i \geq \gamma $ is a valid inequality for $P(R,f)^n$: Note that the inequalities $\pi^Tx \leq \pi_0$ and $\pi^Tx \geq \pi_0 +1$ can be rewritten in terms of $s, s_{k +1}$ variables as $-\pi^TRs -\pi^Tr^{k+1}s_{k +1} \geq -\pi_0 + \pi^Tf$ and $\pi^TRs + \pi^Tr^{k+1}s_{k+1} \geq \pi_0 + 1 - \pi^Tf$ respectively. Therefore, the validity of the inequality $\sum_{i =1}^{k +1}\alpha_is_i \geq \gamma$ is equivalent to existence of $u^1, v^1 \in \mathbb{Q}^{1 \times (g +1)}_{+}$ such that
\begin{eqnarray}
\label{disju1}u^1 \left[\begin{array}{cc} A& A'\\
-\pi^TR&- \pi^Tr^{k+1}\end{array}\right] \leq [\alpha_1\,\alpha_2\,...\,\alpha_{k+1}] \\
\label{disju2}u^1 \left[\begin{array}{c} b\\
-\pi_0 + \pi^Tf\end{array}\right] \geq \gamma \\
\label{disjv1}v^1 \left[\begin{array}{cc} A & A'\\
\pi^TR & \pi^Tr^{k+1}\end{array}\right] \leq [\alpha_1\,\alpha_2\,...\,\alpha_{k+1}] \\
\label{disjv2}v^1 \left[\begin{array}{c} b\\
\pi_0 + 1 - \pi^Tf \end{array}\right] \geq \gamma
\end{eqnarray}

By the induction hypothesis,
\begin{eqnarray}\label{subset}
\textup{Proj}_{s}(P(R, f))^{n -1} \subseteq \{ s \in \mathbb{R}^k_{+} |As \geq b\}.
\end{eqnarray}

Now using (\ref{disju1}), (\ref{disju2}), (\ref{disjv1}), (\ref{disjv2}), and (\ref{subset}), $\sum_{i = 1}^k\alpha_is_i \geq \gamma $ is a valid inequality for $P(R,f)^n$. \hfill $\square$

Next we analyze the split rank for the case when one variable is added to the description of the set without changing the shape of the induced lattice-free set.

\begin{proposition}[Simple Lifting]\label{nvar}
Let $\sum_{i = 1}^k\alpha_i s_i \geq 1$ be a valid inequality of $\textup{conv}(P(R,f))$ of split rank $\eta$. Let $r^{k + 1} \in \textup{cone}_{r \in R}\{r\}$. Then $\sum_{i = 1}^k\alpha_i s_i  + \alpha_{k +1} s_{k +1}\geq \gamma$ is a valid inequality of $\textup{conv}(P([R \,r^{k+1}],f))$ of split rank at most $\eta$ where $\alpha_{k +1}$ is such that $f + \frac{r^{k+1}}{\alpha_{k +1}} \in \partial L_{\alpha}$.
\end{proposition}
\textbf{Proof:} If $\eta = + \infty$, then the result holds. We assume that $\eta$ is finite.

Since $\alpha_{k +1}$ is such that $f + \frac{r^{k+1}}{\alpha_{k +1}} \in \partial L_{\alpha}$, there exist $r^a, r^b \in R$ such that $r^{k +1} = \lambda^a r^a + \lambda^b r^b$ and $\alpha_{k +1} = \lambda^a\alpha^{a} + \lambda^b\alpha^{b}$ where $\lambda^a, \lambda^b \geq 0$. WLOG assume that $a = 1$ and $b = 2$.

Claim: If $(\bar{x}, \bar{s}_1, \bar{s}_2, ..., \bar{s}_k, \bar{s}_{k+1}) \in P([R \textup{  }r^{k +1}],f)^p$, then $(\bar{x}, \bar{s}_1 + \lambda^1\bar{s}_{k+1}, \bar{s}_2+ \lambda^2\bar{s}_{k+1}, \bar{s}_3, ..., \bar{s}_k) \in P(R ,f)^p$. The statement is true for $p = 0$. Assume the claim is true for $p = 1, ..., n$. We need to show that if $(\bar{x}, \bar{s}_1, \bar{s}_2, ..., \bar{s}_k, \bar{s}_{k+1}) \in P([R \textup{  }r^{k +1}],f)^{n +1}$, then $(\bar{x}, \bar{s}_1 + \lambda^1\bar{s}_{k+1}, \bar{s}_2+ \lambda^2\bar{s}_{k+1}, ..., \bar{s}_k) \in P(R ,f)^{n+1}$. As $P([R \textup{  }r^{k +1}],f)^{n +1} \subseteq P([R \textup{  }r^{k +1}],f)^{n}$, we obtain $(\bar{x}, \bar{s}_1, \bar{s}_2, ..., \bar{s}_k, \bar{s}_{k+1}) \in P([R \textup{  }r^{k +1}],f)^{n}$. By the induction argument, $(\bar{x}, \bar{s}_1 + \lambda^1\bar{s}_{k+1}, \bar{s}_2+ \lambda^2\bar{s}_{k+1}, ..., \bar{s}_k) \in P(R ,f)^{n}$. Now consider any disjunction of the form $(\pi^Tx \leq \pi_0) \vee (\pi^Tx \geq \pi_0 +1)$ applied to $P(R ,f)^{n}$. We obtain the following cases (let $P(R,f)^{n}_{\pi, \pi_0} := \textup{conv}( (P(R,f)^{n} \cap \{(x,s) \in \mathbb{R}^2\times \mathbb{R}^k\,|\, \pi^Tx \leq \pi_0\})\cup (P(R,f)^n \cap \{(x,s) \in \mathbb{R}^2\times \mathbb{R}^k\,|\, \pi^Tx \geq \pi_0 + 1\}))$):
\begin{enumerate}
\item $\pi^T\bar{x} \leq \pi_0$ or $\pi^T\bar{x} \geq \pi_0 + 1$. Then $(\bar{x}, \bar{s}_1 + \lambda^1\bar{s}_{k+1}, \bar{s}_2+ \lambda^2\bar{s}_{k+1}, ..., \bar{s}_k) \in P(R,f)^n_{ \pi, \pi_0}$.
\item $\pi_0 < \pi^T\bar{x} < \pi_0 + 1$. Since $(\bar{x}, \bar{s}_1, \bar{s}_2, ..., \bar{s}_k, \bar{s}_{k+1}) \in P([R \textup{  }r^{k +1}],f)^{n +1}$ and $\pi_0 < \pi^T\bar{x} < \pi_0 + 1$, there exist two points $(x^1, s^1_1, s^1_2,s^1_3, ..., s^1_{k+1}), (x^2, s^2_1, s^2_2,s^2_3, ..., s^2_{k+1}) \in P([R \textup{  }r^{k +1}],f)^{n}$ such that
    \begin{itemize}
    \item $\pi^Tx^1 \leq \pi_0$, $\pi^Tx^2 \geq \pi_0 +1$
    \item $(\bar{x}, \bar{s}_1, \bar{s}_2, ..., \bar{s}_k, \bar{s}_{k+1})$ is a convex combination of $(x^1, s^1_1, s^1_2,s^1_3, ..., s^1_{k+1})$ and $(x^2, s^2_1, s^2_2,s^2_3, ..., s^2_{k+1})$.
    \end{itemize}
    This implies however by the induction argument that $( x^1, s^1_1 + \lambda^1s^1_{k+1}, s^1_2 + \lambda^2s^1_{k+1},s^1_3, ...), (x^2, s^2_1 + \lambda^1s^2_{k+1}, s^2_2+ \lambda^2s^2_{k+1},s^2_3, ...) \in P(R ,f)^{n}$ such that
    \begin{itemize}
    \item $\pi^Tx^1 \leq \pi_0$, $\pi^Tx^2 \geq \pi_0 +1$
    \item $(\bar{x}, \bar{s}_1 + \lambda^1\bar{s}_{k+1}, \bar{s}_2+ \lambda^2\bar{s}_{k+1}, ..., \bar{s}_k)$ is a convex combination of $(x^1, s^1_1 + \lambda^1s^1_{k+1}, s^1_2 + \lambda^2s^1_{k+1},s^1_3, ..., s^1_k)$ and $(x^2, s^2_1 + \lambda^1s^2_{k+1}, s^2_2+ \lambda^2s^2_{k+1},s^2_3, ..., s^2_k)$.
    \end{itemize}
    In other words, $(\bar{x}, \bar{s}_1 + \lambda^1\bar{s}_{n+1}, \bar{s}_2+ \lambda^2\bar{s}_{n+1}, ..., \bar{s}_k) \in P(R,f)^n_{\pi, \pi_0}$.
\end{enumerate}

Now we return to the proof of the proposition: Assume by contradiction that the inequality $\sum_{i = 1}^k\alpha_is_i + (\lambda^1\alpha^{1} + \lambda^2\alpha^{2})s_{k+1} \geq \gamma $ has a split rank greater than $\eta$. Therefore, there exists $(\bar{x}, \bar{s}_1, \bar{s}_2, ..., \bar{s}_k, \bar{s}_{k+1}) \in P([R \textup{  }r^{k +1}],f)^\eta$ such that $\sum_{i=1}^{k}\alpha_i\bar{s}_i + (\lambda^1 \alpha_1 + \lambda^2 \alpha_2)\bar{s}_{k+1} < \gamma$. This implies, $\alpha_1(\bar{s}_1 + \lambda^1\bar{s}_{k+1}) + \alpha_2(\bar{s}_2 + \lambda^2\bar{s}_{k +1}) +\sum_{i=3}^{k}\alpha_i\bar{s}_i < \gamma$. However note now from the claim that $(\bar{x}, \bar{s}_1 + \lambda^1\bar{s}_{k+1}, \bar{s}_2+ \lambda^2\bar{s}_{k+1}, \bar{s}_3, ..., \bar{s}_k) \in P(R ,f)^\eta$. This implies that the inequality $\sum_{i = 1}^k\alpha_is_i \geq \gamma $ has a rank greater than $\eta$ wrt $(P(R,f))^0$, a contradiction. \hfill $\square$

Using Proposition \ref{generalprojection} and \ref{nvar}, Lemma \ref{sameshape} can be verified.
\newline\newline\textbf{Proof of Lemma \ref{sameshape}:} Let $R^c = R^a \cup R^b$ and consider the set $P(R^c,f)$. Then the inequality,
\begin{eqnarray}
\label{tempcut}\sum_{r^i \in R^a \cup R^b}\gamma(r^i)s_i \geq 1
\end{eqnarray}
where $\gamma := \phi(L_{\beta})$ is valid for $P(R^c,f)$. Note that by definition $L_{\gamma}\subseteq L_{\beta}$. However, since every column of $R^b$ belongs to $R^c$ and the corresponding coefficients of $\gamma$ and $\beta$ are equal, every vertex (resp. ray) of $L_{\beta}$ is a vertex (resp. ray) of $L_{\gamma}$. Thus, $L_{\gamma}\supseteq L_{\beta}$ or $L_{\gamma}= L_{\beta}$.

Let $\eta_c$ be the split rank of (\ref{tempcut}) wrt $(P(R^c, f))^0$.

Now starting from $P(R^c,f)$ and the inequality (\ref{tempcut}), by the application of Proposition \ref{generalprojection} iteratively for every column of $R^a\setminus R^b$, we obtain that $\eta_c \geq \eta_b$. However, since $L_{\gamma} = L_{\beta}$, by the application of Proposition \ref{nvar}, (starting from $P(R^b,f)$ and $\sum_{i = 1}^{k_2}\beta_is_i \geq 1$ and then simple lifting all the columns in $R^a\setminus R^b$), we obtain $\eta_c \leq \eta_b$. Therefore, $\eta_c = \eta_b$.

Again by application of Proposition \ref{generalprojection} and the fact that $\alpha(r^i) \geq \gamma(r^i)$ $\forall r^i \in R^a$ (since $L_{\alpha} \subseteq L_{\gamma}$), we obtain that $\eta_c \geq \eta_a$. This completes the proof. \hfill $\square$
\newline\newline\textbf{Proof of Proposition \ref{generallifting}:} Proposition \ref{generallifting} is proven by showing that if $\textup{Proj}_{s}(P([r^1, ..., r^k], f))^{\eta} := \{s \in \mathbb{R}^k_{+}\,|\,As \geq b\}$ for some $A \in \mathbb{Q}^{ g \times h}_{+}$ and $b \in \mathbb{Q}^{g \times 1}_{+}$, then $\textup{Proj}_{s, s_{k+1}}(P([r^1, ..., r^k, r^{k+1}], f))^{\eta} \subseteq \{(s,s_{k+1}) \in \mathbb{R}^k_{+}\times \mathbb{R}_{+})\,|\,As + A's_{k+1} \geq b\}$ for some $A' \in \mathbb{Q}^{ g \times 1}_{+}$. This is similar to the proof of Proposition \ref{generalprojection}.

\subsection{For Standardization}\label{standardsec}
The following result allows us to consider `standard' triangles and quadrilateral. Related observations for split cuts were made in Dash et al.~\cite{dashGL:2007}.
\begin{observation}[Integral Translation and Unimodular Bijection]\label{transprop}
Let $w \in \mathbb{Z}^2$ and $M \in \mathbb{Z}^{2 \times 2}$ be a unimodular matrix. Then
\begin{enumerate}
\item A valid inequality $\sum_{i=1}^k\alpha_is_i \geq \gamma$ for $P(R,f)$ is facet-defining for $P(R,f)$ if and only if $\sum_{i = 1}^k\alpha_is_i \geq \gamma$ is valid and facet-defining for $P(MR, M(f + w))$.
\item The split rank of $\sum_{i = 1}^k\alpha_is_i \geq \gamma$ wrt $ P(R,f)^0$ is $\eta$ if and only if the split rank of $\sum_{i = 1}^k\alpha_is_i \geq \gamma$ wrt $P(MR,M(f + w))^0$ is $\eta$.
\end{enumerate}
\end{observation}
\section{Two Variable Problems}\label{twovarsec}
\begin{proposition}\label{2varthm}
Consider a non-trivial facet-defining inequality
\begin{math}
\alpha_1 s_1 + \alpha_2 s_2 \geq 1
\end{math}
for $\textup{conv}(P([r^1, r^2], f))$. Then its split rank is at most two.
\end{proposition}
\begin{figure}
\begin{center}
\includegraphics[width=0.75\linewidth]{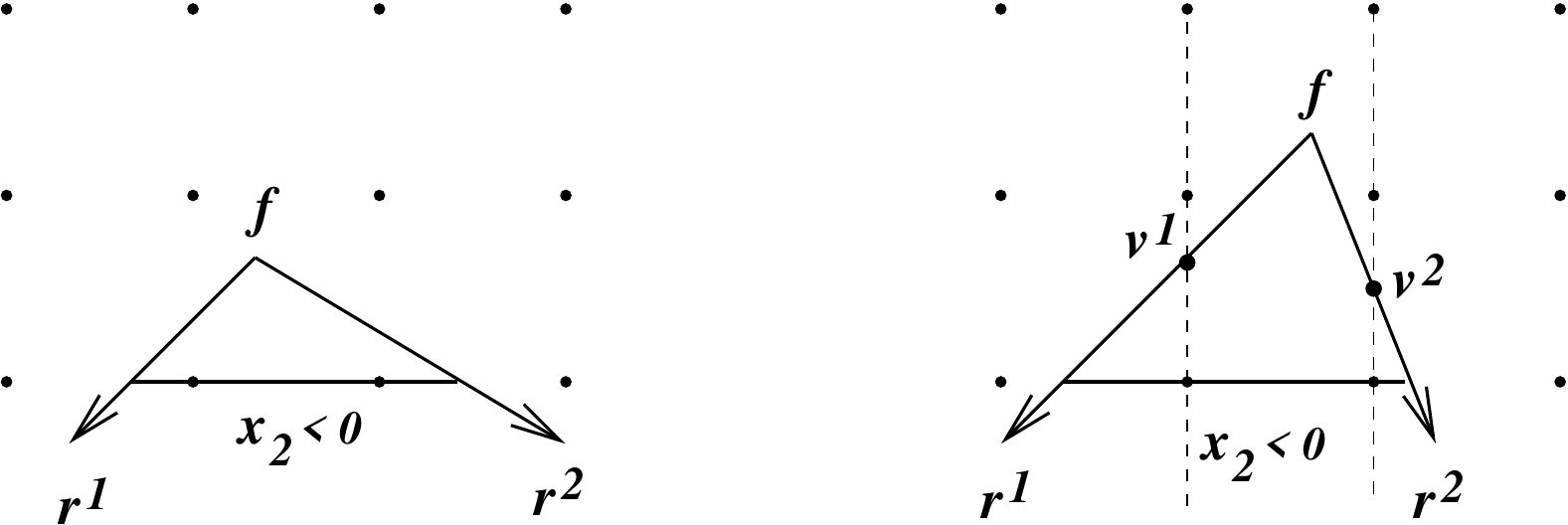}
\caption{The idea behind the proof of Proposition \ref{2varthm}}
\label{FigTwoVar}
\end{center}
\end{figure}
\textbf{Proof:} Since $\alpha_1s_1+\alpha_2s_2\geq 1$ is facet-defining, it satisfies at equality two or more feasible points of $P([r^1, r^2],f)$. By suitable integral translation and unimodular transformation, we can assume that
\begin{enumerate}
\item $f_2 > 0$.
\item $0 \leq f_1 \leq 1$.
\item The $x$ component of the feasible points that satisfy $\alpha_1 s_1 + \alpha_2 s_2 \geq 1$
at equality are $(0,0)$, $(1, 0)$ , ..., $(g, 0)$ where $h \leq 0$ and $g \geq 1$.
\item $r^1_2, r^2_2 \leq 0$.
\item The inequality $\alpha_1 s_1 + \alpha_2 s_2 \geq 1$ is equivalent to $x_2 \leq 0$.
\end{enumerate}
See Figure \ref{FigTwoVar} for an illustration. There are two cases.

If $f_2 < 1$ (see left frame in Figure \ref{FigTwoVar}), then $P([r^1,r^2],f)\cap \{x\mid x_2\geq 1\}=\emptyset$ because $r^1_2,r^2_2\leq 0$. Therefore $x_2\leq 0$ is valid for $\conv((P([r^1,r^2],f) \cap \{x\mid x_2\leq 0\} )\cup(P([r^1,r^2],f)\cap \{x\mid x_2\geq 1\}))$.

Now consider the case where $f_2 \geq 1$ (see right frame in Figure \ref{FigTwoVar}). Let $v^j:= \{x \in \R^2 \mid x=f+r^j s_j, s_j \geq 0\}\cap \{x\in \R^2 \mid x_1=j-1\}, j=1,2.$

Claim: $v^1_2<1$ and {$v^2_2<1$}: Assume $v^1_2\geq 1$. Then, we have $(0, v^1_2),(0, 0)\in \textup{Proj}_x((P(R,f))^0)$. By convexity, we conclude that $(0, 1)\in \textup{Proj}_x((P(R,f))^0)$ which is the required contradiction since $x_2 \leq 0$ is a valid inequality for $\textup{conv}(P(R,f))$. Similarly, we can verify that {$v^2_2<1$}.

Observe that $(x, s_1, s_2):= (z^1, \lambda_1, 0)$ (for a suitable $\lambda_1 > 0$) is the only vertex of the set $Q^{\leq}:=P(R, f)^0 \cap \{(x,s)\in \mathbb{R}^2\times \mathbb{R}^2\,|\,x_1 \leq 0\}$ and $(x, s_1, s_2):= (z^2, 0, \lambda_2)$ is the only vertex of the set $Q^{\geq}:=P(R, f)^0 \cap \{(x,s)\in \mathbb{R}^2\times \mathbb{R}^2\,|\,x_1 \geq 1\}$. The extreme rays of $\conv(Q^{\leq}\cup Q^{\geq})$ are $(r^1,(1, 0))$ and $(r^2,(0, 1))$. As $r^1_2,r^2_2\leq 0$, there exists  $0<\delta \leq 1$ such that $x_2 \leq 1-\delta$ is  valid for $\conv(Q^{\leq}\cup Q^{\geq})$. Define $\bar Q:=\{(x,s)\in \R^2\times \R^3_+\mid x=f+Rs, x_2\leq 1-\delta\}.$ Observe that $x_2\leq 0$ is valid for $\conv(\{\bar Q\cap \{(x,s)\mid x_2\leq 0\} \cup \bar Q\cap \{(x,s)\mid x_2\geq 1\}).$\cqfd

\section{Three Variable Problems}\label{threevarsec}
In this section, we consider the split rank of facet-defining inequalities $\sum_{i = 1}^3\alpha_is_i \geq 1$ for $\textup{conv}(P([r^1, r^2, r^3],f))$ where $\textup{cone}\{r^1, r^2, r^3\} = \mathbb{R}^2$.

\subsection{$T^{2B}$}\label{t2csec}
We prove the following result in this section.
\begin{proposition}\label{theoremt2c}
The split rank of an inequality whose induced lattice-free set is a triangle of type $T^{2B}$ is finite.
\end{proposition}

In Section \ref{standardt2Bsec} we discuss the standard triangle of type $T^{2B}$. In Section \ref{t2cs1} we present some useful definitions and an outline of the proof of Proposition \ref{theoremt2c}. There are four main subcases in the proof of Proposition \ref{theoremt2c} that differ in the details. These cases correspond to Sections \ref{t2cc1} - \ref{t2cc3}.

\subsubsection{Standardization of triangles of type $T^{2B}$}\label{standardt2Bsec}
\begin{center}
\begin{figure}[!h]
\begin{center}
\includegraphics[width=0.4\textwidth]{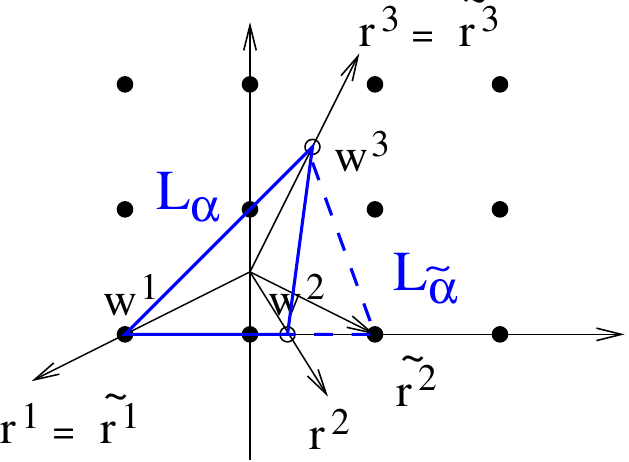}
\end{center}
\caption{The case where $g = 0$}
\label{t2bfig}
\end{figure}
\end{center}
Let $\sum_{i=1}^k \alpha_i s_i \geq 1$ be a facet-defining inequality for $\conv(P([r^1,r^2,r^3],f))$ such that $L_{\alpha}$ is  a triangle of type $T^{2B}$. By a suitable integral translation and unimodular transformation (Dey and Wolsey~\cite{deywolsey2008a}), we can assume that (1) the vertices of $L_{\alpha}$ are (a) $w^1:= (-\delta , 0)$ where $0 < \delta \leq 1$ and $w^1=f+\lambda_1r^1, \lambda_1\geq 0,$ (b) $w^2:= (g + \epsilon, 0)$ where $0 \leq \epsilon < 1$, $g\in \mathbb{Z}_{+}$ and $w^2=f+\lambda_2r^2,\lambda_2\geq 0$, (c) $w^3:= (\bar{x}, \bar{y})$ where $\bar{y} > 1$ and $0 < \bar{x} < 1$, and $w^3=f+\lambda_3r^3,\lambda_3\geq 0$. (2) The side $w^1w^3$ of $L_{\alpha}$ contains the integer point $(0,1)$ in its relative interior. (3) The side $w^1w^2$ of $L_{\alpha}$ contains multiple integer points. (4) The side $w^2w^3$ of $L_{\alpha}$ does not contain any integer point in its relative interior.

However while $w^2_1 = g + \epsilon$ can be less than $1$, it is convenient to work with triangles with $w^2_1 \geq 1$. Consider the case where $g = 0$ (see Figure \ref{t2bfig}). In this case it is possible to consider a different set $P([\tilde{r}^1, \tilde{r}^2, \tilde{r}^3], f)$ and a corresponding inequality $\sum_{i = 1}^3\tilde{\alpha}_is_i \geq 1$ where $\tilde{r}^1 = r^1$, $\tilde{r}^3 = r^3$, $\tilde{r}^2 = (1 - f_1, -f_2)$, and $\tilde{\alpha}_1 = {\alpha}_1$, $\tilde{\alpha}_3 ={\alpha}_3$, and $\tilde{\alpha}_2 =1$. Observe that $L_{\alpha} \subseteq L_{\tilde{\alpha}}$ (see Figure \ref{t2bfig}) and therefore by Lemma \ref{sameshape}, the split rank of $\alpha$ wrt $(P(R,f))^0$ is less than or equal to the split rank of $\tilde{\alpha}$ wrt $(P(\tilde{R},f))^0$.

Hence we consider the standard $T^{2B}$ as presented next.

\begin{proposition}[Standard $T^{2B}$] \label{standardt2C}
Let $\sum_{i=1}^k \alpha_i s_i \geq 1$ be a facet-defining inequality for
$\conv(P([r^1,r^2,r^3],f))$ such that $L_{\alpha}$ is  a triangle of type $T^{2B}$.
By a suitable integral translation and unimodular transformation, we can assume that
\begin{enumerate}
\item The vertices of $L_{\alpha}$ are
\begin{enumerate}
\item $w^1:= (-\delta , 0)$ where $0 < \delta \leq 1$ and $w^1=f+\lambda_1r^1, \lambda_1\geq 0,$
\item $w^2:= (g + \epsilon, 0)$ where $0 \leq \epsilon < 1$, $g\in \mathbb{Z}_{+}\setminus\{0\}$ and
$w^2=f+\lambda_2r^2,\lambda_2\geq 0$,
\item $w^3:= (\bar{x}, \bar{y})$ where $\bar{y} > 1$ and $0 < \bar{x} < 1$,
and $w^3=f+\lambda_3r^3,\lambda_3\geq 0$.
\end{enumerate}
\item The side
$w^1w^3$ of $L_{\alpha}$  contains the integer point $(0,1)$ in its relative interior.
\item The side $w^2w^3$ of $L_{\alpha}$ does not contain any integer point in its relative interior.
\item $r^1_2<0,r^1_1<0,r^2_2<0,r^2_1\geq 0$ and $r^3_2>0$.
\end{enumerate}
\end{proposition}

\subsubsection{Some Definitions and Proof Outline}\label{t2cs1}

Before outlining the proof of Proposition \ref{theoremt2c}, we present a couple of definitions linking a point $x\in \R^2$ to $s \in \mathbb{R}^k_{+}$ such that $x=f+Rs$. These definitions simplify the presentation of the proofs in the remainder of the paper.
\begin{definition}[$\lambda$ Notation]
Let $\bar x\in \R^2$ such that there exists $i\in \{1,\ldots, k\}$ and $\bar \lambda \geq 0$
with $\bar x=f+\bar \lambda r^i.$ For convenience, we denote
$$ \lambda(\bar x) := \bar \lambda.$$
\end{definition}

\begin{definition}[Minimal Representation]
Consider the set $P([r^1, ..., r^k],f)$ and let $\bar{x} \in \mathbb{R}^2$. Then
\begin{enumerate}
\item If $\bar{x} =f+ \lambda_ir^i + \lambda_j r^j$ with $\lambda_i, \lambda_j >0$
$i\neq j$ and $r^i\neq \nu r^j$ for all $\nu\in \R$, we define $\mathcal{M}^{i,j}(\bar{x}) \in \mathbb{R}^k_{+}$ as
    \begin{eqnarray}
    \mathcal{M}^{i,j}_l(\bar{x}) = \left\{\begin{array}{cl} 0 & \textup{if }l \neq i, j \\
    \lambda_i  & \textup{if }l = i \\
    \lambda_j  & \textup{if }l = j \\
    \end{array}\right.
    \end{eqnarray}
\item If $\bar{x} = f + \bar{\lambda}r^i$ with $\bar \lambda \geq 0$ (i.e. $\lambda(\bar{x}) = \bar{\lambda}$) , we define $\mathcal{M}^{i,i}(\bar{x}) \in \mathbb{R}^k_{+}$ as $\mathcal{M}^{i,i}(\bar{x})_l = \bar{\lambda}$ for $l = i$ and $\mathcal{M}^{i,i}(\bar{x})_l = 0$ for $l \neq i$.
\end{enumerate}
\end{definition}
Note that when there are only three variables and $\textup{cone}\{r^1, r^2, r^3\} = \mathbb{R}^2$, every  $\bar{x} \in \mathbb{R}^2$ satisfies exactly one of the two cases in the definition of minimal representation. Moreover, in the first case if $\bar{x} =f+ \lambda_ir^i + \lambda_j r^j$ and $\lambda_i, \lambda_j >0$, then $i$ and $j$ are unique as well. Therefore if $k=3$, we use $\mathcal{M}(\bar{x})$ to represent the unique minimal representation for each vector $\bar{x}$.

\begin{proposition}\label{therelight}
Let $\beta \in \mathbb{R}^k_{+}$ have positive components for some indices $i,j \in \{1,..., k\}$. Then $\sum_{l = 1}^k\beta_l\mathcal{M}^{i,j}_l(\bar{x}) < 1$ if and only if $\bar{x} \in \textup{rel.int}(\textup{conv}\{f, f + \frac{r^i}{\beta_i}, f+  \frac{r^j}{\beta_j}\})$.
\end{proposition}
\textbf{Proof:} We present the proof for the case where $i \neq j$. The proof is similar for the other case. Observe that $\bar{x} \in \textup{int}(\textup{conv}\{f, f + \frac{r^i}{\beta_i}, f+  \frac{r^j}{\beta_j}\})$ iff $x = \mu_0f + \mu_i (f + \frac{r^i}{\beta_i}) + \mu_j(f + \frac{r^i}{\beta_i})= f + \frac{\mu_i}{\beta_i}r^i + \frac{\mu_j}{\beta_j}r^j$ where $1 > \mu_0 >0$. Also since $r^i \neq \nu r^j$ for all $\nu \in \mathbb{R}$, we obtain that $\mathcal{M}^{i,j}(\bar{x})_i = \frac{\mu_i}{\beta_i}$, $\mathcal{M}^{i,j}(\bar{x})_j = \frac{\mu_j}{\beta_j}$, and $\mathcal{M}^{i,j}(\bar{x})_l = 0$ if $l \neq i,j$. Thus, $\bar{x} \in \textup{int}(\textup{conv}\{f, f + \frac{r^i}{\beta_i}, f+  \frac{r^j}{\beta_j}\})$ iff $\sum_{l = 1}^k\beta_l\mathcal{M}^{i,j}_l(\bar{x}) = 1 - \mu_0 <1$.\hfill $\square$

Note that when $k =3$ and $\textup{cone}(r^1, r^2, r^3) = \mathbb{R}^2$, Proposition \ref{therelight} is equivalent to $\sum_{l = 1}^3\beta_i\mathcal{M}^{i,j}_l(\bar{x}) < 1$ if and only if $\bar{x} \in \textup{int}(L_{\beta})$.

In the proof of Proposition \ref{theoremt2c}, we typically consider $(P(R,f))^{0}$ along with one inequality, i.e., the set
\begin{equation}
\label{Qhere}Q:= \{(x,s)\in \R^2\times \R^3_+\mid x = f + r^1s_1 + r^2s_2 + r^3s_3,
\alpha_1s_1 + \alpha_2s_2 + \alpha_3s_3 \geq 1 \}.
\end{equation}
Corresponding to some disjunction $(\pi^Tx \leq \pi^0)\vee (\pi^Tx \geq \pi^0 +1)$ we consider the two sets
\begin{eqnarray}
\label{Qlessgreat}Q^{\leq}:=Q\cap \{(x,s)\mid \pi^T x\leq \pi_0\}, \quad Q^{\geq}:=Q\cap \{(x,s)\mid \pi^T x\geq \pi_0+1\}.
\end{eqnarray}
We would like to prove that an inequality $\beta_1s_1 + \beta_2s_2 + \beta_3s_3 \geq 1$ (where $\beta_1, \beta_2, \beta_3 > 0$), is valid for $\conv(Q^{\leq}\cup Q^{\geq})$. Note that the support of the $s$-component of the vertices of $Q^{\leq}$ and $Q^{\geq}$ is at most 2. More precisely the following observation can be verified.
\begin{observation}\label{plx}
Let $Q^{\leq}$ be as in (\ref{Qlessgreat}). Then the vertices of $Q^{\leq}$ are of the form $(\bar{x}, \mathcal{M}(\bar{x}))$ where $\bar{x}$ is of the form:
\begin{enumerate}
\item $f + \lambda r^i$, $\lambda > 0$, or
\item the intersection points of the boundary of $L_{\alpha}$ and the line segment $\pi_1x_1 + \pi_2x_2 = \pi_0$.
\end{enumerate}
\end{observation}

Note that the extreme rays of $\textup{conv}(Q^{\leq} \cup Q^{\geq})$ are $(r^1, e^1)$, $(r^2, e^2)$, $(r^3,e^3)$. The $s$-component of these rays satisfy $\beta_1s_1 + \beta_2s_2 + \beta_3s_3 \geq 0$. Therefore using Proposition \ref{therelight} and the above observation, checking validity of the inequality $\beta_1s_1 + \beta_2s_2 + \beta_3s_3 \geq 1$ is simplified and is recorded in the next Proposition.

\begin{proposition}\label{3varheart}
Let $Q, Q^{\leq}, Q^{\geq}$ be as in (\ref{Qhere}) and (\ref{Qlessgreat}). Then $\beta_1s_1 + \beta_2s_2 + \beta_3s_3 \geq 1$ is a valid inequality for $\conv(Q^{\leq}\cup Q^{\geq})$ if for every vertex $(\bar{x}, \mathcal{M}(\bar{x}))$ of $Q^{\leq}$ and $Q^{\geq}$, $\bar{x} \notin \textup{int}(L_{\beta})$.
\end{proposition}

Since we will repeatedly reference the $x$-components of the vertices of either $Q^{\leq}$ or $Q^{\geq}$ to check the validity of an inequality, for simplicity we will refer to the $x$-component of the vertices of $Q^{\le}$ and $Q^{\ge}$ as the $x$-vertices.

\emph{Outline of the proof of Proposition \ref{theoremt2c}:}
Apply a sequence of two disjunctions $(x_1\leq 0) \vee (x_1\geq 1)$ and $(x_2\leq 0) \vee (x_2\geq 1)$ successively. At each step, select one inequality valid for $\conv(Q^\leq \cup Q^\geq)$ (ignoring all the other inequalities) and then proceed with the next disjunction. We will show that this procedure converges to the desired inequality in a finite number of steps. Observe that as we keep exactly one inequality at each step, the validity of the inequality that is selected can be checked by the use of Proposition \ref{3varheart}.

We distinguish between four cases that differ slightly in the sequence of disjunctions used for the proof of convergence:
\begin{enumerate}
\item $0 < f_1 \leq 1$ and $r^3_1 < 0$.
\item $f_2 \leq 1$ and $r^3_1 < 0$.
\item $r^3_1 = 0$.
\item $r^3_1 > 0$.
\end{enumerate}
It can be verified that all scenarios are covered in the above four cases. The following notation is used throughout this section.

\begin{notation}\label{allcommonnotation}
\begin{enumerate}[(i)]
\item We define $Q^{[0]}:= (P(R,f))^0.$ Let $\sum_{j = 1}^3\beta^{[i]}_js_j \geq 1$ be the inequality obtained in step $i - 1$. We define $Q^{[i]}:= \{(x,s)\in \R^2\times \R^3_+ \mid x=f+Rs, \sum_{j=1}^3 \beta^{[i]}_j s_j \geq 1\}.$
\item Corresponding to each $r^i$, we define the intersection points
\begin{align*}
v^i&:= \{x \in \R^2 \mid x = f+ \lambda r^i, \lambda \geq 0, x_1 =0 \text{ or } x_1 = 1\}.
\end{align*}
for all $i=1,2,3.$ Note that if $0<f_1<1$, then $v^i$ is uniquely determined for all $i$.
\item If $f_2\geq 1$, we denote $\bar w^2 := \{ x\in \R^2 \mid x=f+\lambda r^2, \lambda \geq 0,
x_2=1\}.$
\item $q = (q_1, 1)$ is the intersection point of the line segment $w^2w^3$
with the line $\{x \in \R^2 \mid x_2=1\}$. (Remember that $w^1$, $w^2$, and $w^3$ are the vertices of $L_{\alpha}$).
\end{enumerate}
\end{notation}
Finally, we introduce a construction that is useful in the presentation of the proof of Proposition \ref{theoremt2c}.

\begin{construction}[$\triangle$]
Let $\{i,j,k\}$ be a permutation of $\{1,2,3\}$. Let $X,Y,Z \in \mathbb{R}^2$ be three affinely independent points such that $X, Y \in f + \textup{cone}\{r^i, r^j\}$ and $Z\in f + \textup{cone}\{r^i, r^k\}$. Suppose that there exists $p^i$, $p^j$, $p^k \in \mathbb{R}^2$ such that
\begin{enumerate}[(i)]
\item $p^i$ is the intersection point of the line $XY$ with the ray $\{x\in \mathbb{R}^2\,|\, f+ \lambda r^i, \lambda \geq 0\}$,
\item $p^j$ is the intersection point of the line $XY$ with the ray $\{x\in \mathbb{R}^2\,|\, f+ \lambda r^j, \lambda \geq 0\}$,
\item $p^k$ is the intersection point of the line $p^iZ$ with the ray $\{x\in \mathbb{R}^2\,|\, f+ \lambda r^k, \lambda \geq 0\}$.
\end{enumerate}
Then we denote $\triangle(XYZ):= \conv\{p^i,p^j,p^k\}$. Note that the ordering of the points $X$, $Y$, and $Z$ in the notation $\triangle(XYZ)$ is not relevant. Therefore, we interchangeably use $\triangle(XZY)$ or $\triangle(ZYX)$ to denote $\triangle(XYZ)$.
\end{construction}

\subsubsection{Case 1: $0 \leq f_1 \leq 1$ and $r^3_1 < 0$.}\label{t2cc1}
\begin{center}
\begin{figure}[htbp]
\begin{center}
\includegraphics[width=0.9\textwidth]{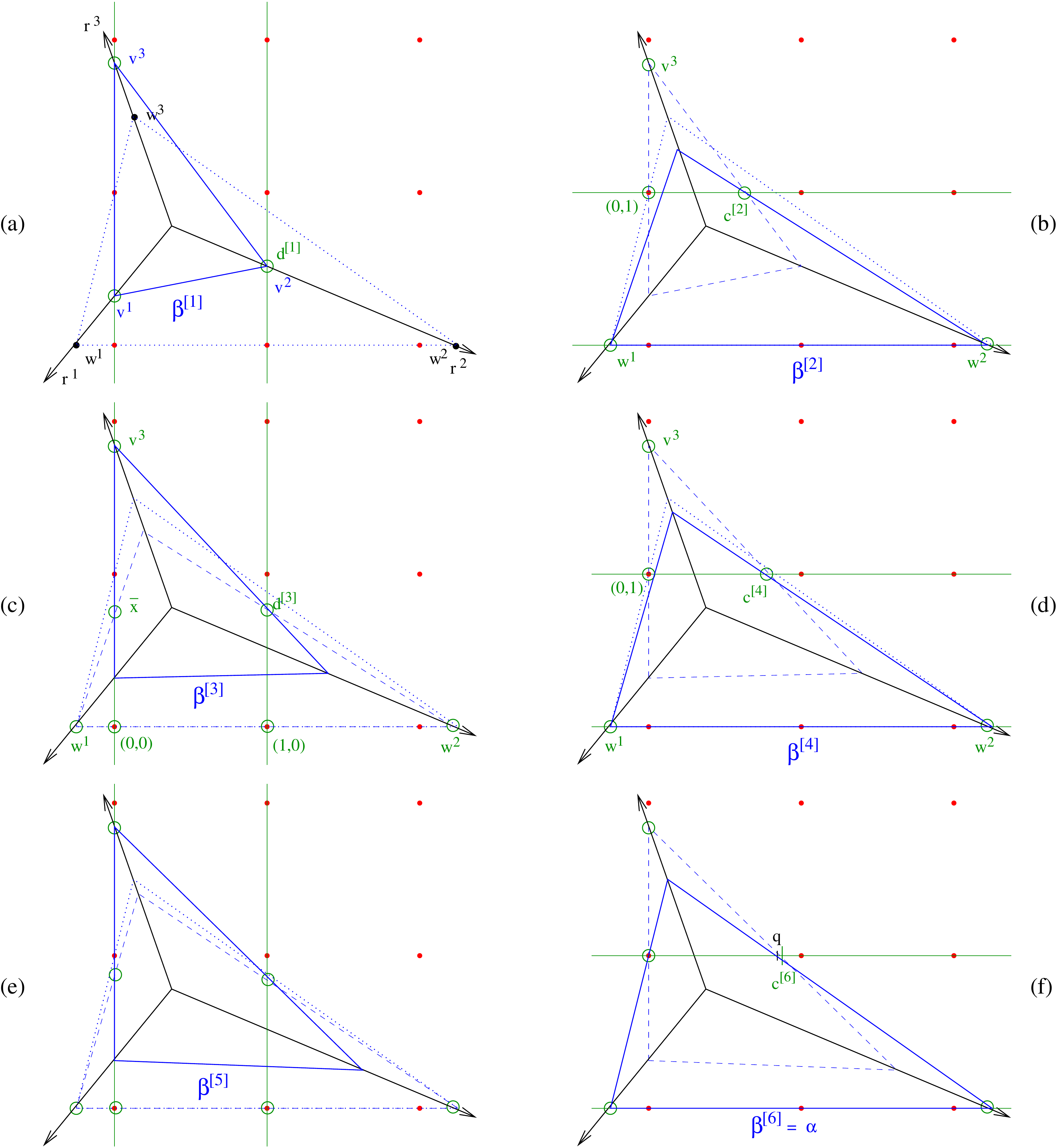}
\end{center}
\caption{In each frame, the dotted triangle is $L_{\alpha}$. The dashed triangle is the induced lattice-free set of the inequality obtained in the previous step. The circles are the $x$-vertices obtained by the application of the disjunction. The solid triangle is the induced lattice-free set of the inequality that is valid for these $x$-vertices.}
\label{lhs_movie}
\end{figure}
\end{center}

Let $\alpha$ be an inequality for $\textup{conv}(P([r^1r^2r^3],f))$ such that $L_{\alpha}$ is a standard $T^{2B}$ triangle, $0 \leq f_1 \leq 1$, and $r^3_1 < 0$. We present a sequence of split disjunctions and the rule for the selection of a valid inequality resulting from the split disjunction that eventually converges to $\alpha$.

\begin{splitproof}\label{proofalgot2Ccase1}
$\ $
\begin{enumerate}
\item \textbf{Initialization Step (Step 0)}:
Let $Q^{[0]}:= (P(R,f))^0.$ Let $Q^{ 0,\leq}:= Q^{[0]} \cap \{(x,s)\in\mathbb{R}^2\times\mathbb{R}^3\,|\,x_1 \leq 0\}$ and $Q^{ 0, \geq}:= Q^{[1]} \cap \{(x,s)\in\mathbb{R}^2\times\mathbb{R}^3\,|\,x_1 \geq 1\}$. The $x$-vertices of $Q^{ 0, \leq}$ and $Q^{0, \geq}$ are  $v^1,v^2,v^3$. Let $\beta^{[1]}:=\phi(\triangle(v^1v^2v^3))$. Define $d^{[1]}:= v^2.$

\item \textbf{Inductive step}:
At the beginning of step $j$, consider the set $Q^{[j]} = \{(x,s)\in \R^2\times \R^3_+ \mid x= f+Rs , \sum_{i=1}^3 \beta^{[j]}_i s_i \geq 1\}.$ At an even step, consider the disjunction $(x_2\leq 0) \vee (x_2 \geq 1)$ while at an odd step, consider the disjunction $(x_1 \leq 0) \vee (x_1 \geq 1).$ We now give the details of each particular step.
\begin{enumerate}
\item \textbf{Step }$\mathbf{2j-1}$: Let $Q^{2j-1, \leq}:= Q^{[2j-1]} \cap \{(x,s)\in\mathbb{R}^2\times\mathbb{R}^3\,|\,x_2 \leq 0\}$ and $Q^{2j -1, \geq}:= Q^{[2j-1]} \cap \{(x,s)\in\mathbb{R}^2\times\mathbb{R}^3\,|\,x_2 \geq 1\}$. The $x$-vertices of $\textup{conv}(Q^{2j-1, \leq} \cup Q^{2j-1, \geq})$ are $(0, 1 ),v^3,c^{[2j]},w^1,w^2$, where $c^{[2j]}$ is the intersection point of $\partial L_{\beta^{[2j+1]}}$ with the line $\{x \in \mathbb{R}^2\,|\,x_2 = 1\}$ different from $(0,1)$. At this stage either $\phi(\triangle(w^1w^2 (0, 1) ))$ or $\phi(\triangle(w^1w^2c^{[2j]}))$ is valid. \textbf{Observe that $\phi(\triangle(w^1w^2 (0, 1)))$ is the goal inequality $\alpha$. Therefore, if it is valid, we have proven that its split rank is at most $2j$}. If not, then set $\beta^{[2j]}:=\phi(\triangle(w^1w^2c^{[2j]})).$

\item \textbf{Step } $\mathbf{2j}$: Let $Q^{2j, \leq}:= Q^{[2j]} \cap \{(x,s)\in \mathbb{R}^2\times\mathbb{R}^3\,|\,x_1 \leq 0\}$ and $Q^{2j, \geq}:= Q^{[2j]} \cap \{(x,s)\in\mathbb{R}^2\times\mathbb{R}^3\,|\,x_1 \geq 1\}$. The $x$-vertices of $\textup{conv}(Q^{2j-1, \leq} \cup Q^{2j-1, \geq})$ are $w^1,w^2,v^3,(0, 1),v^2,d^{[2j+1]}$ where $d^{[2j+1]}$ is obtained as the intersection of $\partial \triangle(w^1w^2c^{[2j]})$ with $\{x \in \mathbb{R}^2\,|\,x_2=1\}$ which is different from $(1,0)$. The inequality $\phi(\triangle(v^1,v^3,d^{[2j+1]}))$ is valid. Set $\beta^{[2j+1]}:=\phi(\triangle(v^1v^3d^{[2j+1]}))$ and proceed to the next step.
\item $j \leftarrow j + 1$. \hfill $\square$
\end{enumerate}
\end{enumerate}
\end{splitproof}
See Figure \ref{lhs_movie} for an illustration of the sequence of inequalities obtained using Disjunction Sequence \ref{proofalgot2Ccase1}.
In Lemmas \ref{welldefinedvertexodd} to \ref{welldefinedineqeven}, we prove that the different steps in Disjunction Sequence \ref{proofalgot2Ccase1} are well-defined, i.e., the proposed points are $x$-vertices and the proposed inequalities are indeed valid. For the sake of clarity we repeat the definition of $c^{[2j]}$ and $d^{[2j + 1]}$ next.
\begin{notation}
Define $d^{[1]}:= v^2$. For $j = 1, 2, ...$
\begin{itemize}
\item $c^{[2j]}\in \R^2$ is the intersection point of $L_{\beta^{[2j-1]}} := \triangle(v^1v^3d^{[2j -1]})$ and the line $\{x\in \R^2 \mid x_2 = 1\}$ which is different from $(0,1)$.
\item $d^{[2j+1]}\in \R^2$ is the intersection point of $L_{\beta^{[2j]}} := \triangle(w^1w^2c^{[2j]})$ and the line $\{x\in \R^2 \mid x_1 =1\}$ which is different from $(1,0)$.
\end{itemize}
\end{notation}

\begin{lemma}[$x$-vertices for step $\mathbf{2j-1,\; r^3_1 < 0}$]\label{welldefinedvertexodd}
The $x$-vertices of $Q^{[2j-1],\leq}:=Q^{[2j-1]}\cap \{(x, s)\in \mathbb{R}^2\times\mathbb{R}^3\,|\,x_2\leq 0\}$ are $w^1$ and $w^2$. The $x$-vertices of $Q^{[2j-1],\geq}:=Q^{[2j-1]}\cap \{(x, s)\in \mathbb{R}^2\times\mathbb{R}^3\,|\,x_2\geq 1\}$ are $( 0,  1 )$, $v^3$, and $c^{[2j]}$. \label{vertices2j-1}
\end{lemma}
\textbf{Proof:} See for example frames (b) and (d) in Figure \ref{lhs_movie}.

By construction, the $x$-vertices of $Q^{[2j-1]}$, namely $v^1,v^3$ and $d^{[2j-1]}$ do not satisfy $x_2\leq 0$. Therefore the $x$-vertices of $Q^{[2j-1],\leq}$ are at the intersection of the rays $\{x\in \mathbb{R}^2\,|\, x = f + \lambda r^1, \lambda \geq 0\}$ and $\{x\in \mathbb{R}^2\,|\, x = f + \lambda r^2, \lambda \geq 0\}$ with the line $\{x\in \mathbb{R}^2\,|\, x_2 = 0\}$. We conclude that the $x$-vertices of $Q^{[2j-1],\leq} $ are $w^1$ and $w^2$.

By construction $v^1$ and $d^{[2j-1]}$ do not satisfy $x_2\geq 1$ whereas $v^3$ does. Therefore $v^3$ is a $x$-vertex of $Q^{[2j-1],\geq}$. All remaining $x$-vertices are at the intersection of $\partial L_{ \beta^{[2j-1]}}$ with $\{x\in \mathbb{R}^2\,|\,x_2=1\}$ i.e. $(0, 1)$ and $c^{[2j]}$.\cqfd

\begin{lemma}[Finding $\mathbf{\beta^{[2j]},\; r^3_1 < 0}$]\label{welldefinedineqodd}
At stage $2j-1$, at least one of the following holds
\begin{enumerate}[(i)]
\item $\triangle(w^1w^2c^{[2j]}) \subseteq \triangle(w^1w^2(0, 1))$ and
$\beta^{[2j]}:=\phi(\triangle(w^1w^2c^{[2j]})) $ is valid for $\conv(Q^{[2j-1],\leq} \cup Q^{[2j-1],\geq}).$
\item $\triangle(w^1w^2c^{[2j]}) \supseteq \triangle(w^1w^2(0, 1))$ and
$\beta^{[2j]}:=\phi(\triangle(w^1w^2(0, 1))) $ is valid for $\conv(Q^{[2j-1],\leq} \cup Q^{[2j-1],\geq}).$
\end{enumerate}
\end{lemma}
\textbf{Proof:} See Figure \ref{lhs_movie}, frames (b) and (d) for example of case (i), and frame (f) for example of case (ii).

The triangles $\triangle(w^1w^2c^{[2j]})$ and $\triangle(w^1w^2 (0, 1))$ share the side $w^1w^2$. By definition, the third vertex of both these triangles lie on the ray $\{x\in \mathbb{R}^2\,|\,  x= f+\lambda r^3, \lambda \geq 0\}$. It follows that, if we compare the two vertices, one of them must be closer to $f$. The triangle for which the third vertex is closer to $f$ is therefore included in the other.

It remains to verify that the inequality corresponding to the included triangle is valid for $\conv(Q^{[2j-1],\leq} \cup Q^{[2j-1],\geq}).$ To do that, it suffices to check that all $x$-vertices given by Lemma \ref{vertices2j-1} do not lie in the interior of the triangle. Observe that it is sufficient to check only whether $v^3$ does not lie in the interior of the smaller triangle as the other $x$-vertices cannot lie in the interior of $\triangle(w^1w^2c^{[2j]})$ and $\triangle(w^1w^2 (0, 1))$.

Claim: $v^3 \not\in \triangle(w^1w^2(0, 1)).$ By hypothesis, $w^1_1<0$, $v^3_1=0$, and $v^3_2>1.$ It follows that $v^3$ lies above the line $w^1(0, 1)$. As the third vertex, (i.e. $w^2$) lies below the same line, it follows that $v^3 \not\in \triangle(w^1w^2(0, 1)).$ This implies that if $\triangle(w^1w^2(0, 1)) \subseteq \triangle(w^1w^2c^{[2j]})$, then no $x$-vertex lies in the interior of this triangle and the corresponding inequality is valid. On the other hand, if $\triangle(w^1w^2c^{[2j]}) \subseteq \triangle(w^1w^2(0, 1))$, as $v^3\not\in \triangle(w^1w^2(0, 1))$, then it also follows that the corresponding inequality is valid.\cqfd

\begin{lemma}[$x$-vertices for step $\mathbf{2j,\;r^3_1 < 0}$]\label{welldefinedvertexeven}
The $x$-vertices of $Q^{[2j],\leq}:=Q^{[2j]} \cap \{(x,s)\in \mathbb{R}^2\times\mathbb{R}^3\,|\,x_1\leq 0\}$ are $w^1$, $(0, 0)$, $v^3$, and a point $\bar{x}$ which satisfies $\bar{x}_1=0$ and $0< \bar{x}_2<1$. The $x$-vertices of $Q^{[2j],\geq}:=Q^{[2j]} \cap \{(x,s)\in \mathbb{R}^2\times\mathbb{R}^3\,|\,x_1\geq 1\}$ are $d^{[2j+1]}$, $w^2$, and $(1, 0)$.
\label{vertices2j}
\end{lemma}
\textbf{Proof:} See for example frames (c) and (e) of Figure \ref{lhs_movie}.

We start by computing all the $x$-vertices of $Q^{[2j],\leq}$. Since $w^1$ is a $x$-vertex of $Q^{[2j]}$ and is valid for $\{x\in \mathbb{R}^2\,|\,x_1 \leq 0\}$, it is an $x$-vertex of $Q^{[2j],\leq}$. The other $x$-vertices come from the intersection of either $\partial L_{ \beta^{[2j]}}$ or the ray $\{x \in \mathbb{R}^2\,|\, x= f + \lambda_3r^3, \lambda_3 \geq 0\}$ with the line $\{x \in \mathbb{R}^2\,|\,x_1=0\}$. In the first case, we obtain $(0, 0)$ and $\bar{x}$. In the second case we obtain $v^3$ as an $x$-vertex.

Consider now $Q^{[2j],\geq}.$ Similar to $w^1$ in the previous case, $w^2$ is an $x$-vertex as it is satisfies $\{x\in \mathbb{R}^2\,|\,x_1\geq 1\}$. The other $x$-vertices come from the intersection of $\partial L_{ \beta^{[2j]}}$ with $\{x\in \R^2\mid x_1=1\}$. We obtain therefore $d^{[2j+1]}$ and $(1, 0).$ \cqfd

\begin{lemma}[Finding $\mathbf{\beta^{[2j+1]},\; r^3_1 < 0}$]\label{welldefinedineqeven}
The inequality $\beta^{[2j+1]}:=\phi(\triangle(v^3(0, 1)d^{[2j+1]}))$ is valid for $\conv(Q^{[2j],\leq} \cup Q^{[2j],\geq})$.
\end{lemma}
\textbf{Proof:} See for example frames (c) and (e) of Figure \ref{lhs_movie}.

It suffices to check that all $x$-vertices given in Lemma \ref{vertices2j} do not lie in the interior of $\triangle(v^3(0, 1)d^{[2j+1]})$. Observe that $((0,0), \mathcal{M}(0,0))$ and $((1,0), \mathcal{M}(1,0))$ can be obtained as convex combination of $(w^1, \mathcal{M}(w^1))$ and $(w^2, \mathcal{M}(w^2))$. Therefore it is not necessary to verify that $(0,0)$ and $(1,0)$ do not lie in the interior of $\triangle(v^3(0, 1)d^{[2j+1]})$. Also it is easily observed that the $x$-vertices $\bar{x}$ and $v^3$ do not lie in interior of $\triangle(v^3(0, 1)d^{[2j+1]})$.

Claim: $w^1\not\in \triangle(v^3(0, 1)d^{[2j+1]})$. The side $(v^3 (0, 1))$ of the triangle satisfies $x_1=0$ and the third vertex $z$ lies on the set $\{x\in \mathbb{R}^2\,|\, x = f + \lambda r^2, \lambda \geq 0 \}$, with $z_1 > f_1$. As $f_1>0$ and $r^2_1 >0$, we conclude that all points in the triangle $\triangle(v^3(0, 1)d^{[2j+1]})$ have a non-negative first coordinate. As $w^1_1 <0$, we obtain $w^1\not\in \triangle(v^3(0, 1)d^{[2j+1]})$.

Claim: $w^2\not\in \triangle(v^3(0, 1)d^{[2j+1]})$. From the previous step, we know that one side of $\triangle(w^1w^2c^{[2j]})$ passes through the points $w^2$ and $d^{[2j+1]}$ and one point $v$ of the form $v=f+\lambda_3 r^3, \lambda_3\geq 0$ such that $\lambda(v) > \lambda(v^3)$. We conclude that the side of $\triangle(v^3(0, 1)d^{[2j+1]})$ that links $v^3$ to $d^{[2j+1]}$ intersects
$\{x\in \R^2\mid x=f+\lambda r^2, \lambda\geq 0\}$ at a point $z$ such that $z=f+ \lambda_2 r^2$, with $\lambda(z) < \lambda(w^2)$. In particular this implies that $z_2 >0$ and that the side of  $\triangle(v^3(0, 1)d^{[2j+1]})$ linking a point on $\{x\in \R^2\mid x=f+\lambda r^1,\lambda\geq 0\}$ to $z$ is always above the axis $x_2=0$. Therefore $w^2 \not\in  \triangle(v^3(0, 1)d^{[2j+1]})$. \cqfd

We now want to understand the convergence of this procedure. To do that, we denote $c^{[0]}=:c$. The next Proposition indicates how $c^{[2j]}$ can be computed recursively.

\begin{lemma}
Denote $w^2 =: (e, 0)$ and $v^3 =: (0, a)$. The previous sequence of split disjunctions
provides the recursion
\begin{align*}
d^{[2j+1]} & = \left( 1, \frac{e-1}{e-c^{[2j]}_1}\right)\\
c^{[2j+2]} & = \left( \frac{a-1}{a-d^{[2j+1]}_2}, 1\right).
\end{align*}
\end{lemma}
\textbf{Proof:}
The point $d^{[2j+1]}$ is obtained by intersecting the line $\{x \in \mathbb{R}^2\,|\,x_1=1\}$ with the line joining $c^{[2j]}$ with $w^2$. The point $c^{[2j+2]}$ is obtained by intersecting the line $\{x \in \mathbb{R}^2\,|\,x_2=1\}$ with the line joining $d^{[2j+1]}$ with $v^3$. By computing the different equations of the lines, we obtain the desired result.\cqfd

Finally, we show that the Disjunction Sequence \ref{proofalgot2Ccase1} converges to the goal inequality $\alpha$. Remember $q$ is the point of the form $f + \lambda_2r^2 + \lambda_3r^3$, $\lambda_2, \lambda_3 >0$ such that $q_2= 1$ and $L_{\alpha} = \triangle(w^1w^2(0,1)) = \triangle(w^1w^2q)$. At the end of step \textbf{$2j -1$} in Disjunction Sequence \ref{proofalgot2Ccase1} if $e^{[2j]}_1 < q_1$, then we select the inequality corresponding to the triangle $\triangle(w^1w^1c^{[2j]})$. If $e^{[2j]}_1 \geq q_1$, then $\beta^{[2j]} := \alpha$. Hence to prove the result we will prove the following: $\textup{lim}_{i \rightarrow \infty}c^{[2j]}_1 > q_1$ in Lemma \ref{seq1}. This shows that in a finite number of iterations, $j^{*}$, the inequality corresponding to $\triangle(w^1w^2(0,1))$ will be valid at the end of step $2j^{*}$ thus completing the proof.

To simplify notation let $c^{[2j]}_1 =: c^j$ and $d^{[2j + 1]}_2 =: d^j$.
\begin{lemma}[Finite Convergence; $r^3_1 < 0$]\label{seq1}
$\textup{lim}_{ i \rightarrow \infty}c^{i} = C$ where $C = \textup{min}\left\{1, \left(1 - \frac{1}{a}\right)e \right\} > q_1$.
\end{lemma}
\textbf{Proof:} Since $v^3_2 >1 $, we obtain $a >1$. Since $w^2_1 > 1$, we obtain $ e > 1$.

Claim: $\textup{min}\left\{1, \left(1 - \frac{1}{a}\right)e \right\} > q_1$. By assumption $0< w^3_1 < 1$. Therefore $\lambda(v^3) > \lambda(w^3)$. This implies that $\left(1 - \frac{1}{a}\right)e  > q_1$. Also by definition of $q$, we obtain $q_1 < 1$.

We assume $0 <c^0 < 1$ and $c^0 < \left(1 - \frac{1}{a}\right)e$ (otherwise there remains nothing to prove).
We now prove this result in the following steps:
\begin{itemize}
\item $C \geq \textup{sup}\{c^i\}$: Note that $c^{0} < 1$ and $c^{0} < \left(1 - \frac{1}{a}\right)e$. Assume by the induction that $c^{k} \leq 1$ and $c^{k} \leq \left(1 - \frac{1}{a}\right)e$ $\forall k \in \{0, ...,n\}$.
    \begin{enumerate}
    \item Now observe that $d^{n} = \frac{e - 1}{e - c^{n}} \leq 1$ since $c^{n} \leq 1$. Therefore,
    \begin{eqnarray}
    c^{n +1} = \frac{a - 1}{ a- d^{n}} \leq 1.
    \end{eqnarray}
    \item Observe that
    \begin{eqnarray}
    c^n &\leq &\left (1 - \frac{1}{a}\right)e \nonumber \\
    e - c^{n} &\geq& \frac{e}{a} \nonumber \\
    \frac{1}{e - c^{n}} &\leq& \frac{a}{e} \nonumber \\
    \frac{ e- 1}{e - c^{n}} &\leq & \left( 1 - \frac{1}{e} \right)a \quad [\textup{since } e > 1 \geq c^n, a > 0] \nonumber \\
    d^{n} &\leq & \left( 1 - \frac{1}{e} \right)a \nonumber \\
    a - d^n & \geq & \frac{a}{e} \nonumber \\
    c^{n+1} = \frac{a -1}{a - d^n} &\leq& \left(1 - \frac{1}{a}\right)e \quad [\textup{since } a > 1 \geq d^n, e > 0].
    \end{eqnarray}
    \end{enumerate}
\item $\{c^{i}\}$ is a non-decreasing sequence:  By algebraic manipulations we obtain,
\begin{eqnarray}
c^{i +1} - c^{i} = \frac{ (1 - c^{i})(a(e - c^i) - e)}{a(e - c^{i}) - (e - 1)}.
\end{eqnarray}
By previous part, $c^{i} \leq 1$ and $c^{i} \leq \left(1 - \frac{1}{a}\right)e$ $\forall i$. Therefore, $a(e - c^{i}) - (e - 1) \geq 0$ and $(1 - c^{i})(a(e - c^i) - e) \geq 0$.
\item If $F = \textup{sup}\{c^i\}$, then $F \geq  C$: Assume by contradiction that $F = \textup{sup}\{c^i\}$ and $F < C$. By definition of $F$, $\forall \epsilon > 0$, there exists $i_{\epsilon}$ such that $c^{i_ \epsilon} \geq F - \epsilon$. Let $\delta^1 = 1 - F > 0$ and $\delta^2 = e -F$. Note that since $F < C$, we have $a\delta^2 > e$. There are two cases:
    \begin{enumerate}
    \item $a\delta^1 < 1$: Then choose any $i$ such that $c^{i} > F - \frac{\delta^1(a \delta^2 - e)}{1 - a\delta^1}$. By assumption $c^i \leq F$. Let $F - c^{i} = \eta$. Then we obtain that
        \begin{eqnarray}
        \frac{\delta^1(a \delta^2 - e)}{1 - a\delta^1} &>& F - c^{i} = \eta \nonumber \\
        \label{start}\Rightarrow \delta^1(a \delta^2 - e) &>& \eta(1 - a \delta^1) \\
        \Rightarrow \delta^1(a \delta^2 + a\eta - e) &>& \eta \nonumber \\
        \Rightarrow \delta^1(a \delta^2 + a\eta - e) + \eta(a\delta^2 + a\eta - e) &>&  \eta(a\delta^2 + a\eta - e + 1) \nonumber
        \end{eqnarray}
        Now note that $a\delta^2 + a\eta - e + 1 > a\delta^2  - e > 0$. Therefore,
        \begin{eqnarray}\label{gotcon}
        \frac{(\delta^1 + \eta)(a \delta^2 + a\eta - e)}{a\delta^2 + a\eta - e + 1} > \eta
        \end{eqnarray}
        Now note that $\frac{(\delta^1 + \eta)(a \delta^2 + a\eta - e)}{a\delta^2 + a\eta - e + 1}$ $= \frac{(1 - F + \eta)(a(e - F + \eta) - e)}{a(e - F + \eta) - (e -1)}$ $ = \frac{(1 - c^i)(a(e - c^i) - e)}{a(e - c^i) - (e -1)} = c^{i+1} - c^i$. Therefore, (\ref{gotcon}) implies that $c^{i+1} - c^{i} > \eta$ or $c^{i+1} > F$, a contradiction.
    \item  $a\delta^1 \geq 1$: Then choose any $c^{i}$. Let $\eta = F - c^{i} \geq 0$ by assumption. As $\delta^1(a \delta^2 - e)  > 0$ and $a\delta^1 \geq 1$, we obtain that $ \delta^1(a \delta^2 - e) > \eta(1 - a \delta^1)$ which is the same as (\ref{start}). Thus again we obtain that $c^{i+1} - c^{i} > \eta$, a contradiction.\hfill $\square$
    \end{enumerate}
\end{itemize}

\subsubsection{Case 2: $f_2 \leq 1$ and $r^3_1<0$.}\label{t2cc2}
If $0 \leq f_1 \leq 1$, then this case is covered in Section \ref{t2cc1}.

Since $r^3_1 < 0$, the remaining case is $f_1 > 1$. The first disjunction $(x_1 \leq 0)\vee(x_1 \geq 1)$ in the Disjunction Sequence \ref{proofalgot2Ccase1} does not yield any new inequalities. In this case we alter Disjunction Sequence \ref{proofalgot2Ccase1} by first starting with the disjunction $(x_2 \leq 0)$ $\vee$ $(x_2 \geq 1)$. The initialization stage is updated as follows: Consider $Q^{0, \leq}:= (P(R,f))^0 \cap \{(x,s)\in\mathbb{R}^2\times\mathbb{R}^3\,|\,x_2 \leq 0 \}$ and $Q^{0 , \geq }:= (P(R,f))^0 \cap \{(x,s)\in\mathbb{R}^2\times\mathbb{R}^3\,|\,x_2 \geq 0 \}$. The $x$-vertices of $Q^{0, \leq}$ are $w^1$ and $w^2$ and the $x$-vertex of $Q^{0,\geq}$ is $\bar{w}^3$, where $\bar{w}^3$ is the intersection point of the ray $\{f + \lambda_3r^3\,|\,\lambda_3 \geq 0\}$ and the line $\{x \in \mathbb{R}^2\,|\, x_2 = 1\}$. The inequality $\phi(\triangle(w^1w^2\bar{w}^3)) =: \beta^{[1]}$ is valid for $\textup{conv}(Q^{0, \leq}\cup Q^{0, \geq})$. The rest of the algorithm is identical to Disjunction Sequence \ref{proofalgot2Ccase1} except that odd steps are now even and the even steps are now odd.

\subsubsection{Case 3: $r^3_1= 0$.}\label{t2cc1a}

In this case, it is easily verified that the split rank is exactly two. We start by considering the disjunction $(x_1 \leq 0)\vee (x_1 \geq 1)$. The inequality $\beta^{[1]}$ that is valid for $\textup{conv}((Q^{0,\leq} \cup Q^{0, \geq}))$ has the induced lattice-free set $\textup{conv}(v^1, v^2) + \textup{cone}(0,1)$. Then considering the disjunction $(x_2 \leq 0)\vee (x_2 \geq 1)$, we obtain the goal inequality. The proof is very similar to the proofs in Section \ref{t2cc1}.

\subsubsection{Case 4: $r^3_1 > 0$.}\label{t2cc3}
We now consider the case where $r^3_1 > 0$. As discussed in the outline of the proof of Proposition \ref{theoremt2c}, the idea of the procedure is essentially the same as in Case 1. We apply the sequence of disjunctions $(x_1\leq 0) \vee (x_1 \geq 1)$ and $(x_2\leq 0) \vee (x_2 \geq 1)$. At each step, we replace all previous inequalities by one valid inequality obtained after the disjunction and proceed. We will prove that after a finite number of steps, this procedure converges to the desired inequality.

The primary difference in this case is that the initialization step is different, where a different rank 2 inequality is added. Moreover the proof of convergence is more involved than the previous cases. We note here that Disjunction Sequence \ref{proofalgot2Ccase1} can be applied to this case. However, we are unable to proof that Disjunction Sequence \ref{proofalgot2Ccase1} converges in finite time in the case where $r^3_1 >0$. On the other hand, it appears that Disjunction Sequence \ref{proofalgot2Ccase3} that we present next, does not seem to apply for the case where $r^3_1 <0$.

\begin{center}
\begin{figure}[!h]
\begin{center}
\includegraphics[height=19cm]{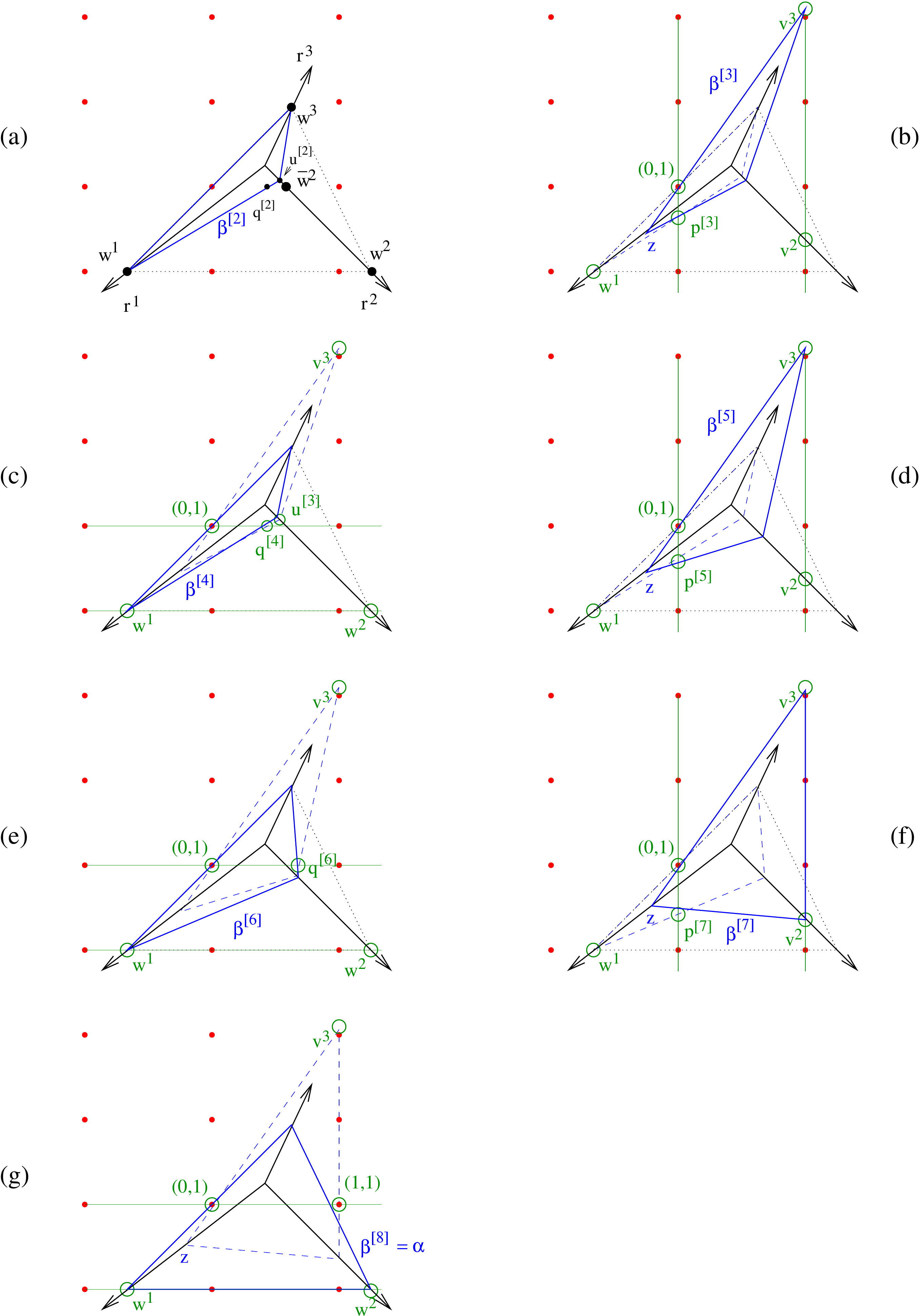}
\end{center}
\caption{In each frame, the dotted triangle is $L_{\alpha}$. The dashed triangle is the induced lattice-free set of the inequality obtained in the previous step. The circles are the $x$-vertices obtained by the application of the disjunction. The solid triangle is the induced lattice-free set of the inequality that is valid for these $x$-vertices.}
\label{rhs_movie}
\end{figure}
\end{center}

As before, let $\alpha$ be the goal inequality such that $L_{\alpha}$ is triangle of type $T^{2B}$.

\begin{splitproof}\label{proofalgot2Ccase3}
$\ $
\begin{enumerate}
\item \textbf{Initialization step}: First consider the two-variable problem $P((r^1,r^3),f)$. By definition the triangle $C:= fw^1w^3$ does not contain any integer point in its interior. Therefore $\phi(C)$ is a valid inequality for $\textup{conv}(P((r^1,r^3),f))$. By Proposition \ref{generallifting}, there also exists $\epsilon>0$ such that, denoting $u^{[2]} :=f+\epsilon r^2$, we obtain that $\beta^{[2]}:=\phi(\triangle(w^1w^3u^{[2]}))$ is a valid inequality for $P(R,f)$. By Proposition \ref{generallifting} and Proposition \ref{2varthm}, we also know that this inequality has a split rank at most two. Let $q^{[2]}$ be the intersection point of $(\triangle(w^1w^3u^{[2]}))$ with the line $\{x \in \mathbb{R}^2\,|\, x_2 = 1\}$. We then directly proceed to step 2 in the inductive process.

\item \textbf{Inductive step}: At the beginning of iteration $j$, we consider the set $Q^{[j]} = \{(x,s)\in \R^2\times \R^3_+ \mid x= f+Rs , \sum_{i=1}^3 \beta^{[j]}_i s_i \geq 1\}.$

\begin{enumerate}
\item \textbf{Step }$\mathbf{2j}$: We consider $Q^{[2j]}$ where $\beta^{[2j]}:= \phi(\triangle(w^1w^3q^{[2j]}))$. The $x$-vertices of $Q^{[2j],\leq}:= Q^{[2j]}\cap \{(x,s) \in \R^2 \times \mathbb{R}^3\mid x_1\leq 0\}$ are $w^1$, $(0,1)$, and $p^{[2j +1]}$ where $p^{[2j+1]}$ is obtained as an intersection of $\triangle(w^1w^3q^{[2j]})$ with $\{x \in \mathbb{R}^2\,|\,x_1=0\}$. The $x$-vertices of $Q^{[2j],\geq}:= Q^{[2j]}\cap \{(x,s) \in \R^2 \times \mathbb{R}^3\mid x_1\geq 1\}$ are $v^3$ and either $v^2$ or two points $z^1,z^2$ which satisfy $z^1_1=z^2_1=1$ and $0<z^1_2<z^2_2<1$. Define $z$ to be the intersection point of the line $v^3 (0, 1)$ with the ray $\{x \in \mathbb{R}^2\,|\, x = f+\lambda r^1, \lambda\geq 0\}$. We now distinguish between two cases.
	 \begin{enumerate}
            \item $p^{[2j+1]}$ is above the line $zv^2$: The inequality $\beta^{[2j+1]}:= \phi(\triangle((0, 1)v^3p^{[2j+1]}))$ is valid for $\textup{conv}(Q^{[2j],\leq}\cup Q^{[2j],\geq})$.
            \item $p^{[2j+1]}$ is below or on the line $zv^2$: The inequality $\beta^{[2j+1]}:= \phi(\triangle(v^2v^3(0 \ 1))$ is valid for $\textup{conv}(Q^{[2j],\leq}\cup Q^{[2j],\geq})$. \textbf{Go to the termination step}.
        \end{enumerate}
    \item \textbf{Step }$\mathbf{2j+1}$: We consider $Q^{[2j+1]}$ where $\beta^{[2j+1]}:= \phi(\triangle((0, 1)v^3p^{[2j+1]}))$. Let $u^{[2j+1]}$ be the vertex of $\triangle((0, 1)v^3p^{[2j+1]})$ that lies on the ray $f+\lambda r^2, \lambda \geq 0$. The $x$-vertices of $Q^{[2j+1],\leq}:=Q^{[2j+1]}\cap \{x\in \R^2\mid x_2\leq 0\}$ are $w^1$ and $w^2$. The $x$-vertices of $Q^{[2j+1],\geq}:=Q^{[2j+1]}\cap \{x\in \R^2\mid x_2\geq 1\}$ are
    \begin{enumerate}
        \item $u^{[2j+1]}_1 \leq \bar w^2_1$: $v^3,(0, 1),\bar w^2,u^{[2j+1]},q^{[2j+2]}$.
        \item Either $u^{[2j+1]}_1 > w^2_1$ or $f_2<1$ (in which case $\bar w^2$ does not exist): $v^3,(0, 1),q^{[2j+2]}$.
    \end{enumerate}

    In both cases, $q^{[2j+2]}$ is obtained as the intersection point of $L_{\beta^{[2j+1]}}$ with $\{x\in \mathbb{R}^2\,|\,x_2 = 1\}$. Remember $q$ is the intersection point of $L_{\alpha}$ with the line $\{x\in \mathbb{R}^2\,|\,x_2=1\}$ (Notation \ref{allcommonnotation}). Now two cases occur. \textbf{Either $q^{[2j+2]}_1 \geq q_1$, in which case, the goal inequality $\alpha$ is valid for $\textup{conv}(Q^{[2j +1], \leq}\cup Q^{[2j +1], \geq})$ and thus proven to be of split rank at most $2j+2$}, or $q^{[2j+2]}_1 < q_1$, in which case the inequality $\phi(\triangle(w^1 (0, 1) q^{[2j+2]})) = \phi(\triangle(w^1 w^3 q^{[2j+2]})) =:\beta^{[2j+2]}$ is valid for $\textup{conv}(Q^{[2j +1], \leq}\cup Q^{[2j +1], \geq})$.
    \item $j \leftarrow j + 1$.
    \end{enumerate}
    \item \textbf{Termination Step}: We consider $Q^{[2j+1]}$ where $\beta^{[2j+1]}=\phi(\triangle(v^2v^3(0, 1)))$. Let $Q^{[2j+1],\leq}:=Q^{[2j+1]}\cap \{(x,s) \in \R^2 \times \mathbb{R}^3\mid x_2\leq 0\}$ and $Q^{[2j+1],\geq}:=Q^{[2j+1]}\cap \{(x,s) \in \R^2 \times \mathbb{R}^3\mid x_2\geq 1\}$. Then the goal inequality $\alpha$ is valid for $\conv(Q^{[2j+1],\leq}\cup Q^{[2j+1],\geq})$. \hfill $\square$
\end{enumerate}
\end{splitproof}
See Figure \ref{rhs_movie} for an illustration of the sequence of inequalities obtained using Disjunction Sequence \ref{proofalgot2Ccase3}.

\noindent In Lemmas \ref{vertices2jRHS} to  \ref{ddsrhs}, we prove that the different steps of the Disjunction sequence \ref{proofalgot2Ccase3} are well-defined i.e., the proposed points are $x$-vertices and the proposed inequalities are indeed valid. For the sake of clarity we repeat the definitions of some of the points introduced in Disjunction sequence \ref{proofalgot2Ccase3} and earlier.
\begin{notation}
$\quad$
\begin{enumerate}
\item Let $j$ be fixed.
\begin{itemize}
\item $p^{[2j+1]}\neq (0,1) $ is the intersection point of $L_{\beta^{[2j]}}:= \triangle(w^1w^3q^{[2j]})$ and the
line $\{x\in \R^2 \mid x_1 = 0\}$,
\item $q^{[2j + 2]}\neq (0,1)$ is the intersection point of $L_{\beta^{[2j+1]}}:= \triangle((0,1)v^3p^{[2j + 1]})$ and the line $\{x\in \R^2 \mid x_2=1\}$,
\item $u^{[2j+1]}$ is the intersection point of $L_{\beta^{[2j+1]}}:= \triangle((0,1)v^3p^{[2j + 1]})$ with the ray $\{x\in \R^2 \mid x= f+ \lambda_2 r^2, \lambda_2\geq 0\}.$
\end{itemize}
\item $z$ is the intersection point of the line $v^3 (0, 1)$ with the ray $\{x\in \R^2\,|\,x = f+\lambda r^1, \lambda\geq 0\}$.
\item $q$ is the intersection point of side $w^2w^3$ of $\partial L_{\alpha}$ with the line $\{x\in \mathbb{R}^2\,|\,x_2=1\}$.
\end{enumerate}
\end{notation}

\begin{lemma}[$x$-vertices for step $\mathbf{2j,\; r^3_1>0}$]
The $x$-vertices of $Q^{[2j],\leq}:=Q^{[2j]} \cap \{(x, s)\in \R^2\times \mathbb{R}^3\mid x_1\leq 0\}$ are $w^1, (0, 1)$, and $p^{[2j+1]}$. The $x$-vertices of $Q^{[2j],\geq}:=Q^{[2j]} \cap \{(x, s)\in \R^2\times \mathbb{R}^3\mid x_1 \geq 1\}$ are $v^3$ and either $v^2$ or two points $z^1,z^2$ which satisfy $z^1_1=z^2_1=1$ and $0<z^1_2<z^2_2<1.$
\label{vertices2jRHS}
\end{lemma}
\textbf{Proof:} See for example frames (b), (d), and (f) of Figure \ref{rhs_movie}.

Consider the $x$-vertices of $Q^{[2j],\leq}$. Clearly $w^1$ is a $x$-vertex as it satisfies $w^1_1<0$. Note that $v^1$ is the intersection point of $\{x\in \mathbb{R}^2\,|\,x = f + \lambda_1 r^1, \lambda_1 \geq 0\}$ and the line $\{x \in \mathbb{R}^2\,|\, x_1 = 0\}$. However, $v^1$ is not an $x$-vertex of $Q^{[2j],\leq}$ as $v^1 \in \textup{int}(\triangle(w^1w^3q^{[2j]}))$. Finally, the intersection points of $L_{\beta^{[2j]}}$ with the line $\{x\in \mathbb{R}^2\,|\,x_1 = 0\}$, namely $(0, 1)$ and $p^{[2j+1]}$, are $x$-vertices.

Consider now $Q^{[2j],\geq}$. Observe that $w^3$ is not a $x$-vertex since $w^3_1<1$. The point $v^3$ is a $x$-vertex. The other $x$-vertices can be obtained as the intersection point(s) of $L_{\beta^{[2j]}}$ with the line $\{x \in \mathbb{R}^2\,|\,x_1 = 1\}$. Therefore $v^2$ is a $x$-vertex iff $v^2 \notin \textup{int}(L_{\beta^{[2j]}})$. If $v^2 \in \textup{int}(L_{\beta^{[2j]}})$, then we obtain two $x$-vertices $z^1$ and $z^2$ which satisfy $z^1_1=z^2_1=1$. Observe that by hypothesis $0<v^2_2<1$ and we must also have $0<z^1_2<z^2_2<1$. \cqfd

\begin{lemma}[Finding $\mathbf{\beta^{[2j+1]}, r^3_1 > 0}$]
If
\begin{enumerate}[(i)]
\item
$p^{[2j+1]}$ lies above the line $zv^2$, then $\phi(\triangle((0, 1) v^3 p^{[2j+1]}))$ is
valid for $\conv(Q^{[2j],\leq} \cup  Q^{[2j],\geq})$,
\item
$p^{[2j+1]}$ is below or on the line $zv^2$, then $\phi(\triangle(v^2v^3(0, 1)))$ is valid
for $\conv(Q^{[2j],\leq} \cup  Q^{[2j],\geq})$.
\end{enumerate}
\label{lemma2jRHS}
\end{lemma}
\textbf{Proof:} See Figure \ref{rhs_movie}, frames (b), (d) for example of case (i) and frame (f) for example of case (ii).
\begin{enumerate}[(i)]
\item We have to check that all the $x$-vertices given by Lemma \ref{vertices2jRHS} do not lie in
the interior of $\triangle(v^3(0, 1)p^{[2j+1]})$. Clearly $(0,1)$, $p^{[2j+1]}$, and $v^3$ lie on the boundary of $\triangle(v^3(0, 1)p^{[2j+1]})$. Therefore, we have to verify that $w^1$, $v^2$ (or $z^1$ and $z^2$) do not lie in the interior of $\triangle(v^3(0, 1)p^{[2j+1]})$.
\begin{itemize}
\item $\mathbf{w^1}$: The line $w^1 (0, 1)$ meets the ray $\{x \in \mathbb{R}^2\,|\,x = f+\lambda_3r^3,\lambda_3\geq 0\}$ at $w^3$ and the ray $\{x=f+\lambda_1r^1\,|\,\lambda_1\geq 0\}$ at $w^1$ and passes through $(0, 1)$. On the other hand, the line $v^3(0, 1)$ meets the ray $\{x\in \mathbb{R}^2\,|\,x= f+\lambda_3r^3, \lambda_3\geq 0\}$ at $v^3$ and the ray $\{x= f+\lambda_1r^1\,|\, \lambda_1\geq 0\}$ at $z$ and passes through $(0, 1)$. As $\lambda(w^3) < \lambda(v^3)$, we conclude that we must have $\lambda(z) < \lambda(w^1)$. As $z$ is a vertex of  $\triangle(v^3(0, 1)p^{[2j+1]})$, we conclude that $w^1 \not\in  \triangle(v^3(0, 1)p^{[2j+1]})$.
\item $v^2$ (or $z^1$ and $z^2$): We verify that if $x$ such that $x_1=1$, then $x$ does not lie in the interior of $\triangle(v^3(0, 1)p^{[2j+1]})$. Let $y$ be the intersection point $zp^{[2j+1]}$ with the ray $\{x = f+\lambda_2r^2\,|\,\lambda_2\geq 0\}$. As $p^{[2j+1]}$ is above the line $zv^2$, we conclude that $\lambda(y) < \lambda(v^2)$. Therefore $y_1<1.$ Hence $\triangle(v^3(0, 1) p^{[2j+1]}) \cap \{x\in \mathbb{R}^2\,|\,x_1=1\} = \{v^3\}.$ Therefore, for all $x\in \R^2$ with $x_1=1$, we have $x\not\in\inter(\triangle(v^3(0, 1)p^{[2j+1]})$.
\end{itemize}
\item From the proof of (i), we can also conclude that $w^1 \not\in\triangle(v^2v^3(0, 1)).$ Since $\triangle(v^2v^3(0, 1))\cap \{x_1=1\} \in \partial \triangle(v^2v^3(0, 1))$, we conclude that there does not exist $x\in \{f+\lambda_2r^2,\lambda_2\geq 0\}$ with $x_1=1$ and $x\in\inter\triangle(v^2v^3(0, 1)).$ It remains to prove that $p^{[2j+1]}\not\in\inter(\triangle(v^2v^3(0, 1)).$ This follows from the fact that $p^{[2j+1]}$ lies below the line $zv^2$, which is a side
of the triangle.
\cqfd
\end{enumerate}

\begin{lemma}[$x$-vertices for step $\mathbf{2j+1,\; r^3_1>0}$)]
The $x$-vertices of $Q^{[2j+1],\leq}:=Q^{[2j+1]}\cap \{(x,s)\in \mathbb{R}^2\times \mathbb{R}^3\,|\,x_2\leq 0\}$ are $w^1$ and $w^2$. The $x$-vertices of $Q^{[2j+1],\geq}:= Q^{[2j+1]}\cap \{(x,s)\in \mathbb{R}^2\times \mathbb{R}^3\,|\,x_2\geq 1\}$ are
\begin{enumerate}[(i)]
\item if $\lambda(u^{[2j+1]}) \leq \lambda(\bar w^2)$, then $v^3,(0, 1),\bar w^2, q^{[2j+2]}$, and $u^{[2j+1]}$.
\item if $\lambda(u^{[2j+1]}) > \lambda(\bar w^2)$ or $f_2 < 1$ ($\bar w^2$ does not exist), then $v^3,(0, 1)$, and $q^{[2j+2]}$.
\end{enumerate}
\label{vertices2j+1RHS}
\end{lemma}
\textbf{Proof:} See Figure \ref{rhs_movie} frame (c) for example of case (i) and frame (e) for example of case (ii).

We first consider $Q^{[2j+1],\leq}$. Observe that we assume that $p^{[2j +1]}$ is above the line $zv^2$, otherwise the Disjunctive Sequence \ref{proofalgot2Ccase3} reaches the termination step. Therefore, the only $x$-vertices are $w^1$ and $w^2$.

We now  consider $Q^{[2j+1],\geq}$. The points $v^3$ and $(0, 1)$ are $x$-vertices of $Q^{[2j+1]}$ and satisfy $x_2\geq 1$. Hence they are $x$-vertices of $Q^{[2j+1],\geq}$. If $\lambda(u^{[2j+1]}) \leq \lambda(\bar w^2)$ (case (i)), then in particular we have $u^{[2j+1]}_2 \geq 1$ and therefore, as it is a vertex for $Q^{[2j+1]}$ and satisfies $x_2\geq 1$, it also is a $x$-vertex of $Q^{[2j+1],\geq}$. Observe that $\bar{w}^2$ is the intersection point of the ray $x \in \mathbb{R}^2\,|\,x = f + \lambda_2r^2, \lambda_2 \geq 0\}$ with the line $\{x \in \mathbb{R}^2\,|\,x_2 = 1\}$ and is therefore a $x$-vertex of $Q^{[2j+1],\geq}$. All other possible $x$-vertices come from the intersection of $\triangle((0, 1) v^3 p^{[2j+1]})$ with the line $\{x\in \mathbb{R}^2\,|\,x_2=1\}$ and yields $(0, 1)$ and $q^{[2j+2]}$.

In case (ii), the proof that $v^3,(0, 1)$, and $q^{[2j+2]}$ are $x$-vertices are the same. If $f_2 < 1$, then the ray $x \in \mathbb{R}^2\,|\,x = f + \lambda_2r^2, \lambda_2 \geq 0\}\cap \{x \in \mathbb{R}^2\,|\,x_2 = 1\} = \emptyset$. Also as $u^{2j +1}_2 < 1$, $u^{2j +1}$ is not an $x$-vertex of $Q^{[2j+1],\geq}$. If $f_2 \geq 1$ and $\lambda(u^{[2j+1]}) > \lambda(\bar w^2)$, then $\bar{w}^2 \in \textup{int}(\triangle((0,1)v^3p^{[2j+1]}))$ and therefore $\bar{w}^2$ is not a $x$-vertex of $Q^{[2j+1],\geq}$. Also as $u^{2j +1}_2 < 1$, $u^{2j +1}$ is not an $x$-vertex of $Q^{[2j+1],\geq}$. \cqfd

Remember $q$ is the intersection point of the line $w^2w^3$ with the line $\{x\in \mathbb{R}^2\,|\,x_2 = 1\}$.

\begin{lemma}[Finding $\mathbf{\beta^{[2j+2]}, r^3_1 > 0}$]$\quad$
\begin{enumerate}[(i)]
\item If $q^{[2j+2]}_1 \geq q_1$, then the goal inequality $\alpha$ is valid for $\conv(Q^{[2j+1],\leq}\cup Q^{[2j+1],\geq})$.
\item If $q^{[2j+2]}_1 < q_1$, then the inequality $\beta^{[2j+2]}:=\phi(\triangle(w^1(0, 1) q^{[2j+2]}))$ is valid for $\conv(Q^{[2j+1],\leq}\cup  Q^{[2j+1],\geq})$.
\end{enumerate}
\end{lemma}
\textbf{Proof:} See Figure \ref{rhs_movie} frame (e) for an example of case (ii).

We have to check that all the $x$-vertices given by Lemma \ref{vertices2j+1RHS} do not lie in
the interior of the corresponding triangle.
\begin{enumerate}
\item[(i)] We claim that we must be in case (b) of Lemma \ref{vertices2j+1RHS}, i.e if $q^{[2j+2]}_1 \geq q_1$, then $\lambda(u^{[2j+1]}) > \lambda(\bar w^2)$. Assume by contradiction that $\lambda(u^{[2j+1]}) \leq \lambda(\bar w^2)$. Then $q^{[2j +1]}_1 \leq \bar{w}^2$. As $\bar{w}^2_1 < q_1$, we obtain a contradiction.

    Observe that $L_{\alpha} = \triangle(w^1w^3q) = \triangle(w^1w^2w^3)$. It is now easy to verify that $v^3, q^{[2j +1]}\notin L_{\alpha}$, and $w^1, w^2, (0,1) \in \partial L_{\alpha}$, thus proving that the goal inequality $\alpha$ is valid for $\conv(Q^{[2j+1],\leq}\cup Q^{[2j+1],\geq})$.
\item[(ii)]  We distinguish between two cases.
\begin{enumerate}[(a)]
\item $q^{[2j+2]}_1 < \bar w^2_1$: (We are in Case (a) of Lemma \ref{vertices2j+1RHS}). As $w^1$, $(0,1)$, $q^{[2j +1]} \in \partial (\triangle(w^1 (0,1) q^{[2j+1]}))$, we need to verify that $v^3$, $w^2$, $\bar{w}^2$, $u^{[2j +1]} \notin \textup{int}(\triangle(w^1 (0,1)q^{[2j +1]}))$. As the line $w^1(0,1)$ meets the ray $\{x \in \mathbb{R}^2\,|\, x = f + \lambda_3 r^3, \lambda_3 \geq 0\}$ at $w^3$ and $\lambda(w^3) < \lambda(v^3)$, we obtain $v^3 \notin \textup{int}(\triangle(w^1 (0,1)q^{[2j +1]}))$.

    Let $u^{[2j +2]}$ be the intersection point of the line $wq^{[2j +2]}$ with the ray $\{x \in \mathbb{R}^2\,|\,x = f + \lambda_2 r^2, \lambda_2 \geq 0\}$. Observe that by construction of $q^{[2j +2]}$ the line $zu^{[2j +1]}$ also passes through $q^{[2j + 2]}$. Therefore as $\lambda(z) < \lambda(w^1)$ we have that $\lambda(u^{[2j +2]}) < \lambda(u^{[2j +1]})$. This proves that $u^{[2j + 1]}\notin \textup{int}(\triangle(w^1 (0,1)q^{[2j +1]}))$. As $\lambda(u^{[2j +1]}) < \lambda(\bar{w}^2) < \lambda(w^2)$, this completes the proof.
    \item $q^{[2j+2]}_1 \geq \bar w^2_1$ or $f_2 < 1$ (then $\bar w^2$ does not exist). Similar to the previous  case, all vertices given in Lemma \ref{vertices2j+1RHS} lie outside of $\textup{int}(\triangle(w^1(0, 1)q^{[2j+2]}))$. \cqfd
\end{enumerate}
\end{enumerate}

\begin{lemma}[Termination step, $\mathbf{r^3_1 > 0}$] Let $\beta^{[2j+1]} = \phi(\triangle(v^2v^3(0, 1)))$. Let $Q^{[2j +1], \leq}= (P(R,f))^0 \cap \{(x,s)\in \mathbb{R}^2 \times\mathbb{R}^3_{+}\,|\, \sum_{i =1}^3\beta^{[2j +1]}_is_i \geq 1, x_2\leq 0\}$ and $Q^{[2j +1], \geq}= (P(R,f))^0 \cap \{(x,s)\in \mathbb{R}^2 \times\mathbb{R}^3_{+}\,|\, \sum_{i =1}^3\beta^{[2j +1]}_is_i \geq 1, x_2\geq 1\}$. Then inequality $\alpha$ is valid for $\conv(Q^{[2j+1],\leq} \cup  Q^{[2j+1],\geq})$.
\label{ddsrhs}
\end{lemma}
\textbf{Proof:} See Figure \ref{rhs_movie} frame (g) for an example.

The vertices of $Q^{[2j+1],\leq} $ are $w^1$ and $w^2$. The vertices of $Q^{[2j+1],\geq} $ are $v^3, (0, 1)$ and $(1, 1)$. All these vertices do not lie in the interior of $L_{\alpha}$, which proves that $\alpha$ valid.\cqfd

Next we are concerned with the convergence of the Disjunction Sequence \ref{proofalgot2Ccase3}. Note that for the step $2j +1$ in Disjunction Sequence \ref{proofalgot2Ccase3}, there are two cases, i.e., $\lambda (u^{[2j +1]}) \leq \lambda(\bar{w}^2)$ or $\lambda(u^{[2j +1]}) > \lambda(\bar{w}^2)$ equivalently $u^{[2j +1]}_1 \leq \bar{w}^2_1$ or $u^{[2j +1]}_1 > \bar{w}^2_1$. Based on these two cases, the proof of convergence is divided into `two phases':

\begin{enumerate}
\item (Phase 1) In the first phase we prove that if $u^{[2]}_1 \leq \bar{w}^2_1$, then after a finite number of iterations $u^{[2j +1]}_1 > \bar{w}^2_1$ holds. This is proven in Lemma \ref{convergerhsphase1}.
\item (Phase 2) Since eventually $u^{[2j +1]}_1 > \bar{w}^2_1$ holds, we assume this to be true. With this assumption, it is shown that there exists a finite $j^{*}$ such that $p^{[2j^{*} +1]}$ is below $zv^2$. This completes the proof since this implies that after iteration $j^{*}$, $\phi(\triangle(v^2v^3(0 \, 1)))$ is valid, allowing the algorithm in Disjunction Sequence \ref{proofalgot2Ccase3} to enter the \textbf{Termination Step.} This is proven in Lemma \ref{convergerhsphase2}.
\end{enumerate}

\begin{center}
\begin{figure}[!h]
\begin{center}
\includegraphics[width=0.6\textwidth]{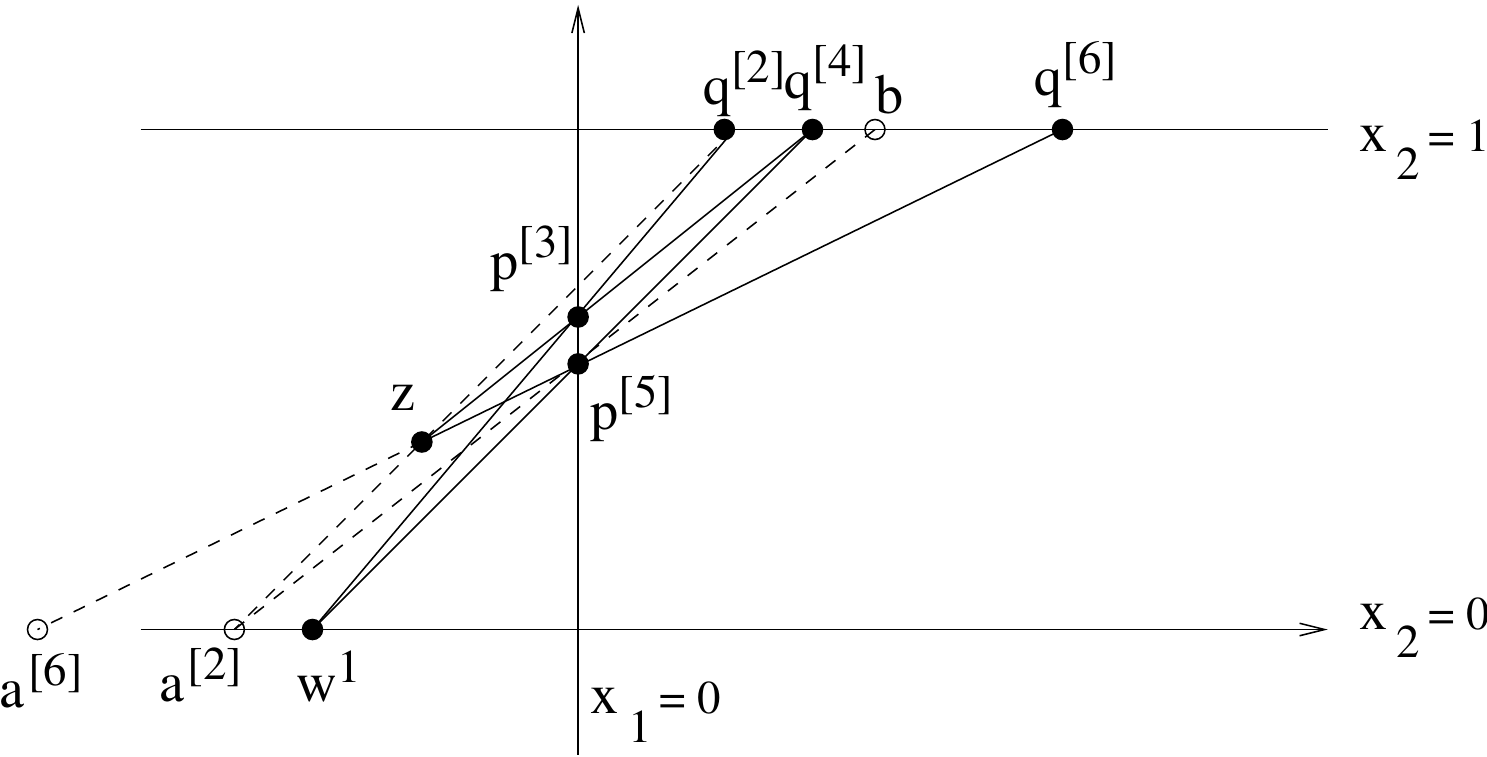}
\end{center}
\caption{Phase 1 convergence proof}
\label{phase1a}
\end{figure}
\end{center}

\begin{lemma}[Finite Convergence - Phase 1]\label{convergerhsphase1} If $u^{[2]}_1 \leq \bar{w}^2_1$, then after a finite number of iterations $j$, $u^{[2j +1]}_1 > \bar{w}^2_1$ holds.
\end{lemma}
\textbf{Proof:} (Refer to Figure \ref{phase1a}.) Let $a^{[2]} = (a^{[2]}_1, a^{[2]}_2)$ be the intersection point of the line $zq^{[2]}$ with the line $\{x \in \mathbb{R}^2\,|\,x_2 = 0\}$. Since $r^2_2 <0$ and $r^2_1 >0$, it can be verified that $a^{[2]}_1 < w^1_1 $. Furthermore $w^1_1 < z_1$ as $\lambda(w^1) > \lambda(z)$. Therefore we obtain $a^{[2]}_1 < w^1_1 < z_1 $. Observe that given $q^{[2j]}$, the point $q^{[2j +2]}$ can be generated as follows:
\begin{enumerate}
\item Join $q^{[2j]}$ and $w^1$ with a line segment. Let this line segment intersect the line $\{x\in \mathbb{R}^2\,|\,x_1 = 0\}$ at $p^{[2j + 1]}$.
\item Construct a line joining $z$ and $p^{[2j + 1]}$. If the half-line $zp^{[2j + 1]}:= \{x\in \mathbb{R}^2\,|\, x = z + \mu(p^{[2j +1]} - z), \mu \geq 0\}$ intersects the line $\{x\in \mathbb{R}^2\,|\,x_2 = 1\}$, then set $\tilde{q}^{[2j + 2]}$ as the intersection point. Else set $\tilde{q}^{[2j +2]}= (+ \infty, 1)$.
\item If $\tilde{q}^{[2j + 2]}_1 < \bar{w}^2_1$, then ${q}^{[2j + 2]}:= \tilde{q}^{[2j + 2]}$. Otherwise $u^{[2j +1]}_1 > \bar{w}^2_1$ holds.
\end{enumerate}
Suppose now that we construct the set of point $q^{[2]}$, $q^{[4]}$, $q^{[6]}$, ... where we set $\tilde{q}^{[2j +2]}$ to $q^{[2j]}$ always (without checking if $\tilde{q}^{[2j + 2]}_1 < \bar{w}^2_1$). Then proving that
\begin{eqnarray*}
\textup{lim}_{j \rightarrow \infty}q^{[2j]}_1 = +\infty
\end{eqnarray*}
proves the result of this lemma. This proves the result since then eventually (in finite number of steps) $\tilde{q}^{[2j]}_1 > \bar{w}^2_1$.

Using the fact that the line $q^{[2j]}w^1$ and $q^{[2j +2]}z$ intersect at $p^{[2j +1]}$ and $\lambda(w^1) > \lambda(z)$, it can be verified that $q^{[2j + 2]}_1 > q^{[2j]}_1$ $\forall j$. Now we verify that $|q^{[2j +2]}_1- q^{[2j]}_1 | < | q^{[2j +4]}_1 - q^{[2j +2]}_1|$ to complete the proof. We verify that $| q^{[2j +2]}_1- q^{[2j]}_1| < | q^{[2j +4]}_1 - q^{[2j +2]}_1|$ for $j = 1$. The proof is identical for any other $j$.

Let $b$ be the intersection point of the line passing through $a^{[2]}$ and $p^{[5]}$ and the line $\{x \in \mathbb{R}^2\,|\,x_2 =1\}$. (See Figure \ref{phase1a}). We claim that $q^{[4]}_1 < b_1 < q^{[6]}_1$. Since $a^{[2]}b$ and $w^1q^{[4]}$ intersects at $q^{[5]}$ and $a^{[2]}_1 < w^1_1$, we must have $q^{[4]}_1 < b_1$. Let $a^{[6]}$ be the intersection point of $zq^{[6]}$ with the line $\{x \in \mathbb{R}^2\,|\, x_2 = 0\}$. Since $q^{[2]}a^{[2]}$ and $q^{[6]}a^{[6]}$ intersect at $z$ and $q^{[2]}_1 < q^{[6]}_1$, we obtain that $a^{[6]}_1 < a^{[2]}_1$. Furthermore, as $a^{[6]}q^{[6]}$ and $a^{[2]}b$ intersect at $p^{[5]}$ we obtain that $b_1 < q^{[6]}_1$.

Therefore, $| b_1 - q^{[4]}_1| < |q^{[6]}_1 - q^{[4]}_1|$. Next we show that $|q^{[4]}_1 - q^{[2]}_1| < | b_1 - q^{[4]}_1|$ to complete the proof. It can be verified that $|b_1 - q^{[4]}_1| = (-a_1)(\frac{1}{p^{[5]}_2} - \frac{1}{p^{[3]}_2})$ and $|q^{[4]}_1 - q^{[2]}_1| = (-w^1_1)(\frac{1}{p^{[5]}_2} - \frac{1}{p^{[3]}_2})$. Since $(-a_1) > (- w^1_1)$, this completes the proof. \hfill $\square$

\begin{center}
\begin{figure}[!h]
\begin{center}
\includegraphics[width=0.5\textwidth]{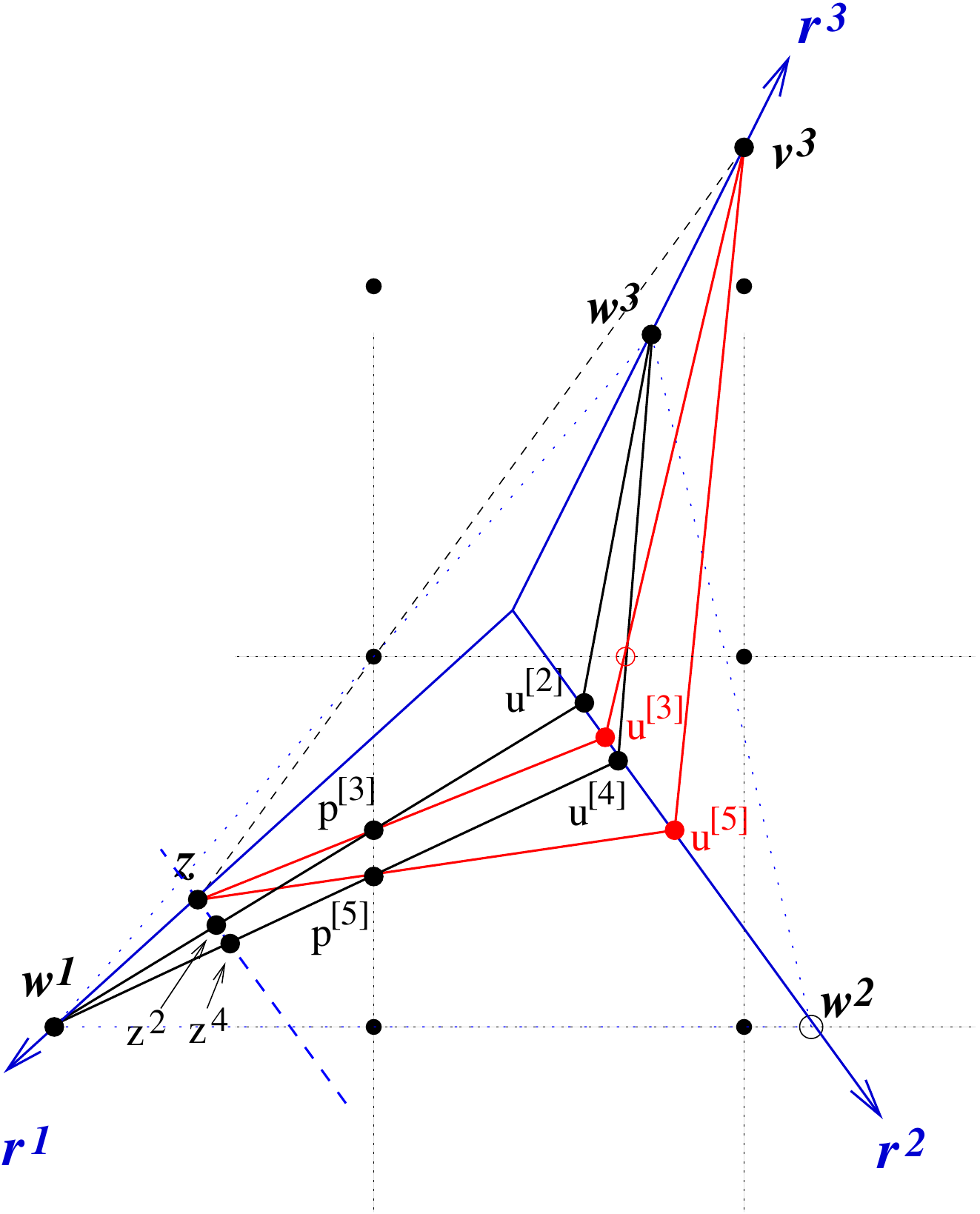}
\end{center}
\caption{Phase 2 convergence proof}
\label{phase2c}
\end{figure}
\end{center}

\begin{lemma}[Finite Convergence - Phase 2]\label{convergerhsphase2}
Let $u^{[2i +1]}_1 > \bar{w}^2_1$. There exists a finite $j^{*} \geq i$ such that $p^{[2j^{*} +1]}$ is below $zv^2$.
\end{lemma}
\textbf{Proof:} (Refer to Figure \ref{phase2c}). Using Disjunction Sequence \ref{proofalgot2Ccase3}, the points $p^{[2j]}$, $q^{[2j]}$, $u^{[2j + 1]}$, and $u^{[2j + 2]}$ are generated as follows:
\begin{enumerate}
\item $p^{[2j + 1]}$ is the intersection point of the line $u^{[2j]}w^{1}$ and the line $\{x\in \mathbb{R}^2\,|\,x_1 = 0\}$.
\item $u^{[2j + 1]}$ is the intersection point of the line $p^{[2j + 1]}z$ and the ray $\{x\in \mathbb{R}^2\,|\,x = f + \lambda r^2, \lambda \geq 0\}$.
\item $q^{[2j + 2]}$ is the intersection point of the line $u^{[2j + 1]}v^3$ and the line $\{x\in \mathbb{R}^2\,|\,x_2 = 1\}$.
\item $u^{[2j + 2]}$ is the intersection point of the line $q^{[2j + 2]}w^3$ and the ray $\{x\in \mathbb{R}^2\,|\,x = f + \lambda r^2, \lambda \geq 0\}$.
\end{enumerate}

Let $\lambda_{k}:= \lambda(u^{[k]})$. We will prove this result by showing that
\begin{eqnarray}\label{limphase2}
\textup{lim}_{k \rightarrow \infty}\lambda_k = +\infty.
\end{eqnarray}
This proves the result since it implies that there exists a finite $j^{*}$ such that $p^{[2j^{*} +1]}$ is below the line $zv^2$.

We prove (\ref{limphase2}) in the following steps:
\begin{itemize}
\item The sequence $\{\lambda_{k}\}_{k =1}^{\infty}$ is non-decreasing: Note that $u^{[k]}$ is generated differently when $k$ is odd and when $k$ is even. Therefore to prove this result, we show that:
\begin{itemize}
\item $\lambda_{2j} < \lambda_{2j + 1}$: Since $\lambda(w^1) > \lambda(z)$ and the lines $w^1u^{[2j]}$ and $zu^{[2j +1]}$ intersect at $p^{[2j +1]}$, we conclude that $\lambda_{2j} = \lambda(u^{[2j]})< \lambda(u^{[2j +1]}) = \lambda_{2j + 1}$.
\item $\lambda_{2j + 1} < \lambda_{2j + 2}$: Since $\lambda(v^3) > \lambda(w^3)$ and the lines $w^3u^{[2j+ 2]}$ and $v^3u^{[2j +1]}$ intersect at $q^{[2j +2]}$, we conclude that $\lambda_{2j + 1} = \lambda(u^{[2j + 1]})< \lambda(u^{[2j +2]}) = \lambda_{2j + 2}$.
\end{itemize}
\item $p^{[2j +1]}_2 > p^{[2j + 3]}_2$: Since from above, $\lambda_{2j + 3} > \lambda_{2j + 1}$ we obtain that $u^{[2j + 3]}_2< u^{[2j + 1]}_2$ and $u^{[2j + 3]}_1 > u^{[2j + 1]}_1$. Since $p^{[2j+3]}$ and $p^{[2j + 1]}$ are the intersections of the lines $zv^{[2j +3]}$ and $zv^{[2j + 1]}$ with the line $x_1 = 0$, we obtain the result.

\item $\lambda_{2j + 1} - \lambda_{2j} < \lambda_{2j+3} - \lambda_{2j+ 2}$. Since the sequence $\{\lambda_{k}\}_{k =1}^{\infty}$ is non-decreasing, this will complete the proof: We present the proof for $k = 1$, the proof is the same for all other values of $k$. Refer to Figure \ref{phase2c}. Construct a ray $\{x\in \mathbb{R}^2\,|\,x = z + \lambda r^2\}$, i.e., parallel to the ray $\{x \in \mathbb{R}^2\,|\,x = f + \lambda r^2\}$. Let $z^2$ and $z^4$ be the intersection of this ray with the line segments $p^{[3]}w^1$ and $p^{[5]}w^1$. Since $p^{[3]}_2 > p^{[5]}_2$ and $r^2_2 < 0$, we obtain that $\epsilon^2 = |zz^2| < |zz^4|= \epsilon^4$. Let $\bar{p}^{[3]}= z + \frac{\eta^2r^2}{|\!|r^2|\!|}$ be the orthogonal projection of $p^{[3]}$ on the line passing though $z$ and $z^2$. Similarly, let $\bar{p}^{[5]} = z + \frac{\eta^4r^2}{|\!|r^2|\!|}$ be the orthogonal projection of $p^{[5]}$ on the line passing though $z$ and $z^2$. Then since $p^{[5]}_2 < p^{[3]}_2$ and $r^3_2 < 0$, we obtain that $\eta^2 < \eta^4$ (Note $\eta^2$ and $\eta^4$ can be negative). Now it can be verified (see Appendix 1) that $|\!|u^{[3]}u^{[2]}|\!| = \epsilon^2(\frac{\gamma}{1 - \zeta\eta^2} -1)$ and $|\!|u^{[5]}u^{[4]}|\!| = \epsilon^4(\frac{\gamma}{1 - \zeta\eta^4} -1)$ where $\gamma$ and $\zeta$ are positive constants and $1 - \zeta\eta^2 >0$ and $1 - \zeta\eta^4 > 0$. Since $\epsilon^2 < \epsilon^4$ and $\eta^2 <\eta^4$ we obtain $|\!|u^{[3]}u^{[2]}|\!|< |\!|u^{[5]}u^{[4]}|\!|$. This proves that $\lambda_{2j + 1} - \lambda_{2j} < \lambda_{2j+3} - \lambda_{2j + 2}$. \hfill $\square$
\end{itemize}

\subsection{$T^{2A}$}\label{t2asec}

\begin{figure}[htbp]
\begin{center}
\includegraphics[width=0.3\linewidth]{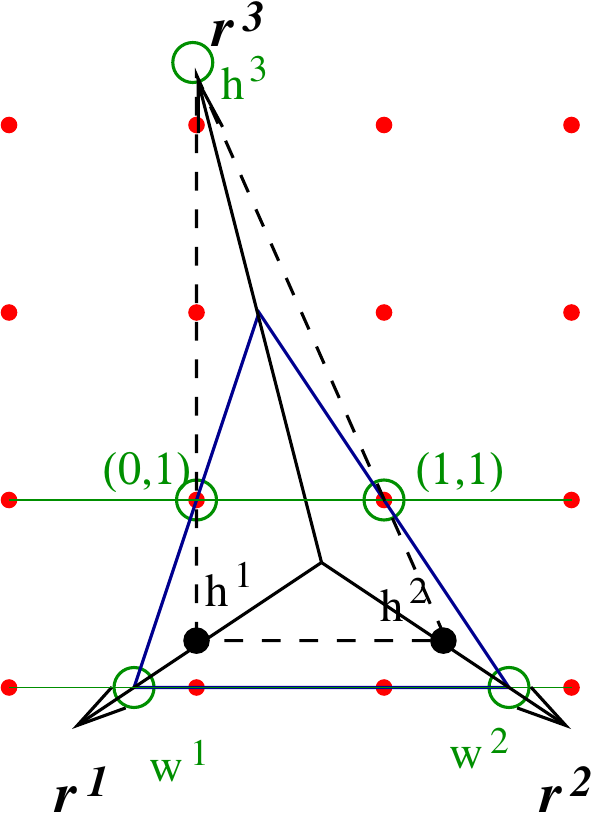}
\caption{A sketch of the proof that an  inequality whose induced lattice-free set is a triangle of type $T^{2A}$ has  finite split rank}\label{triMAX}
\end{center}
\end{figure}

Consider the triangle depicted in Figure \ref{triMAX} in solid lines. The convergence proofs used for Proposition \ref{theoremt2c} would not give the desired triangle in a finite number of steps.

In this section, we show that there exists an inequality of finite split rank (corresponding to the one depicted in dashed lines in Figure \ref{triMAX}) that together with the right split disjunction, provides the desired inequality. The inequality can be constructed as follows.
\begin{construction}[$\wedge^i$]
Let $p^1,p^2\in \Z^2,$ $\{i,j,k\}$ be a permutation of $\{1,2,3\}$ such that $p^1 = f + \lambda_ir^i + \lambda_j r^j$, $\lambda_i , \lambda_j \geq 0$ and $p^2 = f + \mu_ir^i + \mu_k r^k$, $\mu_i , \mu_k \geq 0$. We suppose that $\triangle(p_{\lambda}p^1p^2)$ exists and is lattice-free for some $\lambda >0$ where $p_{\lambda}= f+\lambda r^i$. Let $\bar \lambda := \sup \{\lambda \in \mathbb{R}_{+} \cup \{+\infty\}\mid \triangle(p_{\lambda} p^1p^2) \text{ is lattice free }\}$. If $\bar \lambda=+\infty$, we define $\wedge^i(p^1,p^2)$ as the lattice-free set determined by the two lines parallel to $r^i$ and incident to $p^1$ and $p^2$. If $\bar \lambda \in \R_{+}$, we define $\wedge^i(p^1,p^2):= \triangle(p_{\bar \lambda}p^1p^2)$ with $p_{\bar \lambda}= f+ \bar \lambda r^i$.
\label{wedgeconstruction}
\end{construction}
In Figure \ref{triMAX}, $\wedge^3((0, 1)(1, 1))$ is represented in dashed lines.

\begin{proposition}\label{t2aprop}
The split rank of an inequality whose induced lattice-free set is a triangle of type $T^{2A}$ is finite.
\end{proposition}
\textbf{Proof:} Let $\sum_{i = 1}^3\alpha_is_i \geq 1$ be a facet-defining inequality such that $L_{\alpha}$ is a triangle of type $T^{2A}$. By a suitable integral translation and unimodular transformation (Dey and Wolsey~\cite{deywolsey2008a}), we can assume that \begin{enumerate}
\item The vertices of $L_{\alpha}$ are
\begin{enumerate}
\item $w^1:= (-\delta , 0)$ where $0 < \delta \leq 1$ and $w^1=f+\lambda_1r^1, \lambda_1\geq 0,$
\item $w^2:= (g + \epsilon, 0)$ where $0 \leq \epsilon < 1$, $g\in \mathbb{Z}_{+}$, $g \geq 1$ and $w^2=f+\lambda_2r^2,\lambda_2\geq 0$,
\item $w^3:= (\bar{x}, \bar{y})$ where $\bar{y} > 1$ and $0 < \bar{x} < 1$, and $w^3=f+\lambda_3r^3,\lambda_3\geq 0$.
\end{enumerate}
\item The side $w^1w^3$ of $L_{\alpha}$  contains the integer point $(0,1)$ in its relative interior.
\item The side $w^2w^3$ of $L_{\alpha}$ contains the integer point $(1,1)$ in its relative interior.
\item $r^1_2<0,r^1_1<0,r^2_2<0,r^2_1\geq 0$ and $r^3_2>0$.
\end{enumerate}
\noindent Claim: Let $\beta:= \phi(\wedge^3((0, 1) (1, 1)))$. Then $L_{\beta}$ is either a split set or a triangle of type $T^{2B}$. If $L_{\beta}$ is not a split set, then let $h^i = f + \lambda_ir^i$, $\lambda_i \geq 0$ be the vertices of $L_{\beta}$. We show that $L_{\beta}$ is a triangle of type $T^{2B}$. Since $w^3$ is not an integral point, $\exists$ $\epsilon >0$ such that $\triangle((0,1)(1,1) (f + (\lambda(w^3) + \epsilon)r^3))$ is lattice-free. Therefore $\lambda(h^3) > \lambda(w^3)$. It follows therefore that $\lambda(h^1) <\lambda(w^1)$ and $\lambda(h^2) <\lambda(w^2)$. Hence the side $h^1h^2$ lies completely in the interior of $L_{\alpha}$ and does not contain any integer point. Moreover, by the maximality of $\lambda(h^3)$, either side $h^1h^3$ or side $h^2h^3$ (or both) contains at least two integer points. Thus, $L_{\beta}$ is a triangle of type $T^{2B}$.

Let $Q:= \{(x,s)\in\R^2\times \R^3_+ \mid x=f+Rs, \sum_{i=1}^3 \beta_i s_i \geq 1\}$ and let $Q^{\leq}:= Q\cap \{(x,s)\in \mathbb{R}^2\times \mathbb{R}^3\,|\,x_2\leq 0\}$ and $Q^{\geq}:=Q\cap \{(x,s)\in \mathbb{R}^2\times \mathbb{R}^3\,|\,x_2\geq 1\}$. Since $\beta$ is either a split cut or a cut whose induced lattice-free set is of type $T^{2B}$ it has a finite split rank. Thus proving that $\alpha$ is valid for $\textup{conv}(Q^{\leq} \cup Q^{\geq})$ proves the result.

We first enumerate the $x$-vertices of $Q^{\leq}$. Observe that since $h^1$ and $h^2$ are in the interior of $L_{\alpha}$ the only $x$-vertices of $Q^{\leq}$ are the intersection of the rays $\{x \in \mathbb{R}^2\,|\,x = f + \lambda_ir^i, \lambda_i \geq 0\}$, $i \in \{1, 2\}$ with the line $\{x\in \mathbb{R}^2\,|\, x_2 = 0\}$. These $x$-vertices are $(0, 0)$ and $(1, 0)$. Now consider the $x$-vertices of $Q^{\geq}$. They are $(0, 1), (1, 1)$ (at the intersection of $L_{\beta}$ and $\{x\in \mathbb{R}^2\,|\,x_2=1\})$ and $h^3$ (only if $L_{\beta}$ is not a subset of a split set). Since all $x$-vertices of $Q^{\leq}$ and $Q^{\geq}$ do not lie in the interior of $L_{\alpha}$, the result follows.\cqfd

A class of inequalities known as the sequential-merge inequalities were introduced in Dey and Richard~\cite{dey:2007}. The induced lattice-free set of sequential-merge inequality using two Gomory mixed integer cuts applied to $P(R,f)$ is a triangle of type $T^{2A}$, see Dey and Wolsey~\cite{deywolsey2008a}. By their construction, the split rank of sequential-merge inequalities with two Gomory mixed integer cuts is at most 2. The procedure implied by Proposition \ref{t2aprop} can be shown to also imply a split rank of 2 for these inequalities.
%
%
%
%
\subsection{$T^3$}\label{t1sec}
So far, we have considered a proof technique based on keeping one inequality after each split disjunction. In this section, we need to keep two inequalities before a particular split disjunction is considered. The set reads as
\begin{eqnarray*}
Q:= \{(x,s) \in \mathbb{R}^2 \times \mathbb{R}^2_{+}\, &|&\,(x,s) \in (P(R,f))^0\nonumber \\
&&\beta^1_1s_1 + \beta^1_2s_2 + \beta^1_3s_3 \geq 1 \nonumber \\
&&\beta^2_1s_1+\beta^2_2s_2+\beta^2_3s_3\geq 1 \}.\nonumber
\end{eqnarray*}
Observe that for this set all the vertices are of the form $(x,\mathcal{M}(x))$. In particular, any vertex that is tight for both $\beta^1$ and $\beta^2$ must be of the form $(v,\mathcal{M}(v))$ where $v$ is an intersection point of $\partial L_{\beta^1}$ and $\partial L_{\beta^2}$. The extreme rays of $Q$ are of the form $(r^j, e^j)$, $j \in \{1, 2,3\}$ where $e^j\in \mathbb{R}^3_{+}$ is the unit vector in the direction of the $j^{\textrm{th}}$ coordinate. Since either $\beta^1_j$ or $\beta^2_j >0$ for all $j$ (both $\beta^1$ and $\beta^2$ are not the same split inequality), any $(\bar{x}, \bar{s}) \in Q$ that is tight for both $\beta^1$ and $\beta^2$ must be a convex combination of points of the form $(v^k,\mathcal{M}(v^k))$ where $v^k$ is the intersection point of $\partial L_{\beta^1}$ and $\partial L_{\beta^2}$. This observation is useful in determining the vertices of sets of the form $Q^{\leq}:= Q\cap \{(x,s)\in \mathbb{R}^2\times\mathbb{R}^3\,|\,\pi_1 x_1 + \pi_2 x_2 \leq \pi_0\}$ and $Q^{\geq}:= Q\cap \{(x,s)\in \mathbb{R}^2\times\mathbb{R}^3\,|\, \pi_1 x_1 + \pi_2 x_2 \geq \pi_0 + 1\}$.

\begin{figure}[htbp]
\begin{center}
\includegraphics[width=0.6\linewidth]{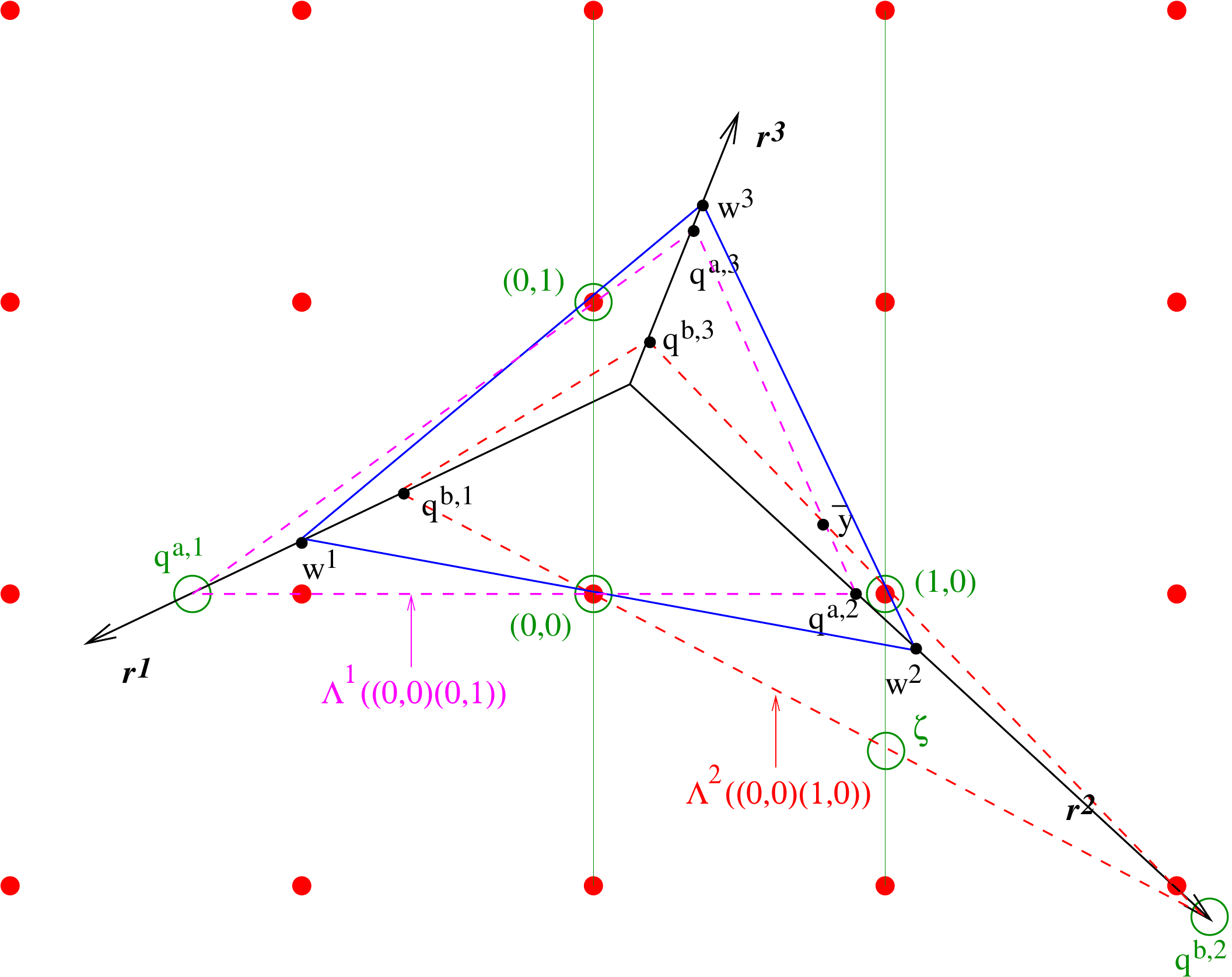}
\caption{A sketch of the proof that a cut whose induced lattice-free set is a triangle of type $T^{3}$ has  finite split rank}
\label{tridis}
\end{center}
\end{figure}
\begin{proposition}\label{t3proof}
The split rank of an inequality whose induced lattice-free set is a triangle of type $T^3$ is finite.
\end{proposition}
\textbf{Proof:} Let $\sum_{i = 1}^3\alpha_is_i \geq 1$ be a facet-defining inequality such that $L_{\alpha}$ is a triangle of type $T^{3}$. By a suitable integral translation and unimodular transformation (Dey and Wolsey~\cite{deywolsey2008a}), we can assume that \begin{enumerate}
\item The vertices of $L_{\alpha}$ are $w^i=f+\lambda_ir^i, \lambda_i\geq 0$ $i \in \{1, 2,3\}$,
\item The sides $w^1w^2$, $w^2w^3$, and $w^3w^1$ of $L_{\alpha}$ contain the integer points $(0,0)$, $(1,0)$, and $(0,1)$ in their relative interior respectively.
\end{enumerate}

Using Construction \ref{wedgeconstruction}, we define $\beta^a:= \wedge^1((0, 0) (0, 1))$ and $\beta^b:= \wedge^2((0, 0)(1, 0))$. Similarly to the proof in Proposition \ref{t2aprop}, it can be verified that $L_{\beta^a}$ and $L_{\beta^b}$ are either a subset of a split set or a triangle of type $T^{2B}$. Hence by Proposition \ref{theoremt2c}, their split rank is finite.

Let $q^{a,i} = f + \lambda_ir^i$, $\lambda_i \geq 0$ be the vertices of $L_{\beta^a}$ and let $q^{b,i} = f + \mu_ir^i$, $\mu_i \geq 0$ be the vertices of $L_{\beta^b}$. In the rest of the proof we assume that $ L_{\beta^a}$ and $ L_{\beta^b}$ are not subsets of a split set (i.e,  $q^{a,1}$, $q^{b,2}$ are well-defined). This is for simplicity and the proof can be modified for the cases where $ L_{\beta^a}$ and $L_{\beta^b}$ are subsets of split sets.

Observe that by construction of $\beta^a$ and $\beta^b$, we obtain that $\lambda(q^{a,1}) > \lambda (w^1) > \lambda(q^{b,1})$ and $\lambda(q^{a,2}) < \lambda(w^2) < \lambda(q^{b,1})$. We first present a key result.

\noindent Claim: $\partial L_{\beta^a}$ and $\partial L_{\beta^b}$ intersect at two points: $(0,0)$ and $(\bar{y}_1, \bar{y}_2)$ where $0 <\bar{y}_1 <1$: Indeed, one point of intersection is $(0, 0)$ by construction. Let us look for other potential intersection points. Since the side $q^{b,1}q^{b,3}$ lies entirely in the interior of $L_{\alpha}$, we can verify that $(\partial (L_{\beta^a}) \cap \{x \in \mathbb{R}^2\,|\, x_1 \leq 0\}) \cap ( (\partial L_{\beta^b}) \cap \{x \in \mathbb{R}^2\,|\, x_1 \leq 0\} ) = (0,0)$. Thus $\bar{y}_1 > 0$. Similarly, $(\partial (L_{\beta^a}) \cap \{x \in \mathbb{R}^2\,|\, x_2 \leq 0\}) \cap ( (\partial L_{\beta^b}) \cap \{x \in \mathbb{R}^2\,|\, x_2 \leq 0\} ) = (0,0)$. Thus $\bar{y}_2 > 0$. Using similar arguments we can verify that $(\bar{y}_1, \bar{y}_2) \in L_{\alpha}$. Together with the fact that $\bar{y}_2 > 0$, we obtain $\bar{y}_1 <1$.

We now consider the set
\begin{align*}
Q=\{(x,s)\in \R^2\times \R^3_+ \mid & (x,s)\in (P(R,f))^0\\
& \sum_{i=1}^3 \beta^a_i s_i \geq 1\\
& \sum_{i=1}^3 \beta^b_i s_i \geq 1\quad \}.
\end{align*}
Let $Q^{\leq}:= Q\cap \{(x,s) \in \mathbb{R}^2\times \mathbb{R}^3\,|\,x_1\leq 0\}$ and $Q^{\geq}:= Q\cap \{(x,s) \in \mathbb{R}^2\times \mathbb{R}^3\,|\,x_1\geq 1\}$. We show that $\alpha$ is a valid inequality for $\textup{conv}(Q^{\leq} \cup Q^{\geq})$, thus proving the result.

Claim: The vertices of $Q^{\leq}$ are $(q^{a,1}, \mathcal{M}(q^{a,1}))$, $((0,0), \mathcal{M}((0,0))$, $(v^3, \mathcal{M}(v^3))$ (if $r^3_1 < 0$, then let $v^3$ be the intersection point of the line $\{x \in \mathbb{R}^2\,|\,x_1 = 0\}$ with the ray $\{x \in \mathbb{R}^2\,|\,x = f + \lambda_3 r^3, \lambda_3 \geq 0\}$) and $((0,1), \mathcal{M}((0,1))$. Let $(\bar{x}, \bar{s})$ be a vertex of $Q^{\leq}$. We consider the following cases:
\begin{enumerate}
\item Support of $\bar{s}$ is 1, where $(\bar{x}, \bar{s})$ is tight for $\beta^a$: $(q^{a,1}, \mathcal{M}(q^{1,1}))$.
\item Support of $\bar{s}$ is 1, where $(\bar{x}, \bar{s})$ is tight for $\beta^b$: No such vertex, since $\lambda(q^{a,1}) >  \lambda(q^{b,1})$ (i.e, the point $(q^{b,1}, \mathcal{M}(q^{b,1}))$ is cut off by the inequality $\sum_{i = 1}^3\beta^a_is_i \geq 1$).
\item Support of $\bar{s}$ is 1, where $(\bar{x}, \bar{s})$ is tight for $\{x \in \mathbb{R}^2\,|\,x_1 = 0\}$: Assume that $r^3_1 < 0$. Then $(v^3, \mathcal{M}(v^3))$ is a vertex.
\item Support of $\bar{s}$ is 2, where $(\bar{x}, \bar{s})$ is tight for $\beta^a$ and $\beta^b$: Note that if $(\bar{x}, \bar{s})$ is tight for $\beta^a$ and $\beta^b$, then it must be a convex combination of the point $((0,0), \mathcal{M}((0,0))$ and the point $((\bar{y}_1, \bar{y}_2), \mathcal{M}((\bar{y}_1, \bar{y}_2))$. By the previous claim, $\bar{y}_1 > 0$. Therefore the only vertex is $((0,0), \mathcal{M}((0,0))$.
\item Support of $\bar{s}$ is 2, where $(\bar{x}, \bar{s})$ is tight for $\beta^a$ and $\{x\in \mathbb{R}^2\,|\,x_1 = 0\}$: $((0,1), \mathcal{M}((0,1))$.
\item Support of $\bar{s}$ is 2, where $(\bar{x}, \bar{s})$ is tight for $\beta^b$ and $\{x\in \mathbb{R}^2\,|\,x_1 = 0\}$: No such vertex. This is because if $\bar{x}$ is the intersection point of $\partial L_{\beta^b}$ and the line $\{x\in \mathbb{R}^2\,|\, x_1 = 0\}$, then $\bar{x}_2 < 1$. Therefore, this point is cut off by the inequality $\sum_{i = 1}^3\beta^a_is_i \geq 1$.
\item Support of $\bar{s}$ is 3, where $(\bar{x}, \bar{s})$ is tight for $\beta^a$, $\beta^b$, and $\{x \in \mathbb{R}^2\,|\,x_1 = 0\}$: Then $(\bar{x}, \bar{s})$ must be convex combination of the points $((0,0), \mathcal{M}((0,0)))$ and $((\bar{y}_1, \bar{y}_2), \mathcal{M}((\bar{y}_1, \bar{y}_2)))$, where $(\bar{y}_1, \bar{y}_2)$ is the intersection point of $\partial L_{\beta^a}$ and $\partial L_{\beta^b}$ different from $(0,0)$. Since from the previous claim, $\bar{y}_1 >0$, we obtain that $\bar{x} = (0,0)$. Therefore, no such vertex.
\end{enumerate}
\noindent Claim:  The vertices of $Q^{\geq}$ are $(q^{b,2}, \mathcal{M}(q^{b,2}))$, $((1,0), \mathcal{M}((1,0))$, $(v^3, \mathcal{M}(v^3))$ ((if $r^3_1 > 0$ then let $v^3$ be the intersection point of the line $\{x \in \mathbb{R}^2\,|\,x_1 = 1\}$ with the ray $\{x \in \mathbb{R}^2\,|\,x = f + \lambda_3 r^3, \lambda_3 \geq 0\}$) and $(\zeta, \mathcal{M}(\zeta))$ where $\zeta \neq (1,0)$ and it is an intersection point of $\partial L_{\beta^b}$ and the line $\{x\in \mathbb{R}^2\,|\,x_1=1\}$. Let $(\bar{x}, \bar{s})$ be a vertex of $Q^{\geq}$. We consider the following cases:
\begin{enumerate}
\item Support of $\bar{s}$ is 1, where $(\bar{x}, \bar{s})$ is tight for $\beta^a$: No such vertex, since $\lambda(q^{a,2}) <  \lambda(q^{b,2})$ (i.e, the point $(q^{a,2}, \mathcal{M}(q^{a,2}))$ is cut off by the inequality $\sum_{i = 1}^3\beta^b_is_i \geq 1$).
\item Support of $\bar{s}$ is 1, where $(\bar{x}, \bar{s})$ is tight for $\beta^b$: $(q^{b,2}, \mathcal{M}(q^{b,2}))$.
\item Support of $\bar{s}$ is 1, where $(\bar{x}, \bar{s})$ is tight for $\{x \in \mathbb{R}^2\,|\,x_1 = 1\}$: Assume that $r^3_1 > 0$. Then $(v^3, \mathcal{M}(v^3))$ is a vertex.
\item Support of $\bar{s}$ is 2, where $(\bar{x}, \bar{s})$ is tight for $\beta^a$ and $\beta^b$: Since $\bar{y}_1 < 1$, no such vertex.
\item Support of $\bar{s}$ is 2, where $(\bar{x}, \bar{s})$ is tight for $\beta^a$ and $\{x \in \mathbb{R}^2\,|\,x_1 = 1\}$: Then $\bar{x}$ must lie in the interior of $L_{\beta^b}$ and therefore there is no such vertex.
\item Support of $\bar{s}$ is 2, where $(\bar{x}, \bar{s})$ is tight for $\beta^b$ and $\{x \in \mathbb{R}^2\,|\, x_1 = 1\}$: $((1,0), \mathcal{M}((1,0))$ and $(\zeta, \mathcal{M}(\zeta))$.
\item Support of $\bar{s}$ is 3, where $(\bar{x}, \bar{s})$ is tight for $\beta^a$, $\beta^b$, and $\{x \in \mathbb{R}^2\,|\,x_1 = 1\}$: Then $(\bar{x}, \bar{s})$ must be convex combination of the points $((0,0), \mathcal{M}((0,0)))$ and $((\bar{y}_1, \bar{y}_2), \mathcal{M}((\bar{y}_1, \bar{y}_2)))$, where $(\bar{y}_1, \bar{y}_2)$ is the intersection point of $\partial L_{\beta^a}$ and $\partial L_{\beta^b}$ different from $(0,0)$. Since from the previous claim, $\bar{y}_1 <1$, we obtain no such vertex.
\end{enumerate}

Finally, observe that all the vertices of $Q^{\leq}$ and $Q^{\geq}$ are of the form $(\bar{x}, \mathcal{M}(\bar{x}))$. Moreover, $\bar{x} \notin \textup{int}(L_{\alpha})$. By Proposition \ref{therelight}, $\alpha$ is therefore valid  for $\textup{conv}(Q^{\leq} \cup Q^{\geq})$, thus completing the proof. \cqfd

The mixing set introduced by G\"unl\"uk and Pochet~\cite{gunluk01} correspond to $P(R,f)$ with a specific $R$ and $f$. The induced lattice-free set of mixing inequalities corresponding to the mixing set with two rows is a triangle of the type $T^3$; see Dey and Wolsey~\cite{deywolsey2008a} and Dey~\cite{deysplitrank}. An upper bound to the split rank of this inequalities is proven to be two in Dash and G\"unl\"uk~\cite{dashgunlukmixingrank}. It can be verified that the split rank implied by the proof of Proposition \ref{t3proof} is also 2. This bound is tight as shown in Andersen et. al~\cite{andersen:2007}, Dey~\cite{deysplitrank}, and Dash and G\"unl\"uk~\cite{dashgunlukmixingrank}.

\section{Four Variable Problems}\label{fourvarsec}
In this section, we consider the split rank of inequalities $\sum_{i = 1}^4\alpha_is_i \geq 1$ for $P([r^1, r^2, r^3, r^4],f)$ where $\textup{cone}\{r^1, ..., r^4\} = \mathbb{R}^2$.

\subsection{$Q^1$}
This class of inequalities corresponds to $L_{\alpha}$ being a quadrilateral with one side containing more than one integer point, two sides containing at least one integer point and the fourth side not containing any integer point in its relative interior.
\begin{figure}[htbp]
\begin{center}
\includegraphics[width=0.5\linewidth]{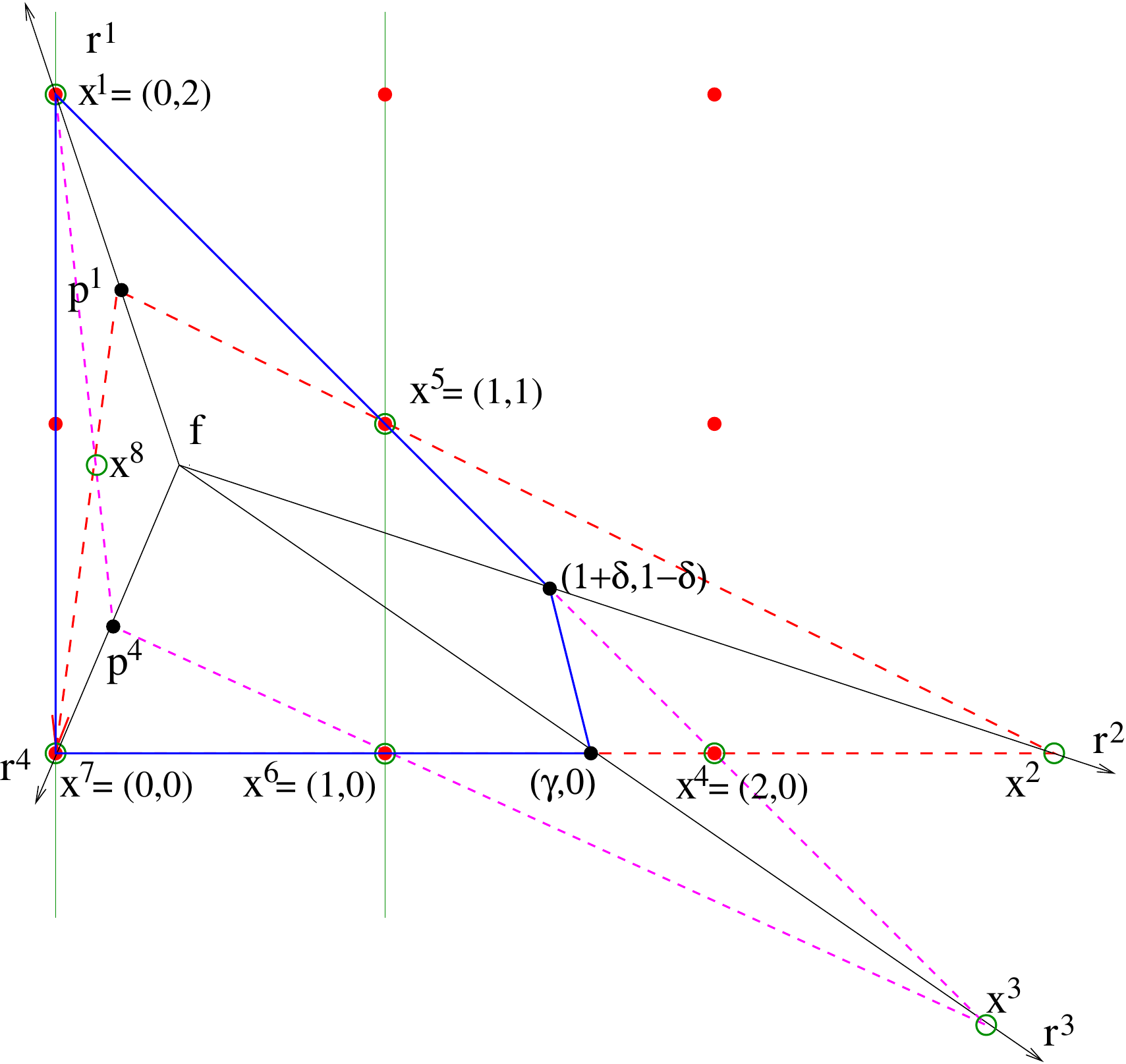}
\caption{A sketch of the proof that a cut whose induced lattice-free set is a quadrilateral of type $Q^{1}$ has  finite split rank}
\label{quadq1}
\end{center}
\end{figure}

\begin{proposition}\label{q1prop}
The split rank of an inequality whose induced lattice-free set is a quadrilateral of type $Q^1$ is finite.
\end{proposition}
\textbf{Proof:} Let $\sum_{i = 1}^3\alpha_is_i \geq 1$ be a facet-defining inequality such that $L_{\alpha}$ is a quadrilateral of type $Q^{1}$. As discussed in Section \ref{presection}, in this case $L_{\alpha}$ is a subset of a lattice-free triangle of type $T^1$ or $T^{2A}$. If it is a subset of a lattice-free triangle of type $T^{2A}$, then using Proposition \ref{t2aprop} and Lemma \ref{sameshape}, the split rank is finite. Therefore we consider the case where $L_{\alpha}$ is a proper subset of lattice-free triangle of type $T^1$.

By a suitable integral translation and unimodular transformation, we can assume that the vertices of $L_{\alpha}$ are the following: $(0,2) = f + \lambda_1 r^1$, $(1 + \delta, 1 - \delta)=f + \lambda_2 r^2$, $(\gamma, 0)= f + \lambda_3 r^3$, $(0,0)= f + \lambda_4 r^4$ where $0< \delta < 1$ and $1 < \gamma  < 2$ and $\lambda_i > 0$. See Figure \ref{quadq1}.

We may assume that $r^2_2 < 0$ and $r^3_1 \geq - r^3_2$. (If $r^2_2 \geq 0$, then consider a new set with all the same data except with a negative value of $r^2_2$ such that the ray $f + \lambda r^2$ passes through a point of the form $\mu (1,1)+ (1- \mu)(2,0)$ where $0 < \mu < 1$. The induced lattice-free set of an inequality $\sum_{i = 1}^4\tilde{\alpha}_is_i \geq 1$ where $\tilde{\alpha}_i = \alpha_i$ for $i\in \{1, 3, 4\}$ and $\tilde{\alpha}_2$ such that $f + \frac{r^2}{\tilde{\alpha}_2}$ lies on the line $\{x\in \mathbb{R}^2\,|\, x_1 + x_2 = 2\}$ is larger that $L_{\alpha}$, and thus by Lemma \ref{sameshape} its split rank is at least as large as the split rank of the original inequality $\alpha$. By similar argument, we may assume that $r^3_1 \geq - r^3_2$).

Now we consider the following two inequalities.
\begin{enumerate}
\item Let $x^2$ be the intersection point of the line $\{x\in \mathbb{R}^2\,|\,x_2 = 0\}$ with the ray $\{x\in \mathbb{R}^2\,|\, x = f + \lambda r^2$, $\lambda > 0\}$. Let $p^1$ be the intersection point of the line $x^2(1,1)$ with the ray $\{x\in \mathbb{R}^2\,|\, x =f + \lambda r^1$, $\lambda > 0\}$. Note now that the triangle with vertices $p^1$, $x^2$, $(0,0)$ is lattice-free (call this triangle $U$). This is because, it is a subset of the set $\{(x_1, x_2)\,|\, 0 \leq x_2 \leq 1\} \cup \{(x_1, x_2)\,|\, 0 \leq x_1 \leq 1\}$. Consider the inequality $\phi(U)$ for $P([r^1, r^2, r^3, r^4],f)$. Observe that since $p^1 \notin \mathbb{Z}^2$, $U$ is a triangle of type $T^{2A}$. Hence by Proposition \ref{t2aprop}, the split rank of $\phi(U)$ is finite. Denote $\beta^1 := \phi(U)$.
\item Let $x^3$ be the intersection point of the line $\{x\in \mathbb{R}^2\,|\,x_1 + x_2 = 2\}$ and the ray $\{x\in \mathbb{R}^2\,|\, x = f + \mu r^3$, $\mu \geq 0\}$. (Since $r^3_1 \geq - r^3_2$, this intersection exists). Let $p^4$ be the intersection point of the line $x^3(1,0)$ with the ray $\{x\in \mathbb{R}^2\,|\, x = f + \mu r^4$, $\mu \geq 0\}$. Note now that the triangle with vertices $(0,2)$, $x^3$ $p^4$ is lattice-free (call this triangle $V$). This is because, it is a subset of the set $\{(x_1, x_2)| 0 \leq x_1 \leq 1\} \cup \{(x_1, x_2)| 1 \leq x_1 + x_2 \leq 2\}$. Consider the inequality $\phi(V)$ for $P([r^1, r^2, r^3, r^4],f)$. Observe that since $p^4 \notin \mathbb{Z}^2$, $V$ is a triangle of type $T^{2A}$. Hence by Proposition \ref{t2aprop}, the split rank of $\phi(V)$ is finite. Denote $\beta^2: = \phi(V)$.
\end{enumerate}

We now consider the following set
\begin{align*}
Q:=\{(x,s)\in \R^2\times \R^4_+ \mid & x = f + Rs\\
& \sum_{i=1}^4 \beta^1_i s_i \geq 1\\
& \sum_{i=1}^4 \beta^2_i s_i \geq 1\quad \}.
\end{align*}
Since every vertex of $Q$ has a support of 2 (for the $s$-variables) and since $r^2_2 < 0$ and $r^3_1 \geq - r^3_2$, it can be verified that the vertices of this system are:
\begin{enumerate}
\item $(x^1,s^1) := ((0,2), \mathcal{M}^{1,1}((0,2)))$
\item $(x^2,s^2) := (x^2, \mathcal{M}^{2,2}(x^2))$
\item $(x^3,s^3) := (x^3, \mathcal{M}^{3,3}(x^3))$
\item $(x^4,s^4) := ((2,0), \mathcal{M}^{2,3}((2,0)))$
\item $(x^5,s^5) := ((1,1), \mathcal{M}^{1,2}((1,1)))$
\item $(x^6,s^6) := ((1,0), \mathcal{M}^{3,4}((1,0)))$
\item $(x^7,s^7) := ((0,0), \mathcal{M}^{4,4}((0,0)))$
\item $(x^8,s^8) := (p^{1,4}, \mathcal{M}^{3,4}(p^{1,4}))$, where $p^{1,4}$ is the intersection point of the lines $p^4(0,2)$ and $p^1(0,0)$.
\end{enumerate}

Let $Q^{\leq}: =Q \cap \{(x, s) \in \mathbb{R}^2 \times \mathbb{R}^4\,|\,x_1 \leq 0 \}$ and $Q^{\geq}: =Q \cap \{(x, s)\in \mathbb{R}^2 \times \mathbb{R}^4\,|\,x_1 \geq 0 \}$.

It can be verified that the vertices of $Q^{\leq}$ are $(x,s) := ((0,2), \mathcal{M}^{1,1}((0,2)))$ and $(x,s) := ((0,0), \mathcal{M}^{4,4}((0,0)))$. By Proposition \ref{therelight}, these points satisfy the goal inequality $\alpha$.

Now consider $Q^{\geq}$. Any vertex of $Q^{\geq}$ is of the form $\sum_{j = 1}^8 \lambda_j (x^j, s^j) + \sum_{k = 1}^4\mu_k(r^k, e^k)$ (where $(r^k, e^k)$ is an extreme ray of $Q$). Also note that any vertex of $Q^{\geq}$ must have a support of at most 3 on the $s$-variables.
\begin{enumerate}
\item Vertices of support 3: Such a vertex of $Q^{\geq}$ is satisfied at equality by the constraints $\sum_{j = 1}^4\beta^1_js_j \geq 1$, $\sum_{j = 1}^4 \beta^2_js_j \geq 1$, and $x_1 \geq 1$. In particular such a vertex is of the form $\lambda_5(x^5,s^5) + \lambda_4(x^4,s^4) + \lambda_6(x^6, y^6) + \lambda_8(x^8,s^8)$ (where $0 \leq \lambda_5, \lambda_4, \lambda_6, \lambda_8 \leq 1$ and $\lambda_5 + \lambda_4 + \lambda_6 + \lambda_8 = 1$) since these are the only vertices of $Q$ that are tight for both the inequalities $\sum_{j = 1}^4\beta^1_js_j \geq 1$ and $\sum_{j = 1}^4 \beta^2_js_j \geq 1$. (Also note that $\mu_k = 0$ for all $k\in \{1, ..., 4\}$, since otherwise the inequalities $\sum_{j = 1}^4\beta^1_js_j \geq 1$ and $\sum_{j = 1}^4 \beta^2_js_j \geq 1$ will not be satisfied at equality). If $\lambda_8 > 0$, to satisfy the constraint $x_1 \geq 1$, we must have that $\lambda_4 >0$. However, this makes the support of the resulting point 4, a contradiction. Thus, any vertex of support 3 of $Q^{\geq}$ is of the form $\lambda_5(x^5,s^5) + \lambda_4(x^4,s^4) + \lambda_6(x^6, y^6)$. Since the points $(x^5,s^5)$, $(x^4, s^4)$, $(x^6,s^6)$ satisfy the goal inequality, all vertices of support 3 satisfy the goal inequality.
\item Vertices of support 2: Let it be of the form $\sum_{j = 1}^8 \lambda_j (x^j, s^j) + \sum_{k = 1}^4\mu_k(r^j, e^k)$. If $\lambda_8 > 0$, then either $\lambda_2 > 0$, $\lambda_4 > 0$, $\lambda_3>0$, or $\mu_2 > 0$ or $\mu_3 > 0$ to satisfy the constraint $x_1 \geq 1$. However, this makes the support 3, a contradiction. Thus, the vertex is of the form $\sum_{j = 1}^7 \lambda_j (x^j, s^j) + \sum_{k = 1}^4\mu_k(r^j, e^k)$. Since the points $(x^j,s^j)$, $j \in \{1, ..., 7\}$ satisfy the goal inequality, all vertices of support 2 satisfy the goal inequality.
\item Vertex of support 1: Proof similar to the above case.
\end{enumerate}
So any vertex of $Q^{\geq}$ is valid for the goal inequality $\alpha$, completing the proof. \hfill $\square$

\subsection{$Q^2$}

\begin{figure}[htbp]
\begin{center}
\includegraphics[width=\linewidth]{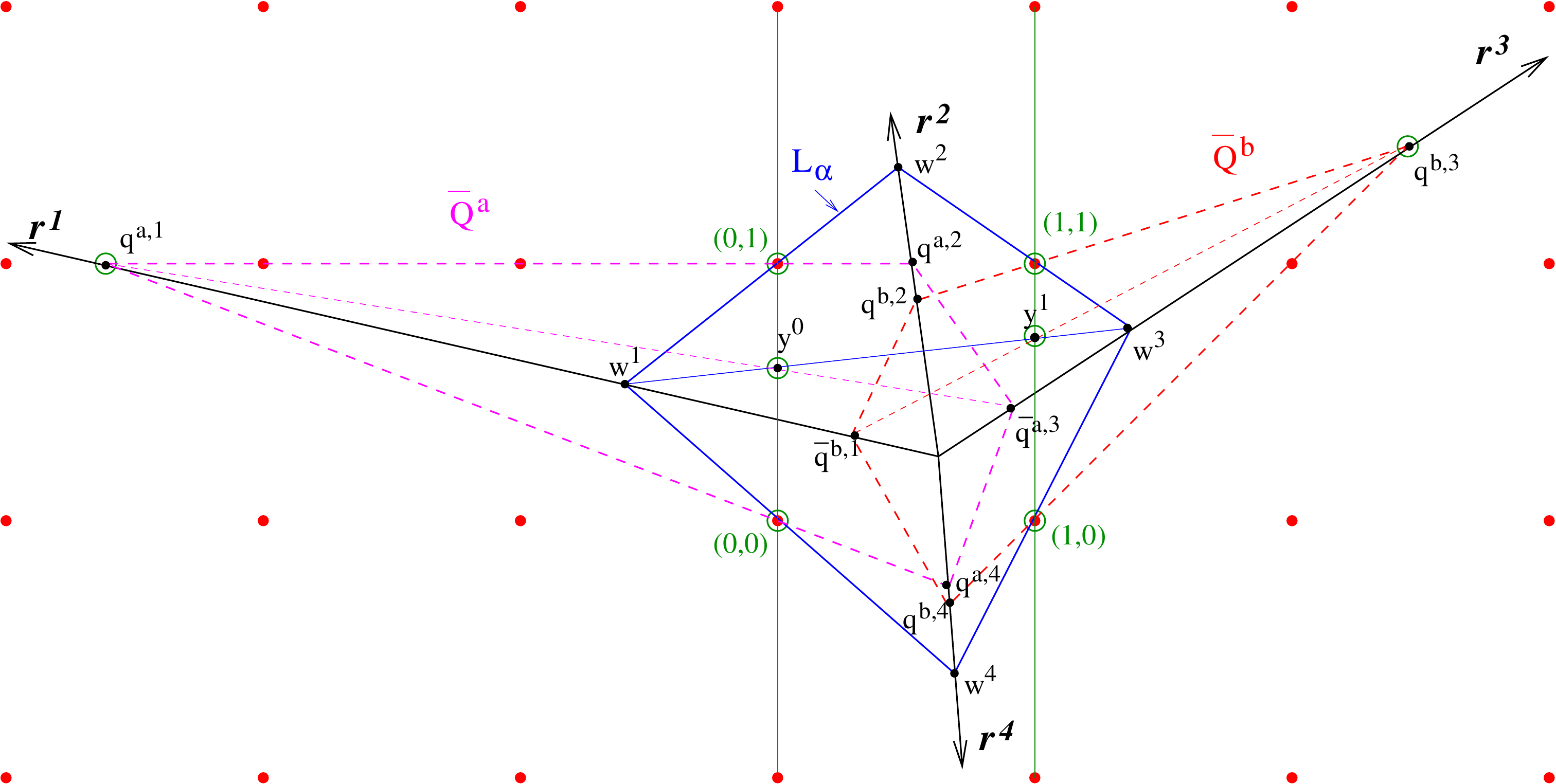}
\caption{A sketch of the proof that a cut whose induced lattice-free set is a quadrilateral of type $Q^{2}$ has  finite split rank}
\label{diq}
\end{center}
\end{figure}

Let $\sum_{i = 1}^3\alpha_is_i \geq 1$ be a facet-defining inequality such that $L_{\alpha}$ is a quadrilateral of type $Q^{1}$. By a suitable integral translation and unimodular transformation, we can assume that the boundary of $L_{\alpha}$ passes through the four integer points $(0, 0), (1, 0), (0, 1), (1, 1)$ where $(0, 0) = f + \lambda_1 r^1 + \lambda_4 r^4$, $\lambda_1, \lambda_4 \geq 0$; $(1, 0)= f + \mu_3 r^3 + \mu_4r^4$ $\mu_3, \mu_4 \geq 0$; $(0, 1) = f + \nu_1 r^1 + \nu_2 r^2$, $\nu_1, \nu_2 \geq 0$; and $(1, 1) = f + \eta_2 r^2 + \eta_3r^3$, $\eta_2,\eta_3 \geq 0$.  Furthermore we may assume that $0 < f_1 < 1$ since otherwise $0 <f_2 <1$ and it is possible to apply a unimodular transformation such that $0 < f_1 < 1$. See Figure \ref{diq} for an illustration.

Before explaining the ideas of the proof, we need the following notations.

\begin{construction}[$\diamond^1$, $\diamond^3$]
Let $p^{\lambda} = f + \lambda r^1$. Let $q^{\lambda}_2$ be the intersection point of the line $p^{\lambda}(0,1)$ with the ray $\{x \in \mathbb{R}^2\,|\, x = \lambda_2 r^2, \lambda_2 \geq 0\}$. Similarly let $q^{\lambda}_4$ be the intersection point of the line $p^{\lambda}(0,0)$ with the ray $\{x \in \mathbb{R}^2\,|\, x = \lambda_4 r^4, \lambda_4 \geq 0\}$. Let $\bar{\lambda} = \textup{sup}\{\lambda\in \mathbb{R}_{+}\cup \{+\infty\}\,|\, \textup{conv}\{p^{\lambda}, q^{\lambda}_2, q^{\lambda}_4\}\textup{ is lattice-free.}\}$. If $\bar{\lambda} = + \infty$, set $q_2 := \textup{lim}_{\lambda \rightarrow +\infty}q^{\lambda}_2$ and $q_4 := \textup{lim}_{\lambda \rightarrow +\infty}q^{\lambda}_4$. Else set $q_2 := q^{\bar{\lambda}}_2$ and $q_4 := q^{\bar{\lambda}}_4$. Let $q_3 = f + \lambda_3r^3$ be the point such that either the line segment $q_2q_3$ or the line segment $q_4q_3$ contains an integer point and the triangle $q_2q_3q_4$ is lattice-free. If $\bar{\lambda} = + \infty$, define $\diamond^1 := \textup{conv}\{q_2, q_3, q_4\} + \textup{cone}(r^1)$. Otherwise define $\diamond^1 := \textup{conv}\{p^{\bar{\lambda}},q_2, q_3, q_4\}$.

Construct $\diamond^3$ symmetrically by exchanging the role of $r^1$ and $r^3$.
\label{wedgequad}
\end{construction}
\begin{proposition}\label{q2prop}
The split rank of an inequality whose induced lattice-free set is a quadrilateral of type $Q^2$ is finite.
\end{proposition}
\textit{Proof:}
We define
\begin{itemize}
\item For $i \in \{1, ...4\}$, let $w^i = f + \sigma_i r^i$, $\sigma_i \geq 0$ be the vertices of $L_{\alpha}$.
\item For $i\in \{0,1\}$, let $y^i$ be the intersection point of the line segment $w^1w^3$ and the line $\{x\in \R^2\mid x_1=i\},$
\item let $Q^a:= \diamond^1$, $\beta^a:=\phi(Q^a)$ and let $Q^b:= \diamond^3$, $\beta^b:=\phi(Q^b)$,
\item We denote by $q^{a,1},q^{a,2},q^{a,3},q^{a,4}$ the vertices of $Q^a$ and by $q^{b,1},q^{b,2},q^{b,3},q^{b,4}$ the vertices of $Q^b$; where $q^{a,i}$ (resp. $q^{b,i}$) lies on the ray $\{x \in \mathbb{R}^2\,|\,f+\lambda_i r^i, \lambda_i \geq 0\}$ (If $r^1_1 = 0$, and $Q^a$ is a subset of the split set, $\{x\in \mathbb{R}^2\,|\, 0 \leq x_2 \leq 1\}$ then $q^{a,1}$ is not defined. Similarly for $q^{a,3}$, $q^{b,1}$, and $q^{b,3}$.)
\end{itemize}

Observe that by construction, $Q^a$ and $Q^b$ are quadrilaterals of type $Q^1$ or triangles of types $T^{2B}$, $T^{2A}$, or subsets of a split set. This is because if $Q^a$ is not a subset of some split set, then by construction either $q^{a,2}$ or $q^{a,4}$ is a vertex of $Q^a$ and both these points lie in the interior of $L_{\alpha}$. Hence they are not integer and $Q^{a}$ cannot be a triangle of type $T^1$. Using a similar argument for $Q^b$, we obtain that the split rank of $\beta^a$ and $\beta^b$ is finite by Proposition \ref{q1prop}.

In the rest of the proof we assume that $Q^{a}$ and $Q^{b}$ are not subsets of a split set (i.e,  $q^{a,1}$, $q^{a,3}$, $q^{b,1}$, and $q^{b,3}$ are well-defined). This is for simplicity and the proof can be modified for the cases where $Q^{a}$ and $Q^{b}$ are subsets of split sets. Also for the purpose of our proof, it is convenient to consider a weaker version of $\beta^a$ for the coefficient of the ray $r^3$ and a weaker version of $\beta^b$ for the coefficient of $r^1$. We define $\bar Q^a:=\conv\{q^{a,1},q^{a,2},\bar q^{a,3},q^{a,4}\}$, $\gamma^a:=\phi(\bar Q^a)$ where $\bar q^{a,3}$ is obtained as the intersection point of the line $q^{a,1}y^0$ with $\{f+\lambda_3 r^3, \lambda_3\geq 0\}$. Similarly $\bar Q^b:=\conv\{\bar q^{b,1},q^{b,2},q^{b,3},q^{b,4}\}$, $\gamma^b:=\phi(Q^b)$ and $\bar q^{b,1}$ is the intersection point of the line $q^{b,3}y^1$ with $\{f+\lambda_1 r^1, \lambda_1\geq 0\}$.

Claim: $\gamma^a$ is weaker than $\beta^a$. Since $\gamma^a_i = \beta^a_i$ for $i \in \{1,2,4\}$, we only need to prove that $\gamma^a_3 \geq \beta^a_3$ or equivalently $\lambda(\bar{q}^{a,3}) \leq \lambda(q^{a,3})$. Observe that by construction $\lambda(\bar{q}^{a,2}) \leq \lambda(w^{2})$ and $\lambda(\bar{q}^{a,4}) \leq \lambda(w^{4})$. This implies $\lambda(q^{a,3}) \geq \lambda(w^3)$. On the other hand $\lambda(q^{a,1}) \geq \lambda(w^1)$ by construction. Since the lines $q^{a,1}\bar{q}^{a,3}$ and $w^1w^3$ meet at $y^{0}$, we conclude that $\lambda(\bar{q}^{a,3}) \leq \lambda(w^3)$ which proves $\lambda(\bar{q}^{a,3}) \leq \lambda(q^{a,3})$.

Symmetrically $\gamma^b$ is weaker than $\beta^b$. Therefore the split rank of $\gamma^a$ and $\gamma^b$ is finite.

We now consider the following set
\begin{align*}
Q:=\{(x,s)\in \R^2\times \R^4_+ \mid & x = f + Rs\\
& \sum_{i=1}^4 \gamma^a_i s_i \geq 1\\
& \sum_{i=1}^4 \gamma^b_i s_i \geq 1\quad \}.
\end{align*}
We claim that $\alpha$ is valid for $\textup{conv}(Q^{\leq} \cup Q^{\geq})$ where $Q^{\leq} := Q \cap \{(x,s) \in \mathbb{R}^2 \times \mathbb{R}^3_{+}\,|\, x_1 \leq 0\}$ and $Q^{\geq} := Q \cap \{(x,s) \in \mathbb{R}^2 \times \mathbb{R}^3_{+}\,|\, x_1 \geq 1\}$. To do this, we directly check that every vertex of $Q^{\leq}$ and $Q^{\geq}$ is valid for $\alpha$. The proofs for $Q^{\leq}$ and $Q^{\geq}$ being completely  symmetric (we exchange the role of $r^1$ and $r^3$ and $\bar Q^{a}$ and $\bar Q^b$), we only prove it for the vertices of $Q^{\leq}.$

We now discuss all the vertices $(\bar{x}, \bar{s})$ of $Q^{\leq}$.
\begin{enumerate}
\item Vertices $(\bar{x}, \bar{s})$ of $Q^{\leq}$, where the support of $\bar{s}$ is 1 and $(\bar{x}, \bar{s})$ is tight for $\gamma^a$.
\begin{itemize}
\item $(q^{a,1}, \mathcal{M}^{1,1}(q^{a,1}))$ is a vertex of $Q^{\leq}$. It is valid for $\alpha$ since $\lambda(q^{a,1}) \geq \lambda(w^1)$ (see Proposition \ref{therelight}).
\item $(q^{a,2}, \mathcal{M}^{2,2}(q^{a,2}))$ is not a vertex of $Q^{\leq}$. Since $\lambda(q^{a,2}) \leq \lambda(w^2)$, $q^{a,2}$ is a convex combination of $f$ with $w^2$. Since $f_1, w^2_1 >0$ we obtain that $q^{a,2}_1>0$. Therefore, $(q^{a,2}, \mathcal{M}^{2,2}(q^{a,2}))$ does not satisfy $x_1 \leq 0$.
\item $(q^{a,3}, \mathcal{M}^{3,3}(q^{a,3}))$ is not a vertex of $Q^{\leq}$. Since it is not valid for $\gamma^{b}$ (since $\lambda(q^{a,1}) \leq \lambda(w^2) < \lambda(\bar{q}^{b,1})$).
\item $(q^{a,4}, \mathcal{M}^{4,4}(q^{a,4}))$ is not a vertex of $Q^{\leq}$. Since $\lambda(q^{a,4}) \leq \lambda(w^4)$, $q^{a,4}$ is a convex combination of $f$ with $w^4$. Since $f_1, w^4_1 >0$ we obtain that $q^{a,4}_1>0$. Therefore, $(q^{a,4}, \mathcal{M}^{4,4}(q^{a,4}))$ does not satisfy $x_1 \leq 0$.
\end{itemize}
\item Vertices $(\bar{x}, \bar{s})$ of $Q^{\leq}$, where the support of $\bar{s}$ is 1 and $(\bar{x}, \bar{s})$ is tight for $\gamma^b$.
\begin{itemize}
\item $(\bar{q}^{b,1}, \mathcal{M}^{1,1}(\bar{q}^{b,1}))$ is not a vertex of $Q^{\leq}$, since $\lambda(\bar{q}^{b,1}) \leq \lambda(w^1) < \lambda(q^{a,1})$ and therefore it is not valid
for $\gamma^a$.
\item $(q^{b,2}, \mathcal{M}^{2,2}(q^{b,2}))$ is not a vertex of $Q^{\leq}$. Since $\lambda(q^{b,2}) \leq \lambda(w^2)$, $q^{b,2}$ is a convex combination of $f$ with $w^2$. Since $f_1, w^2_1 >0$ we obtain that $q^{b,2}_1>0$. Therefore, $(q^{b,2}, \mathcal{M}^{2,2}(q^{b,2}))$ does not satisfy $x_1 \leq 0$.
\item $(q^{b,3}, \mathcal{M}^{3,3}(q^{b,3}))$ is not a vertex of $Q^{\leq}$, since $f_1 >0$ and $r^3_1 >0$.
\item $(q^{b,4}, \mathcal{M}^{4,4}(q^{b,4}))$ is not a vertex of $Q^{\leq}$. Since $\lambda(q^{b,4}) \leq \lambda(w^4)$, $q^{b,4}$ is a convex combination of $f$ with $w^4$. Since $f_1, w^4_1 >0$ we obtain that $q^{b,4}_1>0$.
\end{itemize}

\item Vertices $(\bar{x}, \bar{s})$ of $Q^{\leq}$, where the support of $\bar{s}$ is 1 and $(\bar{x}, \bar{s})$ is tight for $x_1 \leq 0$. Let $v^{i,0}$ be the intersection point of $\{x \in \mathbb{R}^2\,|\, x = f + \lambda_ir^i, \lambda_i \geq 0\}$ with the line $\{x \in \mathbb{R}^2\,|\, x_1 = 0\}$.
\begin{itemize}
\item The point $(v^{1,0}, \mathcal{M}^{1,1}(v^{1,0}))$ is not a vertex of $L_{\alpha}$ since $w^1_1 <0$ and therefore, $\lambda(q^{a,1}) \geq \lambda(w^1) > \lambda(v^{1,0})$. (Thus $(v^{1,0}, \mathcal{M}^{1,1}(v^{1,0}))$ is not valid for $\gamma^a$.
\item The point $(v^{2,0}, \mathcal{M}^{2,2}(v^{2,0}))$ is a vertex for $Q^{\leq}$ (if it exists) since $w^2_1 > 0$ and therefore $\lambda(v^{2,0}) > \lambda(w^2)$. Note also that the previous statement shows that $(v^{2,0}, \mathcal{M}^{2,2}(v^{2,0}))$ is valid for $\alpha$.
\item Similarly we can verify that $(v^{3,0}, \mathcal{M}^{3,3}(v^{1,0}))$ is not a vertex. Also if $(v^{4,0}, \mathcal{M}^{4,4}(v^{4,0}))$ exists, then it is a vertex for $Q^{\leq}$ and also valid for $\alpha$.
\end{itemize}
\item Vertices $(\bar{x}, \bar{s})$ of $Q^{\leq}$, where the support of $\bar{s}$ is 2 and $(\bar{x}, \bar{s})$ is tight for $\gamma^a$ and $\{x \in \mathbb{R}^2\,|\,x_1 = 0\}$.
\begin{itemize}
\item $((0,1), \mathcal{M}^{1,2}((0,1)))$ is a vertex of $Q^{\leq}$ and is valid for $\alpha$.
\item $(y^{0}, \mathcal{M}^{1,3}(y^{0}))$ is a vertex of $Q^{\leq}$ and is valid for $\alpha$. We remark here that by construction, $y^{0}$ is the intersection point of $\textup{conv}(q^{a,1}, \bar{q}^{a,3})$ with the line $\{x \in \mathbb{R}^2\,|\,x_1 = 0\}$. Therefore $(y^{0}, \mathcal{M}^{1,3}(y^{0}))$ is tight for $\gamma^a$.
\item $((0,0), \mathcal{M}^{1,4}((0,0)))$ is a vertex of $Q^{\leq}$ and is valid for $\alpha$.

\item All points of the form $(\bar{x}, \mathcal{M}^{2,3}(\bar{x}))$, $(\bar{x}, \mathcal{M}^{2,4}(\bar{x}))$, and $(\bar{x}, \mathcal{M}^{3,4}(\bar{x}))$ that are tight for $\gamma^a$ satisfy $x_1 > 0$ and are therefore not valid for $Q^{\leq}$. Any such point must be convex combination of two of the three points: $(q^{a,2}, \mathcal{M}^{2,2}(q^{a,2}))$, $(q^{a,3}, \mathcal{M}^{3,3}(q^{a,3}))$, $(q^{a,4}, \mathcal{M}^{4,4}(q^{a,4}))$. Since $\lambda(q^{a,2}) \leq \lambda(w^2)$, $\lambda(q^{a,4}) \leq \lambda(w^4)$, $q^{a,2}$ and $q^{a,4}$ are convex combination of $f$ with $w^2$ and $w^4$ respectively. Since $f_1 , w^2_1, w^4_1 >0$ we obtain that $q^{a,2}_1$, $q^{a,4}_1 >0$. Moreover since $w^3_1 >0$, we obtain $q^{a,3} >0$.
\end{itemize}
\item Vertices $(\bar{x}, \bar{s})$ of $Q^{\leq}$ where the support of $\bar{s}$ is 2 and $(\bar{x}, \bar{s})$ is tight for $\gamma^b$ and $\{x \in \mathbb{R}^2\,|\,x_1 = 0\}$. First consider a point of the form $(\bar{x}, \mathcal{M}^{1,2}(\bar{x}))$. There are two cases. Either $\bar x_2 < 1$ and therefore, $\bar x \in \textup{int}(\textup{conv}\{f, q^{a,1}, q^{a,2}\})$ and therefore $(\bar{x}, \mathcal{M}^{1,2}(\bar{x}))$ is not valid for $\gamma^a$. On the other hand if $\bar x_2 \geq 1$, then $\bar x \notin \textup{int}(\textup{conv}\{f, w^1, w^2\})$ and $(\bar{x}, \mathcal{M}^{1,2}(\bar{x}))$ is valid for $\alpha$. A similar argument holds for points of the form $(\bar{x}, \mathcal{M}^{1,3}(\bar{x}))$ and $\tilde (x, \mathcal{M}^{1,4}(\tilde x))$ that are tight for $\gamma^b$ and
belong to $\{x \in \mathbb{R}^2\,|\,x_1 = 0\}$. Finally, all points of the form $(\bar{x}, \mathcal{M}^{2,3}(\bar{x}))$, $(\bar{x}, \mathcal{M}^{2,4}(\bar{x}))$, and $(\bar{x}, \mathcal{M}^{3,4}(\bar{x}))$ that are tight for $\gamma^b$ satisfy $\bar x_1 > 0$ (proof similar to the previous case) and are therefore not valid for $Q^{\leq}$. The proof is the same as that for the previous case.
\item Vertices $(\bar{x}, \bar{s})$ of $Q^{\leq}$ where the support of $\bar{s}$ is 2 and $(\bar{x}, \bar{s})$ is tight for $\gamma^a$ and $\gamma^b$.
\begin{itemize}
\item If $(\bar{x}, \mathcal{M}^{12}(\bar{x}))$ is tight for $\gamma^a$ and $\gamma^b$, then $(\bar{x}, \mathcal{M}^{12}(\bar{x}))$ is not a vertex of $Q^{\leq}$. Since $\lambda(q^{a,1}) \geq \lambda(\bar{q}^{b,1})$ and $(0,1) \not\in \textup{conv}\{f,\bar{q}^{b,1}q^{b,2}\}$ (as $\lambda(\bar{q}^{b,1}) < \lambda(w^1), \lambda(q^{b,2}) < \lambda(w^1)$), and $(0,1)$ is in the relative interior of the line $q^{a,1},q^{a,2}$, we conclude that any potential intersection point $\bar{x}$ of $\bar{q}^{b,1}q^{b,2}$ and $q^{a,1}q^{a,2}$ satisfies $\bar{x}_1 > 0$ and is therefore not valid for $Q^{\leq}$.
\item If $(\bar{x}, \mathcal{M}^{14}(\bar{x}))$ is tight for $\gamma^a$ and $\gamma^b$, then $(\bar{x}, \mathcal{M}^{14}(\bar{x}))$ is not a vertex of $Q^{\leq}$. A similar argument as above shows that $x_1 > 0$ and is therefore not valid for $Q^{\leq}$.
\item If $(\bar{x}, \mathcal{M}^{13}(\bar{x}))$ is tight for $\gamma^a$ and $\gamma^b$, then $(\bar{x}, \mathcal{M}^{13}(\bar{x}))$ is not a vertex of $Q^{\leq}$. Observe that $\{x \in \mathbb{R}^2\,|\, (x, \mathcal{M}^{1,3}(x)) \textup{ is tight for }\gamma^a\} = \textup{conv}\{q^{a,1}, \bar{q}^{a,3}\}$. Also $\lambda(q^{a,1}) \geq \lambda(w^1)$, $\lambda(\bar{q}^{a,3}) < \lambda(w^3)$, $y^{0} \in \textup{conv}\{q^{a,1}, \bar{q}^{a,3}\}$, and $\sum_{i = 1}^3\alpha_i(\mathcal{M}^{1,3}(y^0))_i = 1$. We therefore obtain the following inclusion, $\{x \in \mathbb{R}^2\,|\, (x, \mathcal{M}^{1,3}) \textup{ is tight for }\gamma^a \textup{ and valid for }\alpha\} = \textup{conv}\{q^{a,1}, y^{0}\} \subseteq \{x\in \mathbb{R}^2\,|\, x_1 \leq 0\}$. Similarly $\{x \in \mathbb{R}^2\,|\, (x, \mathcal{M}^{1,3}) \textup{ is tight for }\gamma^b \textup{ and valid for }\alpha\} = \textup{conv}\{y^1, q^{b,3}\} \subseteq \{x\in \mathbb{R}^2\,|\, x_1 \geq 1\}$. From these observations, we conclude that a point $\bar{x}$ such that $(\bar{x}, \mathcal{M}^{1,3}(\bar{x}))$ is tight for $\gamma^a$ and $\gamma^b$ cannot be valid for $\alpha$ since it would otherwise belong to $\textup{conv}(q^{a,1}, y^{0}) \cap \textup{conv}(y^1, q^{b,3}) = \emptyset$. Therefore, $(\bar{x}, \mathcal{M}^{1,3}(\bar{x}))$ is not valid for $\alpha$ and satisfies $0 < \bar{x}_1 < 1$. However, then $(\bar{x}, \mathcal{M}^{1,3}(\bar{x}))$ is not valid for $Q^{\leq}$.
\item Furthermore all points of the form $(\bar{x}, \mathcal{M}^{2,3}(\bar{x}))$, $(\bar{x}, \mathcal{M}^{2,4}(\bar{x}))$, and $(\bar{x}, \mathcal{M}^{3,4}(\bar{x}))$ that are tight for $\gamma^a$ and $\gamma^b$, satisfy $\bar{x}_1 > 0$ (proof similar to the previous case) and are therefore not valid for $Q^{\leq}$.
\end{itemize}
\item Vertices $(\bar{x}, \bar{s})$ of $Q^{\leq}$ where the support of $\bar{s}$ is 3. Then $(\bar{x}, \bar{s})$ is tight for $\gamma^a$, $\gamma^b$ and $\bar{x}_1 = 0$. Since $(\bar{x}, \bar{s})$ is tight for $\gamma^a$ and $\gamma^b$, it must be a convex combination of points of the form $(\hat{x}, \mathcal{M}^{ij}(\hat{x}))$, $i, j \in \{1,2,3,4\}$, where $(\hat{x}, \mathcal{M}^{ij}(\hat{x}))$ are tight for $\gamma^a$ and $\gamma^b$. However from the previous case, such an $\bar{x}$ satisfies $\bar{x}_1 > 0$. Therefore, such a vertex does not exist.
\hfill $\square$
\end{enumerate}

\bibliographystyle{amsplain}
\bibliography{facetref}
\section*{Appendix 1}
\begin{center}
\begin{figure}[!h]
\begin{center}
\includegraphics[width=\textwidth]{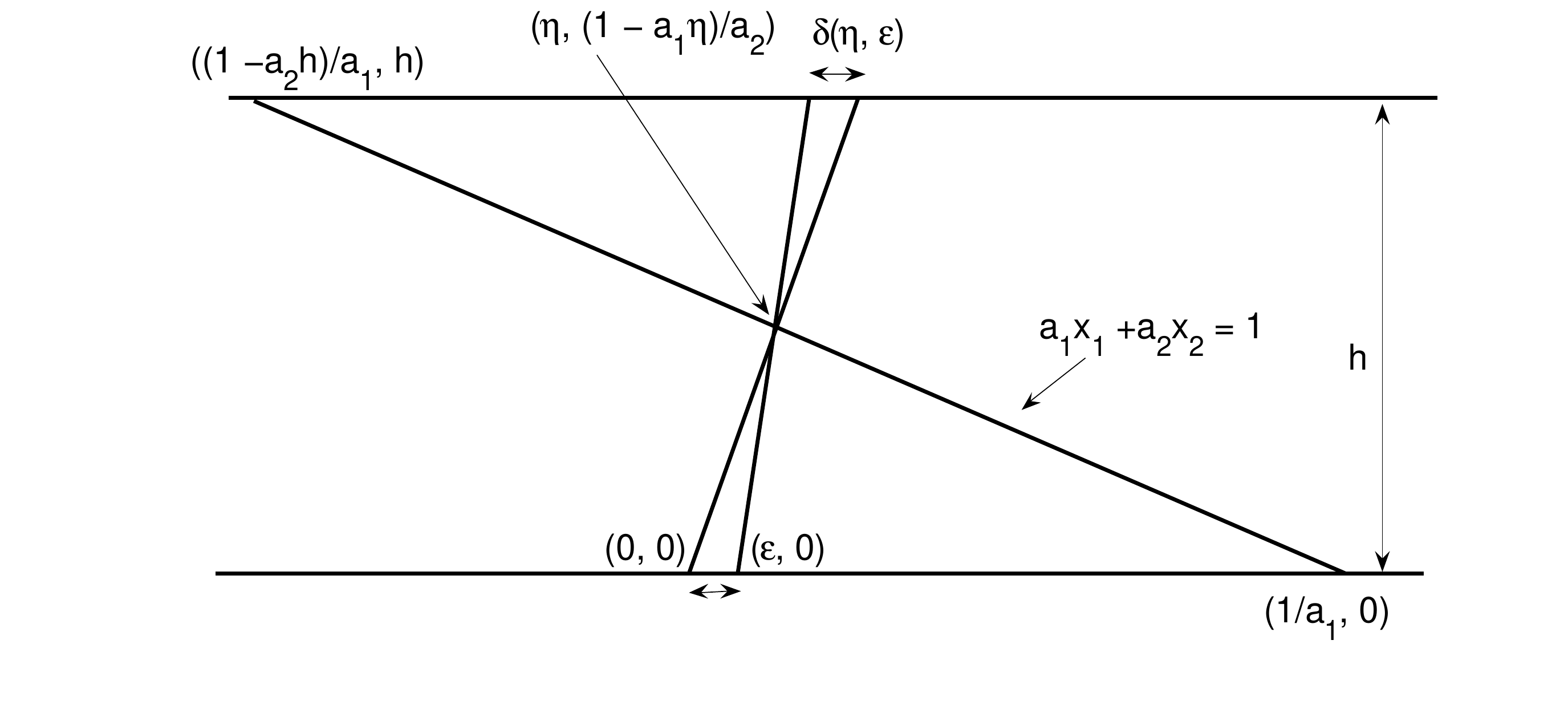}
\end{center}
\caption{For Proof of Observation \ref{phase2nalem}}
\label{phase2na}
\end{figure}
\end{center}

\begin{observation}\label{phase2nalem}
To simplify the computation of $|\!|u^{[3]}u^{[2]}|\!|$, we suitably rotate and translate the points in Figure \ref{phase2c}. In particular, we rotate so that the line $\{x \in \mathbb{R}^2\,|\,x =z+ \lambda r^2, \lambda \in \mathbb{R}\}$ is the line $\{x\in \mathbb{R}^2\,|\,x_2 = 0\}$ and the $\{x \in \mathbb{R}^2\,|\,x= f+ \lambda r^2, \lambda \in \mathbb{R}\}$ is the line $\{x\in \mathbb{R}^2\,|\,x_2 = h\}$ where $h > 0$. We now translate so that the point $z$ becomes $(0,0)$. Under this rotation and translation, the line $\{x \in \mathbb{R}\,|\,x_1 =0\}$ becomes the line $\{x\in \mathbb{R}^2\,|\,a_1x_1 + a_2x_2 = 1\}$ for some $a_1 > 0$ and $a_2 > 0$. Refer to Figure \ref{phase2na}.

The point $(\eta, \frac{1 - a_1\eta}{a_2})$ (representing $p^{[3]}$) lies on the line $\{x\in \mathbb{R}^2\,|\,a_1x_1 + a_2x_2 = 1\}$, where $(1 - a_1\eta) > 0$. The line passing through $(0,0)$ and $(\eta, \frac{1 - a_1\eta}{a_2})$ meets the line $\{x\in \mathbb{R}^2\,|\,x_2 = h\}$ at $A:= (\frac{ha_2\eta}{(1 - a_1\eta)}, h)$ (representing $u^{[3]}$). Another line passing through  $(\epsilon,0)$ (this represents $z^2$, i.e, $\epsilon = |\!|z^2 - z|\!|$) and $(\eta, \frac{1 - a_1\eta}{a_2})$ meets the line $\{x\in \mathbb{R}^2\,|\,x_2 = h\}$ at $B:= (\frac{(1 - a_1\eta)\epsilon - (\epsilon - \eta)a_2h}{(1 - a_1\eta)}, h)$ (representing $u^{[2]}$). Then $|\!|u^{[3]}u^{[2]}|\!| = \delta(\eta, \epsilon) =  |A_1 - B_1| = \frac{ha_2\eta}{(1 - a_1\eta)}- \frac{(1 - a_1\eta)\epsilon - (\epsilon - \eta)a_2h}{(1 - a_1\eta)} = \epsilon\left(\frac{a_2h}{1 - a_1\eta} -1\right).$ \hfill$\square$
\end{observation}

\end{document}